\documentclass[11.5pt]{amsart}
\usepackage{amsfonts}
\usepackage{amssymb}
\usepackage[active]{srcltx}

\usepackage[all,poly,graph]{xy}

\font\got=eufm10
\font\gots=eufm10 at 7pt

\def\g2{ \hbox{\got g}_2}
\def\f4{\hbox{\got f}_4}
\def\d4{\hbox{\got d}_4}
\def\e6c{\hbox{\got e}_6^\mathbb{C}}
\def\fs4{\hbox{\gots f}_4}
\def\F4{\hbox{\got F}_4}
\def\Fs4{\hbox{\gots F}_4}
\def\eseis{\hbox{\got e}_6}
\def\es7{ \hbox{\got e}_7}
\def\es8{ \hbox{\got e}_8}
\def\id{\mathop{\hbox{id}}}

\def\L{{\mathcal L}}

\def\imag{\mathfrak{i}}

\def\T{\mathcal{T}}

\def\Q{ {\mathcal Q}}

\def\W{\mathcal W}
\def\V{\mathcal V}

\def\X{\hbox{\got X}}
\def\Xs{\hbox{\gots X}}

\def\Z{\mathbb Z}
\def\F{\mathbb F}
\def\R{\mathbb R}
\def\aut{\mathop{\rm Aut}}
\def\der{\mathop{\rm Der}}
\def\End{\mathop{\rm End}}

\def\Ad{\mathop{\rm Ad}}
\def\sop{\mathop{\rm Supp}}

\def\Int{\mathop{\rm Int}}
\def\sll{\mathop{\rm sl}}
\def\SL{\mathop{\rm SL}}
\def\SP{\mathop{\rm SP}}
\def\PSL{\mathop{\rm PSL}}
\def\GL{\mathop{\rm GL}}

\def\a{\alpha}
\def\si{\sigma}

\def\span#1{\langle #1 \rangle}

\def\exte{\mathop{\rm ext}}

\def\fix{\mathop{\rm fix}}
\def\tr{\mathop{\rm tr}}
\def\cent{\mathop{\rm cent}}
\def\sop{\mathop{\hbox{Supp}}}
\def\diag{\mathop{\rm diag}}
\def\To{\hbox{\got T}}

\def\No{\hbox{\got N}}

\def\tor#1{\T^{\langle #1\rangle}}
\def\tors#1{\mathfrak{T}^{\langle #1\rangle}}
\def\sor#1{\mathcal{S}^{\langle #1\rangle}}
 \def\hor#1{\mathcal{H}^{\langle #1\rangle}}
 \def\jor#1{\mathfrak{S}^{\langle #1\rangle}}

\newtheorem{te}{Theorem}
\newtheorem{pr}{Proposition}
\newtheorem{lm}{Lemma}

\newtheorem{re}{Remark}
\title[Fine Gradings on     $\frak e_6$]{ Fine Gradings on $\frak e_6$}
\author{Cristina Draper}
\thanks{\ Partially supported by MCYT grant MTM2010-15223  and by   the Junta de Andaluc\'{\i}a PAI
projects  FQM-336,  FQM-2467, FQM-3737.}
\address{Cristina Draper Fontanals: Departamento de
Matem\'atica Aplicada\\ Escuela de las Ingenier\'{\i}as\\ Ampliaci\'{o}n Campus de Teatinos, S/N, 29071 M\'alaga,
Spain\\cdf@uma.es}
\author{Antonio Viruel}
\thanks{\ Partially supported by the MCYT
grant MTM2010-18089 and by the Junta de Andaluc\'{\i}a PAI
projects FQM-2863 and  FQM-213}
\address{Antonio Viruel Arb\'{a}izar:
Departamento de \'Algebra, Geometr\'{\i}a y Topolog\'{\i}a\\
Campus de Teatinos, S/N. Facultad de Ciencias, Ap. 59, 29080
M\'alaga, Spain\\viruel@agt.cie.uma.es}

\usepackage{graphicx}
\begin{document}

\setlength{\unitlength}{0.06in}

\maketitle

\begin{abstract}
There are fourteen fine gradings on the exceptional Lie algebra $\frak e_6$ over an algebraically closed field of zero characteristic.
 We provide their descriptions and a proof that they are all.
\end{abstract}

\section{Introduction}

The concept of gradings on algebraic structures has a long tradition in Mathematics, appearing in the classics as \cite[Chapter III, \S 3]{bourbaki1}. Gradings are the natural ambient for algebraic structures arising from geometrical situations, so all the algebraic structures considered in Algebraic Topology are graded, and gradings of Lie algebras have its source in the original work of Jordan \cite{jordan} towards a mathematical description of Quantum Mechanics. More recent results describe particle physics in terms of gradings of Lie algebras, and these are in the background in many almost any mathematical construction involving an ordered choice of basis.

Motivated by these manifold applications of gradings on Lie algebras, Patera et al.\ initiated in \cite{LGI} a systematic classification of all the possible group gradings on finite-dimensional Lie algebras. For algebraically closed fields of characteristic zero, the gradings on the classical Lie algebras have been studied in \cite{Ivan2,tipoA, d4,Albclasicas} and the gradings in some exceptional ones, namely, $\mathfrak{f}_4$ and $\mathfrak{g}_2$, in \cite{g2,f4,f4spin,otrog2}. Lately, some authors have already studied the case of prime characteristic, \cite{BahKotch} in the classical case (with the exception of $\mathfrak{d}_4$) and \cite{AlbKoche} in $\mathfrak{f}_4$ and $\mathfrak{g}_2$.
Over arbitrary fields, almost nothing is known: for example, over the real field the gradings on simple Lie algebras have been studied only  for the exceptional  algebras $\mathfrak{f}_4$ and $\mathfrak{g}_2$ (\cite{reales}) and in classical ones for some examples in low-dimensions (\cite{LGIII}).

 Summarizing, the classification of fine group gradings on finite-dimensional semi\-simple   Lie algebras is almost complete, remaining just the E-family to be fully studied, although some examples of $\eseis$ have been described (for instance, in \cite{AlbJordangrad}). This work intends to be a first step towards to fill that gap by giving a complete description of the fine group gradings of the smallest representative of Lie algebras of the E-series:

 \begin{te}
There are $14$ fine group gradings on $\eseis$ up to equivalence. The following table describe them in terms of the associated MAD-group $Q$ in $\aut\eseis$, the type of the fine grading and the dimension of the fixed subalgebra $\fix Q=L_e$.
\begin{center}{
\begin{tabular}{|c|c|c|c|}
\hline
 Quasitorus & Isomorphic to & Type &$\dim L_e$\cr
\hline\hline
 $\Q_{1}$ & $ \Z_3^4$   &  $( 72,0,2 )$  &$0$\cr
\hline
 $\Q_{2}$ & $(\F^*)^2\times\Z_3^2$   &  $( 60,9 )$  &$2$\cr
\hline
$\Q_{3}$ & $\Z_3^2\times\Z_2^3$   &  $( 64,7 )$  &$0$\cr
\hline
$\Q_{4}$ & $(\F^*)^2\times\Z_2^3$   &  $( 48,1,0,7 )$ &$2$ \cr
\hline
$\Q_{5}$ & $(\F^*)^6 $   &  $(72,0,0,0,0,1  )$  &$6$\cr
\hline
$\Q_{6}$ & $(\F^*)^4\times\Z_2$   &  $(72,1,0,1  )$  &$4$\cr
\hline
$\Q_{7}$ & $ \Z_2^6$   &  $( 48,1,0,7 )$  &$0$\cr
\hline
$\Q_{8}$ & $\F^*\times\Z_2^4$   &  $(  57,0,7)$  &$1$\cr
\hline
$\Q_{9}$ & $\Z_3^3\times\Z_2$   &  $( 26,26 )$  &$0$\cr
\hline
$\Q_{10}$ & $(\F^*)^2\times\Z_2^3$   &  $(60,7,0,1   )$  &$2$\cr
\hline
$\Q_{11}$ & $\Z_4\times\Z_2^4$   &  $(48,13,0,1     )$ &$0$ \cr
\hline
$\Q_{12}$ & $\F^* \times\Z_2^5$   &  $( 73,0,0,0,1  )$ &$1$ \cr
\hline
$\Q_{13}$ & $ \Z_2^7$   &  $( 72,0,0,0,0,1  )$ &$0$ \cr
\hline
$\Q_{14}$ & $ \Z_4^3$   &  $( 48,15 )$  &$0$\cr
\hline
 \end{tabular}}\end{center}\end{te}\smallskip

The first five MAD-groups contain no outer automorphisms. Note that the type of the fine grading  (how many pieces are of each size) jointly with the dimension of the fixed subalgebra are enough to distinguish the conjugacy
class of the quasitorus inducing such grading.

Due to the dimension of $\eseis$, the techniques previously used for the other Lie algebras seem to be inefficient in this case. The classical Lie algebras were studied by taking into consideration the gradings on the associative matrix algebras in which they live. The gradings on $\mathfrak{g}_2$ are induced from the gradings on the octonions, described in \cite{Alb1}. The gradings on $\mathfrak{f}_4$ have been submitted to several attacks. The first of them  uses this tool: the (fine) grading is produced by a maximal diagonalizable subgroup of the automorphism group of the algebra, and this subgroup lives in the normalizer of a maximal torus. Precisely, if the group is a maximal torus, we obtain the root decomposition. In  \cite{f4}, all the elements in this normalizer are constructed (as matrices of size $52$) and then a detailed analysis of the possible cases is realized.  It turns out very difficult to apply this technique in this case, which would require of more than one hundred thousand   square matrices of size $78$. What has been done is to work, not in the normalizer, but in its projection quotient the torus, isomorphic to the group of automorphisms of the root system (the extended Weyl group). This allows to work with   matrices of size $6$, which can be implemented in any computer. Even more, we use the computer only in two proofs  during the  paper. The remaining computations, although long, are made by hand, and the representatives or the orbits of the Weyl group are extracted from \cite{Atlas}.

  Another technique used  in this work is to take advantage of the knowledge of the elementary $p$-groups and their centralizers in the complex case. Thus, there is an   effort to use these results for algebraically closed fields of characteristic zero. This saves some computational work. Some basic facts about the structure of the MAD-groups are studied in order to make use of such elementary groups, although most of the useful information is extracted and generalized from \cite{f4,d4,f4spin}.

About the results we would like to mention the great amount of fine gradings on $\eseis$. If we take into account that in $\mathfrak{g}_2$ there are $2$ fine gradings, and in $\mathfrak{f}_4$ there are $4$   ones, we find that the number in this case is much bigger, proportionally speaking. Of course the reason is a greater symmetry in $\eseis$ than in the two previously mentioned cases: the presence of outer automorphisms is an signal of such   symmetry. Besides we find some remarkable gradings:
It was thought that every outer fine grading on a simple Lie algebra of type $A$ was induced by a quasitorus of automorphisms containing
an outer order two automorphism. Elduque proved that this is not true in
  \cite{Albclasicas}, in which he gave a revised version on the fine gradings on the simple classical Lie algebras with new arguments (and some new fine grading).
  The same fact happens to the gradings on $\eseis$. There is an outer fine grading such that the MAD-group producing the grading does not contain
  outer involutions. This is an unexpected $\mathbb{Z}_4^3$-grading and it turns out  a nice  symmetry based on the number $4$.
Over other finite groups, a $\mathbb{Z}_2^7$-grading and a $\mathbb{Z}_3^4$-grading appear.

The structure of the work is as follows. In Section~\ref{sec_generalidades}, we recall some generalities about gradings and we study the quasitori of the automorphism group producing the gradings. We deep into the structure of a MAD-group, that is, the maximal abelian diagonalizable group which is producing a fine grading. Besides, in this section about generalities, we recall the elementary $p$-groups of $\Int\e6c$ and try to translate the result to a more general field. All this will be used in
Section~\ref{sec descripcionesinternos}, after describing some examples of fine gradings on $\eseis$, to prove that they are all the possible cases of fine gradings produced by inner automorphisms. This only needs a technical result, the fact that every MAD-group contains a minimal elementary non-toral $p$-group (even more, of minimum rank). This technical result is proved in
Section~\ref{sec_demostraciones}, in which the Weyl group is recalled, their orbits under conjugation exhibited, and the result is got by means of a technical play in which   the representatives of the orbits in this Weyl group and the elements in a maximal torus fixed by them have an important role.
Afterwards, in
Section~\ref{sec descripcionesexternos} we describe all the possible fine gradings, up to equivalence, produced by a group of automorphisms not all of them inner. We consider the gradings obtained extending gradings on $\mathfrak{f}_4$ and  the ones obtained extending gradings on $\mathfrak{c}_4$, and last, we show an example of (outer) grading which is not in any of the two previous situations. Finally
Section~\ref{sec_demostracionesouter} is devoted to prove that there are no more fine gradings that seen before. Again we make use of a technical result, proved with a computer. It also contains the descriptions of the MAD-groups in computational terms in Equation~(\ref{eq_losMADsversioncomputacional}).

At the end of the paper we have added   an Appendix, with natural representatives of the orbits used through the work.


\section{Generalities}\label{sec_generalidades}

\subsection{Notions about gradings}\label{subsec_nocionesgenerales}

  If $L$ is a finite-dimensional Lie  algebra and $G$ is an
abelian group, a decomposition $\Gamma:L=\oplus_{g\in
G}L_g$ is said to be a \emph{$G$-grading} whenever for all $g,h\in G$,
$L_gL_h\subset L_{g+h}$ and $G$ is generated by the set
$\sop(\Gamma):=\{g\in G\colon L_g\ne0\}$, called the \emph{support} of
the grading. The subspaces $L_g$ are referred to as the \emph{homogeneous components} of the grading.

Given two gradings $\Gamma:L=\oplus_{g\in G}L_g$ and $\Gamma':L=\oplus_{g'\in G'}L_{g'}$ over two abelian groups $G$ and $G'$,
$\Gamma$ is said to be a \emph{refinement} of $\Gamma'$ (or $\Gamma'$ a \emph{coarsening} of $\Gamma$) if for any $g\in G$ there is  $g'\in G'$ such that $L_g\subseteq L_{g'}$. That is, any homogeneous component of $\Gamma'$ is the direct sum of several homogeneous components of $\Gamma$. The refinement is \emph{proper} if there are $g\in G$ and $g'\in G'$ such that $L_g\subsetneq L_{g'}$. A grading is said to be \emph{fine} if it admits no proper refinements.

 The gradings $\Gamma$ and $\Gamma'$ are said to be \emph{equivalent} if the sets
of homogeneous subspaces are the same up to isomorphism, that is,
if there are an automorphism $f\in\aut L$ and a bijection between
the supports $\alpha\colon \sop(\Gamma)\to \sop(\Gamma')$ such that
$f(L_g)=L_{\alpha(g)}$ for any $g\in \sop(\Gamma)$.
The \emph{type} of a grading $\Gamma$   (following \cite{Hesse})  is the sequence of numbers $(h_1,\ldots,h_r)$ where $h_i$ is the number of homogeneous components of dimension $i$, $i=1,\ldots,r$, with $h_r\ne 0$. Thus $\dim L=\sum_{i=1}^r ih_i$. This sequence is of course an invariant up to equivalence.
Our
objective is to classify fine gradings up to equivalence, because any grading is obtained as a coarsening of some fine grading.

  Given a grading $\Gamma:{
{L}}=\oplus_{g\in G}{ {L}}_g$, consider the free abelian group
${\tilde G}$ generated by $\sop(\Gamma)$ and subject to the
relations $g_1+g_2=g_3$ if $0 \neq [{L}_{g_1},{L}_{g_2}] \subset
{L}_{g_3}.$ The group ${\tilde G}$ is called    the \emph{universal grading group} of $\Gamma$.
Note that we have    a ${\tilde G}$-grading
$\tilde\Gamma: { {L}}=\oplus_{{\tilde g}\in {\tilde G}}{
{L}}_{\tilde{g}}$ equivalent to $\Gamma$,
where ${ {L}}_{\tilde{g}}$ is the sum of the
homogeneous spaces $L_g$ of $\Gamma$ such that the class of $g$ in
${\tilde G}$ is ${\tilde g}.$
Besides ${\tilde G}$ has the
following universal property: given any coarsening $L=\oplus_{h\in
H}{ {L}}_h$ of $\tilde\Gamma$, there exists a unique group epimorphism
$\alpha\colon {\tilde G} \to H$ such that $L_h=\oplus_{{\tilde g} \in
\alpha^{-1}(h)}{ {L}}_{{\tilde g}}$.

The ground field $\F$ will be supposed to be algebraically closed
and of characteristic zero  throughout this work. In this context, the
group of automorphisms of the algebra $L$ is an algebraic linear
group.  There is a deep connection between gradings on $L$ and
quasitori of the group of automorphisms $\aut L$, according to
\cite[\S 3, p.\,104]{enci}.
If $L=\oplus_{g\in G}L_g$ is a $G$-grading, the map
$\psi\colon\X(G)=\hom(G,\F^\times)\to\aut L$ mapping each
$\alpha\in\X(G)$ to the automorphism $\psi_\alpha\colon L\to L$
given by $L_g\ni x\mapsto \psi_\alpha(x):=\alpha(g)x$ is a group
homomorphism. Since $G$ is finitely generated, then $\psi(\X(G))$
is an algebraic quasitorus. And conversely, if $Q$ is a quasitorus and
$\psi\colon Q\to\aut L$ is a homomorphism, $\psi(Q)$ is  formed
by semisimple automorphisms
 and we have a $\X(Q)$-grading
$L=\oplus_{g\in{\Xs}(Q)}L_g$ given by $ L_g=\{x\in L\colon
\psi(q)(x)=g(q)x\ \forall q\in Q\} $, with $\X(Q)$  a finitely
generated abelian group.

If $\Gamma:L=\oplus_{g\in G}L_g$ is a
$G$-grading, then  the set of automorphisms of $L$ such that every
$L_g$ is contained in some eigenspace is an abelian group formed
by semisimple automorphisms, called the \emph{diagonal group} of the grading, and   denoted by Diag$(\Gamma)$.
If $\psi\colon\X(G)\to\aut L$  is the related homomorphism, then the group Diag$(\Gamma)$ contains
 $\psi(\X(G))$, and both groups coincide when $G$ is the universal grading group of $\Gamma$.

The grading is fine if and only if Diag$(\Gamma)$ is a maximal abelian
subgroup of semisimple elements, usually called a \emph{ MAD-group
}(\lq\lq maximal abelian diagonalizable" group). It is convenient to
observe that the number of conjugacy classes of MAD-groups in
$\aut L$ agrees with the number of equivalence classes of fine
gradings on $L$, and that if $Q$ is a MAD-group, then $\X(Q)$ is the universal group of the induced fine grading.

A grading is \emph{toral} if it is a coarsening of the root decomposition of a semisimple Lie algebra $L$ relative to some Cartan
subalgebra. In other words, if the grading is produced by
a quasitorus contained in a torus of the automorphism group of $L$.
If $L_e$ denotes the identity component of a grading on $L$, such grading is toral if and only if $L_e$ contains a Cartan subalgebra
of $L$ (\cite[Subsection~2.4]{g2}). As $L_e$ in any case
 is a reductive subalgebra (\cite[Remark 3.5]{Hesse}), the grading is toral if and only the rank of $L_e$ coincides with the rank of $L$.

If $Q$ is a MAD-group of $\aut L$, then    the homogeneous component $L_e$, that is,
the subalgebra fixed by $Q$, is an abelian subalgebra whose dimension coincides with the dimension of $Q$
 (by definition, the
  dimension of the maximal torus contained in $Q$) according to  \cite[Corollary~5]{f4}.
%
  Moreover, for any
quasitorus $Q$ of $\aut L$, the dimension of $L_e$ is, at least, the dimension of $Q$.
 Indeed, take $Q'$   a MAD-group containing $Q$, then $L_e$ contains $\fix Q'$
  and $\fix Q'$  is a subalgebra of dimension equal to $\dim Q'$, as above.
  Hence $\dim L_e\ge\dim \fix Q'\ge\dim Q$.
  In particular, if $L_e=0$ (if the grading is \emph{special}, by using terminology of
  \cite{Hesse}), the grading is produced by a finite quasitorus.
  The converse is true, in the case of a fine grading, again by  \cite[Corollary~5]{f4}: if the MAD-group  producing the grading
  (the diagonal group)
  is finite, then the grading is special.
  This kind of gradings has an extra-property:

 \begin{lm}\label{le_sobresemisimplesenfinas}
Every homogeneous element in a fine grading on a simple Lie algebra   produced by a finite quasitorus is semisimple.
\end{lm}

\textbf{Proof.}
 If $L=\oplus_{g}{L_g}$ is a grading, for any nonzero element $x\in L_g$, there is $x_s\in L_g$ semisimple
and $x_n\in L_g$ nilpotent such that $x=x_s+x_n$ and $[x_s,x_n]=0$ according to \cite[Th.~3.3]{enci}
(its semisimple and nilpotent parts also belong to $L_g)$.  Hence either every homogeneous element is semisimple or there
exists a homogeneous nilpotent element. The latter case does not happen if the grading is in the conditions of the lemma, because if $x\in L_g$
is nilpotent,  there are a semisimple
element $h\in L_e$ and  a nilpotent one $y\in L_{-g}$   such that $[h,x]=2x$, $[h,y]=-2y$ and $[x,y]=h$  by \cite[Th.~3.4]{enci}, so that
$L_e\ne0$.
 $\square$\smallskip

The general fact is that every homogeneous element in a fine grading on a simple Lie algebra is either semisimple or nilpotent
(\cite[Proposition~10]{f4}), which implies that it is very easy to find bases formed by semisimple or nilpotent elements.

It is worth pointing out that in this case of a grading  $L=\oplus_{g\in G}{L_g}$ on a simple Lie algebra $L$,
if $k$ is the Killing form of $L$ (non-degenerate), then $k(L_g,L_h)=0$ if $g,h\in G$ with $g+h\ne e$. Thus $L_{-g}$ can be identified with $L_g^*$.


\subsection{Some  techniques for group gradings}\label{subsec_detallestecnicos}

 We bring during this and the next subsections some key results from the structure of a MAD-group, extracted mainly from \cite{d4,f4,f4spin}.

 As in \cite[Section~5]{f4},   the Borel-Serre Theorem tell us that
  every quasitorus of $\mathfrak{G}:=\aut\eseis$ normalizes some of the maximal tori of   $\aut\eseis$.
 More precisely, the quasitorus is the product of a torus $T$ with a finite group, and we can take a maximal torus containing $T$
 such that the quasitorus is contained in its normalizer:

 \begin{lm}\label{le_puedoajustarlaparatedeltoro}\cite[Lemma~3]{f4spin}
If $H_1 $ is a toral subgroup of $\mathfrak{G}$ and $ H_2$ is a diagonalizable subgroup of $\mathfrak{G}$   which commutes with $H_1$, then there is 
a maximal torus $T$ of $\mathfrak{G}$ such that $H_1\subset T$ and $H_2$ is contained in the normalizer $\No(T)$.
\end{lm}

Assume we have fixed $\T$ a maximal torus of $\aut\eseis$. If $f\in\No(\T)$, then $f\T f^{-1}=\T$ and we define
\begin{equation}\label{eq_eltsuperf}
\begin{array}{l}
\tor{f}:=\cent_{\mathfrak{G}}(f)\cap\T=\{t\in\T\mid ftf^{-1}=t\}\\
Q(f):=\span{f,\tor{f}}
\end{array}
\end{equation}
where we use the notation $\span{S_1,\dots,S_l }$ for the quasitorus of $\aut\eseis$
generated by $S_1\cup\dots\cup S_l\subset\aut\eseis$ (the closure, with the Zarisky topology, of the group generated by them).
 The quasitori of the form $Q(f)$ have proved to be    relevant for the study of the MAD-groups in \cite{f4}, since every MAD-group of $\aut\f4$ is $Q(f)$ for certain
 $f\in\aut\f4$. This also happens for the MAD-groups of $\aut\d4$ containing outer automorphisms of order $3$, although not for   the MAD-groups containing outer automorphisms of order $2$ (for more details about these MAD-groups, see \cite{d4}).
 Through this paper we will also find several MAD-groups in the set $\{Q(f)\mid f\in\aut\eseis\}$, 
namely, all the five MAD-groups contained in $\Int\eseis$ (Section~\ref{sec_demostraciones}) and six of the MAD-groups not
contained in $\Int\eseis$ (Section~\ref{sec_demostracionesouter}).

\begin{lm}\label{le_comoeselTsuperf}
If  $f\in\No(\T)$ has   order $r\in\mathbb{N}$, then the set $\tor{f}$  is equal to $\sor{f} \cdot \hor{f}$ for the
 subtorus $\sor{f}=\{(tf)^r\mid t\in\T\}$ of $\T$ and a (finite) subgroup $\hor{f}\subset\{t\in \T\mid t^r=1_\mathfrak{G}\}$ such that $\sor{f} \cap \hor{f}=\{1_\mathfrak{G}\}$.
\end{lm}

\textbf{Proof.}
All is proved in \cite[Lemma~6]{f4spin}, except for the fact that the subtorus $S=\{(tf)^r\mid t\in\T\}\subset\sor{f}$ fills the whole
$\sor{f}$. If $s\in\sor{f}$, take $t\in\sor{f}$ such that $t^r=s$ (a torus contains roots of all its elements). As $\sor{f}\subset\tor{f}$, then $tf=ft$, so that $s=t^r=(tf)^r\in S$.
 $\square$\smallskip

 Although $\hor{f}$ is not determined by $f$ and $\T$, we will use such notation for any group satisfying that  $\tor{f}=\sor{f} \cdot \hor{f}$
 and $\sor{f} \cap \hor{f}=\{1_\mathfrak{G}\}$.

\begin{re}\label{re_paramover}
A way for identifying $Q(f)$ with $Q(g)$ for $f,g\in\No(\T)$ with the same projection on the Weyl group ($f\T=g\T$)
was developed  in \cite{f4}. Consider the group $\jor{f}:=\{f^{-1}t f t^{-1}\mid
t\in\T\}\subset\T$.
Denote by $\Ad F\colon \mathfrak{G}\to \mathfrak{G}$  the conjugation $\Ad F(h)=FhF^{-1}$ if $F\in \mathfrak{G}$.
If $s=f^{-1}t f t^{-1}\in\jor{f}$ ($t\in\T$), then $\Ad t(f)=fs$, so that $Q(f)$ is conjugated to $Q(fs)$, since $\Ad t$ does not move the torus $\T$.
As obviously $Q(f)=Q(ft)$ if $t\in\tor{f}$, then $Q(f)$ is conjugated to $Q(ft)$ for all $t\in\tor{f}\jor{f}$.
A sufficient   condition (also necessary) for having $\tor{f}\jor{f}=\T$ is that $\tor{f}\cap\jor{f}$ is finite, following the arguments in
the proof of \cite[Proposition~6]{f4}. This condition is easy to check in practice.
\end{re}

\begin{re}  \label{re_sobreordendelaextension}
If $f\in\No(\T)$, then $\tor{f}$ is finite if and only if all the elements $ft$ have the same order, independently of the element  $t\in\T$.
 This result and its  proof are almost the same than those ones in \cite[Lemma~1]{f4}, where we were working with the projections
 on the Weyl group.
 \end{re} 


 \subsection{Structure of a MAD-group of $\aut\eseis$}\label{sub_estructuraMADs}


 \begin{lm}\label{necesarioprimodivisor} (\cite[Lemma~4]{f4spin})
 If a prime $p$ does not divide the order of the (extended) Weyl group of $L$ for $L$ a simple Lie algebra,
 then every abelian $p$-group $H\le \aut L$ is toral.
 \end{lm}

 \begin{lm}\label{le_algunodelosfactoresnotoral} (\cite[Corollary~1]{f4spin})
 Any non-toral quasitorus of $\aut L $ for $L$ a simple Lie algebra   contains a non-toral $p$-group for some prime $p$.
 \end{lm}

 Now, take into account that
 \begin{itemize}
 \item the cardinal of the Weyl group of $\eseis$ is $3^42^65$;
 \item    there
 are no non-toral $5$-groups of $\aut\eseis$ (see Lemma~\ref{le_nohaydel5} afterwards or \cite[Lemma~10.3]{Griess} for the complex case);
 \end{itemize}
 to conclude that any non-toral quasitorus of $\aut\eseis$   contains either a non-toral $2$-group or a non-toral $3$-group.

 Our purpose is to go further:

 \begin{pr}\label{pr_contienepgrupoeltal}
If $Q$ is a MAD-group of $\aut\eseis$ contained in $\Int\eseis$, different from a maximal torus,
then $Q$
contains a non-toral group isomorphic to either $\Z_2^3$ or $\Z_3^2$.
\end{pr}

Hence any MAD-group contained in $\Int\eseis$ contains a  non-toral elementary $p$-group.
The point is that  it is not true that,  for any simple
  Lie algebra $L$, any non-toral quasitorus of $\aut L$  contains a  non-toral elementary $p$-group for some prime $p$.
   The condition of being maximal is necessary, as the following example shows:\smallskip

\noindent \textbf{ Example.}
By using the notations   in \cite{f4},
$$
Q=\langle\{t_{-1,1,-1,1}, t_{1,-1,-1,1},\widetilde{\sigma}_{105}t_{1,1,1,i}\}\rangle\le\aut\f4
$$
is a non-toral quasitorus isomorphic to $\mathbb{Z}_2^2\times \mathbb{Z}_8$ which   does not contain any $2$-elementary non-toral subquasitorus.  \smallskip

\noindent Examples of this situation in $\Int\eseis$
will appear in Section~\ref{sec_demostraciones} for $p=3$ and in Section~\ref{sec_demostracionesouter} for $p=2$.


  This makes necessary an ad-hoc proof of
    Proposition~\ref{pr_contienepgrupoeltal} (for $L=\eseis$ and $Q$ a MAD-group), which we have to postpone until Section~\ref{sec_demostraciones}.
    The purpose of the first part of this paper is to find the MAD-groups of $\mathfrak{G}_0:=\Int \eseis$ by using that they must live in the centralizers of such subgroups of type $\mathbb{Z}_2^3$ or $\mathbb{Z}_3^2$. But all the non-toral elementary $p$-subgroups in
    $\mathfrak{G}_0$ are well-known in the literature for the case $\mathbb{F}=\mathbb{C}$. In order to distinguish the complex case from the abstract one, we will denote  $\e6c$ the complex Lie algebra of type $E_6$ and $\mathfrak{G}^{\mathbb{C}}=\aut\e6c$
    and $\mathfrak{G}_0^{\mathbb{C}}=\Int\e6c$ the corresponding algebraic groups (in fact, Lie groups).

\begin{te}\label{teo_E6-data}
Let $p$ be a prime and $Q\leq \mathfrak{G}_0^{\mathbb{C}}=\Int \frak{e}_6^{\mathbb{C}} $ be a non-toral
elementary abelian $p$-subgroup. Then either $p=2$ or $p=3$, and $Q$ is
(up to conjugacy) one of the following subgroups:
\begin{enumerate}
\item If $p=2$:
\begin{itemize}
\item $V_2^3\cong \Z_2^3$ such that
$\cent_{\mathfrak{G}_0^{\mathbb{C}}}(V_2^3)=V_2^3\times\PSL(3)$,

\item $V_2^4\cong \Z_2^4$ such that
$\cent_{\mathfrak{G}_0^{\mathbb{C}}}(V_2^4)=V_2^3\times\GL(2)$,

\item $V_2^5\cong \Z_2^5$ such that
$\cent_{\mathfrak{G}_0^{\mathbb{C}}}(V_2^5)=V_2^3\times (\mathbb{C}^*)^2$.

\end{itemize}

\item If $p=3$:
\begin{itemize}
\item $V_3^{2a}\cong \Z_3^2$ such that
$\cent_{\mathfrak{G}_0^{\mathbb{C}}}(V_3^{2a})=V_3^{2a}\times\PSL(3)$,

\item $V_3^{2b}\cong \Z_3^2$ such that
$\cent_{\mathfrak{G}_0^{\mathbb{C}}}(V_3^{2b})=V_3^{2b}\times G_2 $,

\item $V_3^{3a}\cong \Z_3^3$ such that
$\cent_{\mathfrak{G}_0^{\mathbb{C}}}(V_3^{3a})=V_3^{2a}\times((\mathbb{C}^*)^2\rtimes \Z_3)$,

\item $V_3^{3b}\cong \Z_3^3$ such that
$\cent_{\mathfrak{G}_0^{\mathbb{C}}}(V_3^{3b})=V_3^{3b} \,\Z_3^3$,

\item $V_3^{3c}\cong \Z_3^3$ such that
$\cent_{\mathfrak{G}_0^{\mathbb{C}}}(V_3^{3c})=V_3^{2b}\times\SL(2)$,

\item $V_3^{3d}\cong \Z_3^3$ such that
$\cent_{\mathfrak{G}_0^{\mathbb{C}}}(V_3^{3d})=V_3^{2b}\times\GL(2)$,

\item $V_3^{4a}\cong \Z_3^4$ such that
$\cent_{\mathfrak{G}_0^{\mathbb{C}}}(V_3^{4a})=V_3^{4a}$,

\item $V_3^{4b}\cong \Z_3^4$ such that
$\cent_{\mathfrak{G}_0^{\mathbb{C}}}(V_3^{4b})=V_3^{2b}\times (\mathbb{C}^*)^2$.
\end{itemize}
\end{enumerate}
\end{te}

\noindent We are using the notation $\cent_AB=\{x\in A\mid xb=bx\,\forall b\in B\}$ for the centralizers ($B\le A$), and
$A=B\times C$ if  $A=BC$ for $B$ and $C$ normal subgroups of $A$ such that $B\cap C=\{1_A\}$. As usual, if only  $B$ is a normal
subgroup, the used notation is that one for semidirect product, $A=B\rtimes C$.

\begin{proof}
Either $p=2$ or $p=3$, as above.
Now, the data involving non-toral elementary abelian $2$-subgroups
in $\Int \e6c $ can be read in \cite[pg.\
78]{Moller-fundamenta-II}, where $V_2^3$, $V_2^4$ and $V_2^5$ are
denoted by $\langle V_2, c\rangle$, $\langle V_3, c\rangle$ and
$\langle V_4, c\rangle$ respectively. The reader may also consult
\cite[Theorem 8.2 and its proof, pg. 279]{Griess} where non-toral
elementary abelian $2$-subgroups in $\aut \e6c $ are
identified with non-toral elementary abelian $2$-subgroups in
$\aut \frak{f}_4^{\mathbb{C}} $. The non-toral   abelian $2$-subgroups in
$\aut \frak{f}_4 $ are  described in
\cite{f4}.

The data related to non-toral elementary abelian $3$-subgroups in
$\Int\e6c$ can be read in \cite[Theorem 8.10, pg.\
148]{Viru}.
\end{proof}


The important point is that the $p$-subgroups of $\mathfrak{G}_0^{\mathbb{C}}$ (respectively, of $\mathfrak{G}^{\mathbb{C}}$)
are
in bijective correspondence with the $p$-subgroups of $\mathfrak{G}_0$ (respectively, of $\mathfrak{G}$)
if the transcendence degree of the field extension  $\mathbb{F}\vert\mathbb{Q}$ is infinite.
But, even if this is not the case, we will find only one non-toral subgroup of type $\mathbb{Z}_2^3$ and
just two non-toral groups of type $\mathbb{Z}_3^2$ in $\mathfrak{G}_0$. This will be consequence of the following result, communicated to the authors by A.~Elduque,
which of course could be applied to any other simple Lie algebra.

\begin{pr}\label{pr_pasoacposarbitrarios}
There is an injective map
$$
\Omega\colon\{\text{finite subgroups of $\mathfrak{G}$}\}\rightarrow\{\text{finite subgroups of $\mathfrak{G}^{\mathbb{C}}$}\}.
$$
Furthermore, if $P_1$ and $P_2$ are two finite  subgroups of $\mathfrak{G} $,
$P_1$ and $P_2$ are conjugate in $\mathfrak{G} $ if and only if  $\Omega(P_1)$
and $\Omega(P_2)$ are    conjugate subgroups of $\mathfrak{G}^{\mathbb{C}}$.
\end{pr}

\textbf{Proof.}
Take $P=\{f_1,\dots,f_s\}$ a finite subgroup of $\mathfrak{G}$.
Fix $B=\{b_1,\dots,b_{78}\}$ a basis of $\eseis$ such that
 $[b_i,b_j]\in\sum_{b\in B}\mathbb{Q}b$ for all $i,j$ (for instance, take a Chevalley basis).
 Thus $\e6c\cong \hbox{\got e}_6^\mathbb{Q}\otimes_{\mathbb{Q}}\mathbb{C}$
 and $\eseis\cong \hbox{\got e}_6^\mathbb{Q}\otimes_{\mathbb{Q}}\mathbb{F}$,
 for $ \hbox{\got e}_6^\mathbb{Q}=\sum_{b\in B}\mathbb{Q}b$ the $\mathbb{Q}$-Lie algebra of type $E_6$ spanned by $B$.
 Note that, for each $f_i\in P$ and each $j\le78$, $f_i(b_j)=\sum_{k=1}^{78}c_{ijk}b_k$ for some $c_{ijk}\in\mathbb{F}$.
 Take $S=\{c_{ijk}\mid i\le s;\, j,k\le78\}$, which is  a finite subset of $\mathbb{F}$.
 Take $\mathbb{K}=\mathbb{Q}(S)$, a subfield of $\mathbb{F}$.
 Let $n$ be the transcendence degree of the field extension $\mathbb{K}\vert\mathbb{Q}$, which is obviously finite (as $S$).
 That means that there is an algebraically independent subset $\{X_1,\dots, X_n\}\subset \mathbb{K}$
 such that the extension $\mathbb{K}\vert\mathbb{Q}(X_1,\dots, X_n)$ is a finite (algebraic) extension.
 Every element of $\mathbb{K}$ is a root of some non-zero polynomial with coefficients in $\mathbb{Q}(X_1,\dots, X_n)$,
 concretely we can take $\{a_1,\dots,a_m\}\subset\mathbb{K}$ such that
 $\mathbb{K}=(\dots(\mathbb{Q}(X_1,\dots, X_n)(a_1))\dots)(a_m)$ and $p_i$ is the minimal polynomial of $a_i$ with coefficients in the field $(\dots(\mathbb{Q}(X_1,\dots, X_n)(a_1))\dots)(a_{i-1})$.
 Now recall that $\mathbb{C}\vert\mathbb{Q}$ has infinite transcendence degree, therefore there exists $\{x_1,\dots, x_n\}$
 an algebraically independent subset  of $\mathbb{C}$. It is clear that the fields $\mathbb{Q}(X_1,\dots, X_n)$ and $\mathbb{Q}(x_1,\dots, x_n)$ are isomorphic. Now we construct $K$ a subfield of $\mathbb{C}$ isomorphic to $\mathbb{K}=\mathbb{Q}(S)$, simply by adjoining to $\mathbb{Q}(x_1,\dots, x_n)$ elements $\{y_1,\dots,y_m\}\subset\mathbb{C}$ such that  $p_i$ is just the minimal polynomial of $y_i$ with coefficients in the field $(\dots(\mathbb{Q}(x_1,\dots, x_n)(y_1))\dots)(y_{i-1})\le\mathbb{C}$.
 Let $\Psi\colon\mathbb{K}\to K$ be a field isomorphism between $\mathbb{K}$ the subfield of $\mathbb{F}$ and $K$ the subfield of $\mathbb{C}$.
 Now define $\tilde f_i\colon\e6c\to\e6c$ by linearity, with $\tilde f_i(b_j)=\sum_{k=1}^{78}\Psi(c_{ijk})b_k\in\e6c$.
 It is clear that $\tilde f_i\in \mathfrak{G}^{\mathbb{C}}$ (it is an automorphism). Hence  $\Omega(P):=\{\tilde f_1,\dots,\tilde f_s\}$ is the desired subgroup of
 $\mathfrak{G}^{\mathbb{C}}$.

 Now assume we have $P_1=\{f_1,\dots,f_z\}$ and $P_2=\{g_1,\dots,g_z\}$ two finite subgroups of $\mathfrak{G}$.
 Thus $f_i(b_j)=\sum_{k=1}^{78}c_{ijk}b_k$ and $g_i(b_j)=\sum_{k=1}^{78}d_{ijk}b_k$ for some $c_{ijk},d_{ijk}\in\mathbb{F}$, and more concretely, these scalars  $c_{ijk}$ and $d_{ijk}$ live  in an extension $\mathbb{K}$ of $\mathbb{Q}$ such that there is an isomorphism $\Psi\colon\mathbb{K}\to K$ onto a subfield of $\mathbb{C}$. The fact of being $P_1$ and $P_2$ conjugate in $\mathfrak{G}$ is equivalent to the existence of $\{a_{jk}\mid j,k\le78\}\subset \mathbb{F}$ verifying certain polynomial equations. Namely, if $\varphi\colon\eseis\to\eseis$ is the automorphism such that $\varphi f_i\varphi^{-1}=g_i$ and  $A=(a_{jk})\in \rm{Mat}_{78\times78}(\mathbb{F})$ is its matrix relative to the basis $B$, the conditions for $\varphi$ are equivalent to the existence of solutions of the following polynomial equations, for each $i,j,l$:
 $$
 \begin{array}{l}
 \sum_k c_{ijk}a_{kl}=\sum_k d_{ikl}a_{jk}\\
 \sum_k \alpha_{ijk}a_{kl}=\sum_{s,t}a_{is}a_{jt}\alpha_{stl}\\
 (\mathop{\hbox{\rm det}}(A)+1)(\mathop{\hbox{\rm det}}(A)-1)=0
 \end{array}
 $$
if $\alpha_{ijk}\in\mathbb{Q}$ are such that $[b_i,b_j]=\sum_k\alpha_{ijk}b_k$.
So we have $I\triangleleft \mathbb{K}[Y_{ij}\mid i,j\le78]$ an ideal of polynomials such that the conjugacy of $P_1$ and $P_2$ is equivalent to the existence in $\mathbb{F}^r$ ($r=78^2$) of some common zero of all the polynomials in $I$.
But, as $\mathbb{F}$ is algebraically closed, by weak Nullstellensatz, this is equivalent to the fact that $1\in I$.
 By similar arguments the conjugacy of $\Omega(P_1)$ and $\Omega(P_2)$ is equivalent to the existence in $\mathbb{C}^r$   of some common zero of all the polynomials in $I$ (passing through $\Psi$)
and again  this is equivalent to the fact that $1\in I$.
\hfill  $\square$ \smallskip   \

In the following section we will provide descriptions of three quasitori of types $\mathbb{Z}_2^3$, $\mathbb{Z}_3^2$ and $\mathbb{Z}_3^2$  (two of them 
non-conjugate), which will be unique by Theorem~\ref{teo_E6-data} jointly with the above proposition.

  To finish this section, we would like to remark some   results that, though simple, will be highly useful during this paper:

  \begin{lm}\cite[Theorem 8.2.(3)]{Viru}\label{le_kaspersingeneralizar}
Let $\mathcal{G}$ be a linear algebraic group over an algebraically
closed field. Assume that $\mathcal{G}$ is a connected reductive group such that its commutator
subgroup   is simply connected. If $Q$ is a subquasitorus of
$\mathcal{G}$   generated by at most two elements,
then $Q$ is toral.
\end{lm}

This cannot be applied to our context, since $\Int\eseis$ has fundamental group $\Z_3$, so we will generalize it a little bit.

\begin{lm}\label{necesito3factores}
Let $\mathcal{G}$ be a linear algebraic group over an algebraically closed field.
Assume that $\mathcal{G}$ is a connected reductive group such that its commutator
group has fundamental group $\Z_n$. If $Q$ is a quasitorus of $\mathcal{G}$
generated by at most two elements and the order of $Q$ is prime to
$n$, then $Q$ is toral.
\end{lm}

\textbf{Proof.}
Let $\mathcal{G}$ be a linear algebraic group over an algebraically closed field,
and let $\mathcal{G}'$ denote its commutator subgroup $[\mathcal{G},\mathcal{G}]$. Then, the short exact
sequence of groups $\mathcal{G}'\to \mathcal{G} \to \mathcal{G}/\mathcal{G}'$ is indeed a fibration of connected
topological spaces that induces a short exact sequence of fundamental
groups $\pi_1(\mathcal{G}')\to \pi_1(\mathcal{G}) \to \pi_1(\mathcal{G}/\mathcal{G}')$. Therefore $\pi_1(\mathcal{G}')\cong
\Z_n$ can be identified with a subgroup of $\pi_1(\mathcal{G})$, and we denote by
$p\colon\widetilde{\mathcal{G}}\to \mathcal{G}$ the $n$-sheeted cover of $\mathcal{G}$ associated to
that subgroup. Then $\widetilde{\mathcal{G}}$ is again a connected reductive linear
algebraic group whose commutator subgroup is denoted by $\widetilde{\mathcal{G}}'$,
and $p$ is an algebraic epimorphism that identifies $\widetilde{\mathcal{G}}$ with a
central extension of $\mathcal{G}$ by $\Z_n$. We now claim that $\widetilde{\mathcal{G}}'$ is
simply connected. Indeed, since $p$ is a group epimorphism, its
restriction $p|_{\widetilde{\mathcal{G}}'}:\widetilde{\mathcal{G}}'\to \mathcal{G}'$ is so, and the
snake lemma shows that $Z_n\cong Z(\widetilde{\mathcal{G}})\subset \widetilde{\mathcal{G}}'$
is in the kernel of $p|_{\widetilde{\mathcal{G}}'}$ what shows that $\widetilde{\mathcal{G}}'$
is the universal cover of $\mathcal{G}'$.

Now, let $Q$ be a quasitorus of $\mathcal{G}$ generated by two elements and define
$\widetilde{Q}=p^{-1}(Q)$. Then $\widetilde{Q}$ is a the quasitorus of
$\widetilde{\mathcal{G}}$, that fits in a central short exact sequence of discrete
groups $\Z_n\to \widetilde{Q} \to Q$. Since $|Q|$ is coprime with $n$,
then $H^*(Q;\Z_n)=0$ and the previous exact sequence splits, that is
$\widetilde{Q}\cong\Z_n\times Q$ and $\widetilde{Q}$ can be generated by
two elements (as $Q$ is so).Then the results follows from Lemma~\ref{le_kaspersingeneralizar}
applied to $\widetilde{Q}$ and $\widetilde{\mathcal{G}}$.
\hfill  $\square$ \smallskip   \

Consequently there are no non-toral $2$-groups of $\Int\eseis$
with less than $3$ factors.

\begin{lm} \label{le_toroporalgotoral} (\cite[Lemma~2]{d4})
If $L$ is a simple Lie algebra, $T$ is a torus of $\aut L$ and $H$ is a toral subgroup of
$\aut L$ commuting with $T$, then $HT$ is toral.
\end{lm}

\begin{re}\label{re_cuandoQ(f)toralseguro}
An immediate consequence of Lemma~\ref{le_toroporalgotoral} is that,
if $\T$ is a maximal torus of $\aut\eseis$, and $f\in\No(\T)$ is an inner automorphism such that
$\tor{f}$ is a torus, then $Q(f)$ is toral. Moreover,
if $\tor{f}\cong(\F^*)^l\times\Z_2$, then $Q(f)$ is also toral,
by  applying Lemma~\ref{necesito3factores} and Lemma~\ref{le_toroporalgotoral}.
\end{re}


\section{Description of the inner gradings}\label{sec descripcionesinternos}

 We call an \emph{inner} grading that one produced by a quasitorus contained in $\Int\eseis$, the identity component of $\aut\eseis$.
 And we call an \emph{outer} grading a not-inner grading.
 Of course the fine inner gradings are produced by MAD-groups of $\aut\eseis$ contained in $\Int\eseis$. For describing them,
 first we fix some notation about the finite order inner automorphisms.

\subsection{Inner automorphisms of finite order}\label{subsec_autinternos}

Recall that the finite order automorphisms of the simple Lie algebras are completely described in \cite[Chapter~8]{Kac}.
The $\Z_m$-inner gradings on $\eseis$ can be obtained by assigning weights $\bar p=(p_0,\dots,p_6)$ ($p_i\in\Z_{\ge0}$)
to the nodes of the extended affine diagram of $\eseis$, $E_6^{(1)}$,
such that $\sum_{i=0}^6 p_in_i=m$, for $n_i$ the label of the corresponding node,
that is, $n_0=1$ and $\alpha_0=-\sum n_i\alpha_i$ denotes the opposite of the maximal root, for $\{\alpha_i\}_{i=1}^6$ a set of simple roots.
\vskip0.8cm

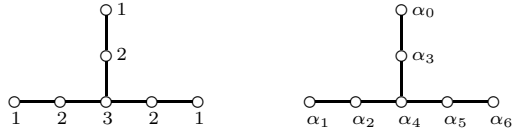
\begin{center}{
\begin{picture}(25,5)(4,-0.5)
\put(5,0){\circle{1}} \put(9,0){\circle{1}} \put(13,0){\circle{1}}
\put(17,0){\circle{1}} \put(21,0){\circle{1}}
\put(13,4){\circle{1}} \put(13,8){\circle{1}}
\put(5.5,0){\line(1,0){3}}
\put(9.5,0){\line(1,0){3}} \put(13.5,0){\line(1,0){3}} \put(13,4.5){\line(0,1){3}}
\put(17.5,0){\line(1,0){3}} \put(13,0.5){\line(0,1){3}}
\put(4.7,-2){$\scriptstyle 1$} \put(8.7,-2){$\scriptstyle 2$}
\put(12.7,-2){$\scriptstyle 3$} \put(16.7,-2){$\scriptstyle 2$}
\put(20.7,-2){$\scriptstyle 1$} \put(13.9,3.6){$\scriptstyle 2$} \put(13.9,7.6){$\scriptstyle 1$}
\end{picture}
\begin{picture}(25,5)(4,-0.5)
\put(5,0){\circle{1}} \put(9,0){\circle{1}} \put(13,0){\circle{1}}
\put(17,0){\circle{1}} \put(21,0){\circle{1}}
\put(13,4){\circle{1}} \put(13,8){\circle{1}}
\put(5.5,0){\line(1,0){3}}
\put(9.5,0){\line(1,0){3}} \put(13.5,0){\line(1,0){3}} \put(13,4.5){\line(0,1){3}}
\put(17.5,0){\line(1,0){3}} \put(13,0.5){\line(0,1){3}}
\put(4.7,-2){$\scriptstyle \alpha_1$} \put(8.7,-2){$\scriptstyle \alpha_2$}
\put(12.7,-2){$\scriptstyle \alpha_4$} \put(16.7,-2){$\scriptstyle \alpha_5$}
\put(20.7,-2){$\scriptstyle \alpha_6$} \put(13.9,3.6){$\scriptstyle \alpha_3$} \put(13.9,7.6){$\scriptstyle \alpha_0$}
\end{picture}
}\end{center}\vskip0.3cm

\noindent The subalgebra fixed by this automorphism is reductive of rank $6$, and the Dynkin diagram of its semisimple part is just the one obtained when removing the nodes with non-zero weights from $E_6^{(1)}$.

In particular, there are $5$ conjugacy classes of order three automorphisms.
We say that an order three automorphism $f$ is of type $3B$, or  that $f\in\,3B$, if $f$ is 
obtained with the choice of weights $(0,0,0,0,0,1,1)$ (for some Cartan subalgebra and some set of simple roots).
 Thus the fixed subalgebra is of type $ \mathfrak{a}_5$ plus a one-dimensional center, whose dimension is $36$.
The same notations will be used for the remaining order three automorphisms, according to the following table:\vskip0.2cm

\begin{center}{
\begin{tabular}{|c|c|c|c|}
 Type & Fixed subalgebra & dim &   $\bar p $\cr
\hline
 $3B$ & $\frak{a}_5+Z$  & $36$  & $(0,0,0,0,0,1,1)$  \cr
\hline
 $3C$ & $3\frak{a}_2$  & $24$  &  $(0,0,0,0,1,0,0)$  \cr
\hline
 $3D$ & $\frak{d}_4+2Z$  & $30$  &   $(1,1,0,0,0,0,1)$ \cr
\hline
 $3E$ & $\frak{a}_4+\frak{a}_1+Z$  & $28$  &  $(1,0,1,0,0,0,0)$ \cr
\hline
 $3F$ & $\frak{d}_5+Z$  & $46$  &  $(2,1,0,0,0,0,0)$ \cr
\hline
 \end{tabular}}\end{center}\vskip0.2cm

\noindent We have chosen the names of the types of the automorphisms according to \cite[Table~VI]{Griess}.
Besides, during this paper we use the notations $\mathfrak{a}_l$, $\mathfrak{b}_l$, $\mathfrak{c}_l$, $\mathfrak{d}_l$,
$\mathfrak{g}_2$, $\mathfrak{f}_4$ and $\mathfrak{e}_l$ for the simple finite-dimensional Lie algebras, instead of capital letters, to avoid confusions.

%
%
In the same way, there are two conjugacy classes of order two (inner) automorphisms:\vskip0.2cm

\begin{center}{
\begin{tabular}{|c|c|c|c|}
 Type & Fixed subalgebra & dim &   $\bar p$\cr
\hline
 $2A$ & $\frak{a}_5+\frak{a}_1$  & $38$  &  $(0,0,1,0,0,0,0)$ \cr
\hline
 $2B$ & $ \frak{d}_5+Z$  & $46$  & $(1,1,0,0,0,0,0)$ \cr
\hline
 \end{tabular}}\end{center}\vskip0.2cm

We observe that the type of any order three or two inner automorphism of $\eseis$
is determined by the dimension of its fixed subalgebra.


\subsection{A $\Z_3^4$-grading}\label{sub Z34gr}

Let $V_1=V_2=V_3=V$ be a  three-dimensional vector space. Let $B_V=\{u_0,u_1,u_2\}$ be a basis of $V$. Take
\begin{equation*}
    \mathcal{L}=\sll(V_1)\oplus\sll(V_2)\oplus\sll(V_3)\oplus\  V_1\otimes V_2\otimes V_3\oplus\   V_1^*\otimes V_2^*\otimes V_3^*
\end{equation*}
with the product given as in \cite[Chapter~13]{Adams}, which is a simple Lie algebra of type $ \mathfrak{e}_6$. Note that we have a $\Z_3$-grading on $\eseis\equiv\mathcal{L} $ by doing 
\begin{equation}\label{eq_modeloAdams}
\begin{array}{l}
\mathcal{L}_{\bar0}=\sll(V_1)\oplus\sll(V_2)\oplus\sll(V_3),\\
\mathcal{L}_{\bar1}= V_1\otimes V_2\otimes V_3,\\
\mathcal{L}_{\bar2}=V_1^*\otimes V_2^*\otimes V_3^*.
\end{array}
\end{equation}
Consider   $F_1\in\aut\frak{e}_6$ the order three grading automorphism, that is, $F_1\vert_{ \mathcal{L}_{\bar i}}=\omega^i\id$ for $\omega$ a primitive cubic root of the unit and $i=0,1,2$. Take $F_2$ the automorphism which \emph{permutes the $V_i$-components}, that is, the only one verifying $F_2(u\otimes v\otimes w)=v\otimes w\otimes u$ for all $u,v,w\in V$. As $F_2(\mathcal{L}_{\bar i})\subset \mathcal{L}_{\bar i}$, the automorphisms $F_1$ and $F_2$ commute.

Now, if $A\in\GL(3)$ and $f_A\in\GL(V)$ is the endomorphism whose associated matrix relative to $B_V$ is $A$,  we call $\Psi(A)$ the only automorphism of $\eseis$ whose action in $\mathcal{L}_{\bar1}$ is $u\otimes v\otimes w\mapsto f_A(u)\otimes f_A(v)\otimes f_A(w)$. Note that the uniqueness is a consequence of the fact that $\mathcal{L}_{\bar2}=[\mathcal{L}_{\bar1},\mathcal{L}_{\bar1}]$
 ($[\mathcal{L}_{\bar1},\mathcal{L}_{\bar1}]$ is an $\mathcal{L}_{\bar0}$-submodule of  $\mathcal{L}_{\bar2}$, which is irreducible)
 and $\mathcal{L}_{\bar0}=[\mathcal{L}_{\bar1},\mathcal{L}_{\bar2}]$. It is then easy to check that the only possible extension
 preserves the bracket.
Thus we have a well-defined map
\begin{equation}\label{eq_lapsidemodeloAdams}
\Psi\colon \GL(3)\to\rm{cent}_{\aut \mathfrak{e}_6}(F_1,F_2).
\end{equation}
with kernel $\{I_3,\omega I_3,\omega^2 I_3\}$, where  $I_n$ will denote the identity matrix of size $n$ throughout the text.
Take the following invertible matrices
$$
b:=\begin{pmatrix}1&0&0\cr 0&\omega&0 \cr 0&0&\omega^2
\end{pmatrix},\qquad
c:=\begin{pmatrix}0&0&1\cr 1&0&0 \cr 0&1&0
\end{pmatrix},
$$
 and $F_3:=\Psi(b)$ and $F_4:=\Psi(c)$ the related order three automorphisms.
As $bc=\omega cb$, thus $F_3$ and $F_4$ commute.
Consider the subgroup of automorphisms
$$
\Q_1:=\langle F_1,F_2,F_3,F_4\rangle\le\Int\eseis,
$$
which is isomorphic, as abstract group, to $\Z_3^4$.
If we denote a homogeneous component of the $\mathbb{Z}_3^4$-grading induced by $\Q_1 $ on $\mathcal{L}$ by
\begin{equation}\label{eq_componenteshomogdelaZ3ala4}
L_{(\bar i,\bar j,\bar k,\bar l)}=\{x\in\mathcal{L}\mid F_1(x)=\omega^{i}x,F_2(x)=\omega^{j}x,F_3(x)=\omega^{k}x,F_4(x)=\omega^{l}x\}
\end{equation}
for $i,j,k,l\in\{0,1,2\}$, and denote by $u_{ijk}$ the element $u_i\otimes u_j\otimes u_k$ in $\mathcal{L}_{\bar1}$,
then it is easy
to check that:
\begin{itemize}
\item
$\sum_{j} L_{(\bar 1,\bar j,\bar 0,\bar 0)}=\span{u_{000}+u_{111}+u_{222},u_{012}+u_{120}+u_{201},u_{021}+u_{210}+u_{102}}$
is fixed by $F_2$, so that it coincides with $L_{(\bar 1,\bar 0,\bar 0,\bar 0)}$, while
 $L_{(\bar 1,\bar j,\bar 0,\bar 0)}=0$ if $j=1,2$.
\item
 $\sum_{j} L_{(\bar 1,\bar j,\bar 0,\bar 1)}=\span{u_{000}+\omega^2u_{111}+\omega u_{222},
 u_{012}+\omega^2 u_{120}+\omega u_{201},u_{021}+\omega^2u_{210}+\omega u_{102}}$ is not fixed by $F_2$, which acts in the three
 elements with eigenvalues $1$, $\omega$ and $\omega^2$ respectively. The same happens with
 $\sum_{j} L_{(\bar 1,\bar j,\bar 0,\bar 2)}$, by duality.
 \item
$\sum_{j} L_{(\bar 1,\bar j,\bar 1,\bar 0)}=\span{x_0:=u_{001}+\omega^2u_{112}+\omega u_{220},
 x_1:=u_{010}+\omega^2 u_{121}+\omega u_{202},x_2:=u_{100}+\omega^2u_{211}+\omega u_{022}}$ is not fixed by $F_2$,
 but
 this automorphism maps each $x_i$ into $x_{i+1}$, so that it acts with eigenvalues $1,\omega,\omega^2$ in
 $x_0+x_1+x_2$, $x_0+\omega^2x_1+\omega x_2$ and $x_0+\omega x_1+\omega^2x_2$ respectively.
%
\end{itemize}
On the other hand, $\span{F_3,F_4}$ breaks each copy $\sll(V_i)$ of  $\mathcal{L}_{\bar0}$
in $8$ pieces of dimension one (the   non-toral $\Z_3^2$-grading on $\frak{a}_2$ usually called \emph{Pauli grading}, see for instance \cite{tipoA}), so,
$\sum_{j} L_{(\bar 0,\bar j,\bar k,\bar l)}$ has dimension three  if $(\bar k,\bar l)\ne (\bar 0,\bar 0)$ and breaks into three
one-dimensional pieces when $F_2$ is applied.
Therefore all the homogeneous components have dimension one, except the following cases:
 $$
 \begin{array}{l}
\dim L_{(\bar 0,\bar j,\bar 0,\bar 0)}=0\quad\text{ for any $j\in\{0,1,2\}$,}\\
\dim L_{(\bar i,\bar j,\bar 0,\bar 0)}=0\quad\text{  for any $i,j\in\{1,2\}$,}\\
\dim L_{(\bar i,\bar 0,\bar 0,\bar 0)}=3\quad\text{  for any $i\in\{1,2\}$.}
\end{array}
$$
Consequently the type of this grading is $(72,0,2)$.

In order to prove that we have found our first fine grading, note the following result:

\begin{lm}\label{le_soloelP1}
 The  non-toral quasitorus  $\mathcal{P}_1:=\span{F_1,F_2}$ has centralizer  $\cent_{\mathfrak{G}}(P_1)\cong P_1\times\PSL(3)$.
 So it has type $V_3^{2a}$   with the notations in Theorem~\ref{teo_E6-data}.
\end{lm}

\textbf{Proof.}
First we are going to prove that $\cent_{\mathfrak{G}}(\mathcal{P}_1)\cong \mathcal{P}_1\rm{Im}\,\Psi$, for the map $\Psi$ defined in Equation~(\ref{eq_lapsidemodeloAdams}).
It is clear that $\fix \mathcal{P}_1=\{x^1+x^2+x^3\mid x\in\sll(V)\}\cong\sll(V)$ is an algebra of type $ \mathfrak{a}_2$, where if $x\in\sll(V)$,
we are denoting by $x^i$ the element $x$ in $\sll(V_i)$, for $i=1,2,3$.
Take $F$   any automorphism
belonging to $\cent_{ \mathfrak{G}}(\mathcal{P}_1)$. This $F$ preserves all the homogeneous components of the $\mathbb{Z}_3^2$-grading
$\L=\oplus L_{(\bar i,\bar j)}$ produced by $\mathcal{P}_1$. Note that $\dim L_{(\bar 1,\bar 0)}=\dim L_{(\bar 2,\bar 0)}=11$ ($F_2\in3D$ and fixes a subalgebra of dimension $30$ isomorphic to $\mathfrak{d}_4+2Z$) but   $\dim L_{(\bar i,\bar j)}=8$ for the six remaining homogeneous components, which are $L_{(\bar 0,\bar 0)}=L_e$-irreducible modules of adjoint type.

In particular the map $F$ leaves $L_e=\fix \mathcal{P}_1 \cong\sll(V)$ invariant. Thus $F\vert_{L_e}\in\aut\sll(V)\cong\PSL(V)\rtimes\mathbb{Z}_
2$, concretely $\aut\sll(V)=\{\Ad f,\theta\Ad f \mid  f\in\GL(V) \}$,
for $\Ad f (x)=fxf^{-1}$  and for the outer order $2$ automorphism of $\sll(V)$ given by $\theta(x)=-x^t$.
If $F\vert_{L_e}\in\{\Ad A,\theta\Ad A\}$  for $A\in \GL(3)$ (identified with $\GL(V)$ by means of $B_V$), by replacing $F$ with
$F\Psi(A^{-1})$ we can assume that $F\vert_{L_e}\in\{\id,\theta\}$.

 Assume that $F\vert_{L_e}=\id$ and we are going to check that then $F\in\langle F_1,F_2\rangle$.
As $L_{(\bar 0,\bar 1)}$ is $L_e$-irreducible, the restriction $F\vert_{L_{(\bar 0,\bar 1)}}\in\hom_{L_e}(L_{(\bar 0,\bar 1)} , L_{(\bar 0,\bar 1)})=\mathbb{F}\id$, according to Schur's Lemma. So there is $\a\in\mathbb{F}^*$ such that $F\vert_{L_{(\bar 0,\bar 1)}}=\a\id$.
But $[[L_{(\bar 0,\bar 1)},L_{(\bar 0,\bar 1)}],L_{(\bar 0,\bar 1)}]=L_e$ ($F_2$ produces a $\mathbb{Z}_3$-grading on
$\L_{\bar0}$ with fixed subalgebra of type $ \mathfrak{a}_2$), so $\a^3=1$, and changing, if necessary, $F$ with either $FF_2$ or $FF_2^2$, we can assume that $F\vert_{L_{(\bar 0,\bar 1)}}=\id$. Hence $F\vert_{\L_{\bar0}}=\id$ (recall that $ {\L_{\bar0}}=  \oplus_j L_{(\bar 0,\bar j)}  $).
Again we can apply  Schur's Lemma, because $\L_{\bar1}$ is an  $\L_{\bar0}$-irreducible module, thus $F\vert_{\L_{\bar1}}=\beta\id$
for some nonzero scalar $\beta\in\mathbb{F}$. By similar arguments, $\beta^3=1$, and we change, if necessary, $F$ with $FF_1$ or $FF_1^2$
to obtain that $F\vert_{\L_{\bar1}}=\id$. From that we conclude that (after multiplying by elements in $\langle F_1,F_2\rangle$), $F=\id_\L$.

Now assume that we have $F\in \cent_{ \mathfrak{G}}(\mathcal{P}_1)$ such that    $F\vert_{L_e}=\theta$.
  Thus, $F^2\vert_{L_e}=\theta^2=\id$. As above, this implies that $F^2\in \mathcal{P}_1$, so that
$ \span{F_1,F_2,F}$ is isomorphic as abstract group to $\mathbb{Z}_3\times\mathbb{Z}_6$. Take $G\in \mathfrak{G}$ of order $6$ such that
$\span{G,F_2}=\span{F_1,F_2,F}$. Now $\fix (G\vert_{\fix F_2\cong\mathfrak{d}_4+2Z})=\fix\span{F_1,F_2,F}=\fix F\vert_{\fix \mathcal{P}_1\cong\sll(V)}$
is identified with the algebra $\fix \theta=\{x\in\sll(3)\mid x=-x^t\}=\mathfrak{so}\,(3)\cong\mathfrak{a}_1$. But $G\vert_{\mathfrak{d}_4}$ is an automorphism of $\mathfrak{d}_4$ of order $r$ a divisor of $6$, so that the possibilities for fixed subalgebras of rank different from $4$ are just $2\mathfrak{a}_1$ and $\mathfrak{a}_1+Z$ if $r=6$, $\mathfrak{g}_2$ and $\mathfrak{a}_2$ if $r=3$ and $\mathfrak{b}_3+Z$ if $r=2$, of course not contained in an algebra of type $ \mathfrak{a}_1$. We have got a contradiction in this case.

Note that  we have really proved $\cent_{ \mathfrak{G}}(\mathcal{P}_1)=\mathcal{P}_1\Psi(\SL(3))$. But this is a direct product,
in spite of the fact $F_1\in\rm{Im}\,\Psi$ (if $\xi$ is a 9th root of unit such that $\xi^3=\omega^2$, then $\Psi(\xi^2 I_3)=F_1$).
Indeed, take $A\in\SL(3)$ and $n_1,n_2\in\{0,1,2\}$
such that $\Psi(A)=F_1^{n_1}F_2^{n_2}$.
 On one hand, the fact $\Psi(A)\vert_{\L_{\bar0}}=\id$ implies that $A$ commutes with $\sll(3)$ and hence is in the center of gl$(3)$, this center equal to $\mathbb{F}I_3$, so that there is $\gamma\in\F$ such that $A=\gamma I_3$. As $\det A=1$, then $\gamma^3=1$.
On the other hand, the elements  $u_{ijk}\in\L_{\bar1}=V^{\otimes3}$ must be eigenvectors of $F_1^{n_1}F_2^{n_2}$, thus $n_2=0$. But $\Psi(A)(u_{ijk})=\gamma^3 u_{ijk}=u_{ijk}=F_1^{n_1}(u_{ijk})$, so $\omega^{n_1}=1$ and $n_1=0$. In conclusion, $\cent_{ \mathfrak{G}}(\mathcal{P}_1)=\mathcal{P}_1\times\Psi(\SL(3))\cong\mathbb{Z}_3^2\times \PSL(3)$ since the kernel of $\Psi$ is $\{I_3,\omega I_3,\omega^2 I_3\}$.
\hfill  $\square$ \smallskip   \

Hence $\Q_1$ is a MAD-group of $\aut\eseis$, isomorphic to $\Z_3^4$, since $\langle b,c\rangle$ is a MAD-group of $\PSL(3)$.
Moreover, this grading is fine not only as a group-grading, but it is possible to check that it is also fine as a general grading (as partition into subspaces such that the product of two of them is contained in some other).

Now we would like to pay attention on some subquasitori of $\Q_1$ for further use.

\begin{re}\label{re_unZ3cubotoral}
The subquasitorus $\span{F_2,F_3,F_4}\cong\mathbb{Z}_3^3$ of $\Q_1$ is toral.
To check it, recall that if $L_g$ is a homogeneous component of the grading induced by $\Q_1$,
then $L_g$ and $L_{-g}$ are composed by semisimple elements by Lemma~\ref{le_sobresemisimplesenfinas}.
As $[L_g,L_{-g}]\subset L_e=0$, then
the elements in $L_g\oplus L_{-g}$ are semisimple too. Hence, the subalgebra fixed by
$\span{F_2,F_3,F_4}$, that is   $\mathfrak{h} =L_{(\bar 1,\bar 0,\bar 0,\bar 0)}\oplus L_{(\bar 2,\bar 0,\bar 0,\bar 0)}$,
is a toral subalgebra. As $ \mathfrak{h}$ has dimension $6$, it is a Cartan subalgebra.
\end{re}

\begin{re}\label{re_laJOrdangrading}
Observe that:\begin{itemize}
\item $\span{F_1,F_3,F_4}\cong\Z_3^3$ is a nontoral quasitorus of $\Q_1$, the Jordan subgroup in \cite{Alek} (appearing also in \cite[Chapter\,3, \S3.13]{enci}).
    \item It is a minimal non-toral elementary $3$-group of $\mathfrak{G}_0$. In particular it does not contain any non-toral group isomorphic to $\mathbb{Z}_3^2$.
        \end{itemize}
As $\sum_{j}\dim L_{(\bar 0,\bar j,\bar 0,\bar 0)}=0$,
$\sum_{j}\dim L_{(\bar i,\bar j,\bar 0,\bar 0)}=3+0+0=3$ if $i=1,2$,
and
 $\sum_{j}\dim L_{(\bar i,\bar j,\bar k,\bar l)}=1+1+1=3$ if $(k,l)\ne(0,0)$,
 this quasitorus induces a non-toral grading of type $(0,0,26)$ with $L_e=0$, which is just the Jordan grading in
 \cite[Main Theorem~(vi)]{AlbJordangrad}.
 The second item is a consequence of the fact that $L_g\oplus L_{-g}$ is a Cartan subalgebra for all $0\ne g\in\Z_3^3$,
by reasoning as in Remark~\ref{re_unZ3cubotoral}.
 \end{re}

We will use these notations during the whole section.


\subsection{A $\Z_3^2\times\Z^2$-grading}\label{sub Z32Zala2gr}


Lemma~\ref{le_soloelP1} suggests us another interesting grading.
Take for any scalars $\alpha,\beta\in\F^*$, the automorphism
$T_{\alpha,\beta}=\Psi\left(p_{\alpha,\beta}\right)$, for
$$
p_{\alpha,\beta}=\begin{pmatrix}\alpha&0&0\cr 0&\beta&0 \cr 0&0&\frac1{\alpha\beta}
\end{pmatrix}.
$$
The quasitorus
$$
\Q_2:=\langle\{F_1,F_2,T_{\alpha,\beta}\mid \alpha,\beta\in\F^*\}\rangle\le\Int\eseis
$$
 is
isomorphic to $\Z_3^2\times(\F^*)^2$.

\begin{lm}\label{le_elotroZ34}
The quasitorus $R_1=\span{F_1,F_2,T_{\omega,1},T_{\xi,\xi}}\cong\mathbb{Z}_3^4$  for $\xi$ a ninth root of the unit such that $\xi^3=\omega^2$ is of type $V_3^{4b}$, with the notations in Theorem~\ref{teo_E6-data}. Its centralizer in $\mathfrak{ G}$ is just $\Q_2$.
\end{lm}

\textbf{Proof.}
As $\mathcal{P}_1\subset R_1$, then $\cent_{\mathfrak{G}}R_1\subset \cent_{\mathfrak{G}}\mathcal{P}_1=\mathcal{P}_1\times \Psi(\SL(3))\cong\mathbb{Z}_3^2\times\PSL(3)$.
More concretely, $ \cent_{\mathfrak{G}}R_1=\mathcal{P}_1\times\cent_{\Psi(\SL(3))}\span{T_{\omega,1},T_{\xi,\xi}}$.
But if $A\in\SL(3)$ commutes with $p_{\omega,1}$, then $A$ is diagonal and hence it belongs to $\{p_{\alpha,\beta}\mid \alpha,\beta\in\F^*\}$.
\hfill  $\square$ \smallskip   \

Consequently $\mathcal{Q}_2$ is again a MAD-group of $\aut\eseis$, because $\cent_{\mathfrak{G}} \mathcal{Q}_2\subset\cent_{\mathfrak{G}} R_1=\mathcal{Q}_2$, that is, it is self-centralizing.

Our description makes it easy the computation of the simultaneous diagonalization.
If we denote by  $L_{(\bar i,\bar j,\gamma)}$ the set of elements of $\L $  in which $F_1$ acts  with eigenvalue  $\omega^i$, $F_2$ with eigenvalue $\omega^j$ and
 $T_{\alpha,\beta}$ with eigenvalue $\gamma$, then
we check that all these homogeneous components are zero except for
\begin{equation}\label{dimesniones}
\begin{array}{l}
\dim L_{(\bar i,\bar j,1)}=2,\\
\dim L_{(\bar i,\bar j,\gamma)}=1\quad\text{if $\gamma\in\{(\alpha^2\beta)^{\pm1},(\alpha\beta^2)^{\pm1},(\alpha/\beta)^{\pm1}\}$},\\
\dim L_{(\bar 1,\bar 0,\gamma)}=1\quad\text{if $\gamma\in\{\alpha^3,\beta^3,1/(\alpha^3\beta^3) \}$},\\
\dim L_{(\bar 2,\bar 0,\gamma)}=1\quad\text{if $\gamma\in\{1/\alpha^3,1/\beta^3, \alpha^3\beta^3 \}$},
\end{array}
\end{equation}
for all $i,j\in\{0,1,2\}$, so that   $\Q_2$  produces a fine grading of type $(60,9)$.

\begin{lm}\label{le_elP1P2}
 The quasitorus    $\mathcal{P}_2:=\span{F_1T_{\xi,\xi},F_2}$ fixes a subalgebra of type $ \mathfrak{g}_2$, so it
 is  non-toral of type    $V_3^{2b}$, with the notations in Theorem~\ref{teo_E6-data}.

Note that $\Q_2=\mathcal{P}_1\times T_2=\mathcal{P}_2\times T_2$, for $T_2=\span{T_{\alpha,\beta}\mid \alpha,\beta\in\F^*}$.
\end{lm}

Before proving this lemma, observe something about the isomorphy classes of the automorphisms in $\mathcal{P}_2$,
as well as in several quasitori.

\begin{re}\label{re_composiciones}
As $\dim\fix F_2=\sum_{i,k,l}\dim L_{(\bar i,\bar 0,\bar k,\bar l)}=24+2\cdot3=30$  and
$\dim\fix F_3=\sum_{i,j,l}\dim L_{(\bar i,\bar j,\bar 0,\bar l)}=18+2\cdot3=24$, then $F_2\in\,3D$
and $ F_3 \in\,3C$.
Similar arguments tell us  that all the  non-trivial automorphisms  in $\mathcal{P}_2$ are of the class $3D$
and that the only automorphisms in $\Q_1$ which are not  of type $3C$
are $F_2F$ and $F_2^2F$
for all $F\in \span{ F_3,F_4}$.
Therefore:
\begin{itemize}
\item  $\mathcal{P}_2$ is of type $D^8$ (this notation means that it contains $8$ automorphisms of type $3D$, of course joint with the identity),
\item $\mathcal{P}_1$ is of type $C^6D^2$,
\item $\mathcal{Q}_1$ is of type $C^{62}D^{18}$,
\item $\span{F_2,F_3,F_4}$ (the Jordan group in Remark~\ref{re_laJOrdangrading})
is of type $C^{26}$.
\end{itemize}
In particular this provides a direct way
of knowing when a non-toral subgroup $Q$ of two commuting order three automorphisms 
is conjugated to  either $\mathcal{P}_1$ or $\mathcal{P}_2$.
\end{re}\smallskip

\textbf{Proof.}
Taking into account Equation~(\ref{dimesniones}), $\dim\fix \mathcal{P}_2=\dim L_{(\bar 0,\bar 0,1)}+\dim L_{(\bar 1,\bar 0,\omega^2)}+\dim L_{(\bar 2,\bar 0,\omega)}=4+5+5=14$. We try to know more about this subalgebra.
As $F_2$ is of type $3D$, the subalgebra $\fix F_2$ is of type $\mathfrak{d}_4$ summed with a two-dimensional center. Now $F_1T_{\xi,\xi}$ preserves this subalgebra and its derived subalgebra, that is, $\mathfrak{d}_4$, producing a $\mathbb{Z}_3$-grading on
$\mathfrak{d}_4$. This implies that the restriction $F_1T_{\xi,\xi}\vert_{\mathfrak{d}_4}$ must fix a subalgebra of some of the types
$\{\mathfrak{a}_2,\mathfrak{g}_2,3\mathfrak{a}_1+Z,\mathfrak{a}_3+Z\}$, of dimensions $\{8,14,10,16\}$ respectively. As
$\fix \mathcal{P}_2=\fix F_1T_{\xi,\xi}\vert_{\mathfrak{d}_4+2Z}$ is equal to $\fix F_1T_{\xi,\xi}\vert_{\mathfrak{d}_4}$ summed with some abelian subalgebra of dimension either 0, or 1 or 2,
by dimension count the only possibility is that $\fix \mathcal{P}_2$ is a subalgebra of type $ \mathfrak{g}_2$.
Hence $\mathcal{P}_2$
is  non-toral of type    $V_3^{2b}$, since a quasitorus of $\aut\eseis$ is non-toral when its fixed subalgebra has rank different from $6$.
 \hfill  $\square$ \smallskip   \

In the complex case, Theorem~\ref{teo_E6-data} tells us that $\cent_{\mathfrak{G_0}^{\mathbb{C}}}\mathcal{P}_2\cong\mathcal{P}_2\times G_2$ (for $G_2$ the automorphism group of the octonion algebra), but we cannot translate this result directly to arbitrary fields by applying Proposition~\ref{pr_pasoacposarbitrarios}, since   $G_2$ is of course not finite. Besides we are interested in the centralizers in $\mathfrak{G}$, not in $\mathfrak{G}_0$. In order to understand why it works for an arbitrary $\F$, see first the action of  $\fix \mathcal{P}_2$ on the remaining homogeneous components.

 \begin{re}\label{re_sobrelasZ3cuadradogradsnotorales}
 As $\mathcal{P}_2$ is of type $D^8$,  the $\Z_3^2$-grading $\Gamma_0$ induced by $\mathcal{P}_2$ must have all the non-identity components of the same dimension, $(78-14)/8=8$ (alternatively see Equation~(\ref{dimesniones})). As $L_e$ is isomorphic to $\der\mathcal{C}$, the Lie algebra of derivations of an octonion algebra $\mathcal{C}$, and it does not act trivially on any homogeneous component, then all of them are isomorphic as $L_e$-modules to the $\der\mathcal{C}$-module $\mathcal{C}$ (sum of $\mathcal{C}_0$, the irreducible
  set of zero trace  octonions, with the trivial module $\F1$).

  In particular we can compute the type of the grading induced by $\Q_2$
 without doing the simultaneous diagonalization: as the root decomposition of $\frak{g}_2$ is of type $(12,1)$ and breaks each module
 $\mathcal{C}$ in $(6,1)$, hence the grading induced by $\Q_2$ has type $(12,1)+8(6,1)=(60,9)$, as we already knew.
 \end{re}

 Checking that $\cent_{\mathfrak{G}}\mathcal{P}_2\cong\mathcal{P}_2\times G_2$ is equivalent to checking that:
 \begin{itemize}
 \item[a)] If $F\in\mathfrak{G}$ preserves the homogeneous components of $\Gamma_0$ and $F\vert_{L_e}=\id$, then $F\in\mathcal{P}_2$.
 \item[b)] For any $f\in\aut\mathcal{C}=G_2$, there is an extension $\tilde f\in\aut\L$ such that $\tilde f(d)=fdf^{-1}$ if $d\in\der\mathcal{C}=L_e$ and $\tilde f$ preserves the homogeneous components of $\Gamma_0$.
 \end{itemize}
 If  we identify $\rho\colon\L\rightarrow\L':=\der\mathcal{C}\oplus\mathcal{C}^{(g_1)}\oplus\dots\oplus\mathcal{C}^{(g_8)}$ by means of the $L_e$-isomorphisms of modules of each component ($\mathcal{C}^{(g_i)}$'s are several copies of $\mathcal{C}$, indexed in a set with eight elements, $\mathbb{Z}_3^2\setminus\{e\}$), this map $\rho$ is by construction an $L_e$-isomorphism of modules, which allows to endow $\L'$ with a Lie algebra structure of type $E_6$, when asking for $\rho$ to be a Lie homomorphism.
 Taking into consideration that $\dim\hom_{\der\mathcal{C}}(\mathcal{C}_0\otimes\mathcal{C}_0,\mathcal{C}_0)=\dim\hom_{\der\mathcal{C}}(\F\otimes\mathcal{C}_0,\mathcal{C}_0)=1$, there must exist  some fixed nonzero scalars $\a_{ij},\beta_{ij},\gamma_i\in\F^*$, with $\a_{ij}=\a_{ji} $, such that,
 \begin{equation*}
 \begin{array}{l}
 {[x^{(g_i)},y^{(g_j)}]_{\L'}}=\a_{ij}(xy-yx)^{(g_i+g_j)} \qquad\text{if $g_i\ne 2g_j$},   \\
 {[1^{(g_i)},y^{(g_j)}]_{\L'}}= \beta_{ij}y^{(g_i+g_j)} \qquad\text{if $g_i\ne 2g_j$},\\
 {[x^{(g_i)},y^{(2g_i)}]_{\L'}}= \gamma_i D_{x,y}=\gamma_i([l_x,l_y]+[l_x,r_y]+[r_x,r_y])\in\der\mathcal{C},\\
 {[1^{(g_i)},y^{(2g_i)}]_{\L'}}=0,
 \end{array}
 \end{equation*}
 for any
 $x,y\in \mathcal{C}_0$, where  $l_x$ and $r_x$ denote respectively the left and right multiplication operators on $\mathcal{C}$.
 The scalars can be determined by passing through $\rho$ or simply by using the Jacobi identity
 (this provides an interesting model of a Lie algebra of type $E_6$, not to be developed in this paper). Now, the map $\tilde f$ which
 maps $d\mapsto fdf^{-1}$ if $d\in\der\mathcal{C}$  and $x^{(g_i)}\mapsto f(x)^{(g_i)}$ if $x\in\mathcal{C}$,
 is the required automorphism in item b).
 Finally, the statement in a) is proved with arguments as used in Lemma~\ref{le_soloelP1}, with caution, because $\mathcal{C}$ is not  $\der\mathcal{C}$-irreducible.


\subsection{A $\Z_3^2\times\Z_2^3$-grading}\label{sub Z32Z23}

Let $G_1\in 2A$ any order two automorphism of $\mathcal{M}=\eseis$ fixing an algebra of type $ \mathfrak{a}_5\oplus  \mathfrak{a}_1$. Let
$\mathcal{M}=\mathcal{M}_{\bar0}\oplus  \mathcal{M}_{\bar1}   $ the $\Z_2$-induced grading.
Hence there exist
$U$ and $W$
 vector spaces of dimensions $2$ and $6$ respectively such that
 $$\begin{array}{l}
 \mathcal{M}_{\bar0}=\sll(W)\oplus\sll(U),\\
 \mathcal{M}_{\bar1}=\wedge^3W\otimes U.
 \end{array}
 $$
 We will introduce some notations. If $f$ is an automorphism of the vector space $W$,
 we denote by $f^{\wedge3}$ the automorphism of the vector space $\wedge^3W$ mapping
$w^{1}\wedge w^{2}\wedge w^{3}$ ($w^{i}\in W $) into $f(w^{1})\wedge f(w^{2})\wedge f(w^{3})$. And, if $f\in\aut A$ and $g\in\aut B$ ($A$ and $B$ vector spaces),
we denote by $f\otimes g$ the automorphism of $A\otimes B$ such that $(f\otimes g)(a\otimes b)=f(a)\otimes g(b)$ for any $a\in A, b\in B$.
Now fix $B_U=\{u_0,u_1\}$ and $B_W=\{w_i\mid i=0,\dots,5\}$ two   bases of $U$ and $W$ respectively, and take
 $H_1$ and $H_2$ the only automorphisms of $\mathcal{M}$ whose restrictions to  $\mathcal{M}_{\bar1}   $ are
 $$
 H_1\vert_{\mathcal{M}_{\bar1}}=\begin{pmatrix}b&0\cr 0&-b
\end{pmatrix}_W^{\wedge3}\otimes\begin{pmatrix}1&0\cr 0&-1
\end{pmatrix}_U
$$
and
$$
 H_2\vert_{\mathcal{M}_{\bar1}}= \begin{pmatrix}0&c\cr c&0
\end{pmatrix}_W^{\wedge3}\otimes\begin{pmatrix}0&1\cr 1&0
\end{pmatrix}_U,
 $$
where we are identifying the automorphisms of $W$ (respectively $U$) with their matrices relative to $B_W$ (respectively $B_U$).

Note that $H_1$ and $H_2$ are order six automorphisms commuting with $G_1$ and between them, so that we
can consider the quasitorus
$$
\Q_3:=\langle\{H_1,H_2,G_1\}\rangle\le\Int\eseis
$$
  isomorphic, as abstract group, to $\Z_3^2\times\Z_2^3$.
  Let us prove that  $\Q_3$ is another MAD-group.
  For that, let us look at their $p$-subgroups.

\begin{lm}\label{le_elP3}
$\mathcal{P}_3=:\span{H_1^3,H_2^3,G_1 }$ is a non-toral quasitorus isomorphic to $\Z_2^3$, hence
   of type $V_2^{3}$.
\end{lm}

\textbf{Proof.}
An element fixed by $G_1$ belongs to $\sll(U)\oplus\sll(W)$. If $x_U\in\sll(U)$ is fixed by $H_1^3$ and $H_2^3$, then $x_U=0$.
Now, if we write $x_W\in\sll(W)$ in square matrix blocks as $\tiny{\begin{pmatrix}A&B\cr C&D
\end{pmatrix}}$, the fact that $x_W$ commutes with
$\tiny{\begin{pmatrix}I_3&0\cr 0&-I_3
\end{pmatrix} =\begin{pmatrix}b&0\cr 0&-b
\end{pmatrix}^3}$
forces $B$ and $C$ to be zero,
and the fact that $x_W$ commutes with $\tiny{ \begin{pmatrix}0&I_3\cr  I_3&0
\end{pmatrix}=\begin{pmatrix}0&c\cr  c&0
\end{pmatrix}^3}$ forces $A=D$.
As $0=\tr(A)+\tr(D)=2\tr(A)$, hence the fixed subalgebra
$\fix(\mathcal{P}_3)=\left\{\tiny{\begin{pmatrix}A&0\cr 0&A
\end{pmatrix}_W}\mid A\in \sll(3)\right\}$
is isomorphic to an algebra of type $ \mathfrak{a}_2$ and $\mathcal{P}_3$ is non-toral.
\hfill  $\square$ \smallskip

Following similar arguments to Lemma~\ref{le_soloelP1}, it is not difficult to find the
 centralizer  $\cent_{\mathfrak{G}}\mathcal{P}_3\cong \mathcal{P}_3\times\PSL(3)$.
 The idea is to consider the well-defined map
 $$
 \Psi'\colon\SL(3)\rightarrow\rm{cent}_{\mathfrak{G}}\mathcal{P}_3
 $$
 where, if $A\in\SL(3)$, $\Psi'(A)$ is the only automorphism of $\mathcal{M}$ whose restriction to $\mathcal{M}_{\bar1}$ is
 \begin{equation}\label{eq_elpsiprima}
 \Psi'(A) \vert_{\mathcal{M}_{\bar1}}=\begin{pmatrix}A&0\cr 0&A
\end{pmatrix}_W^{\wedge3}\otimes I_V
\end{equation}
again with the identifications between automorphisms of $W$ and $U$ and matrices relative to $B_W$ and $B_U$.

\begin{lm}\label{le_elP4}
$\mathcal{P}_4:=\span{H_1^2,H_2^2}$ is a non-toral quasitorus isomorphic to $\Z_3^2$
   of type $V_3^{2b}$.
\end{lm}

\textbf{Proof.}
We have to find the fixed part by the automorphisms which are extensions of
$
 \tiny{\begin{pmatrix}b^2&0\cr 0& b^2
\end{pmatrix}_W^{\wedge3}\otimes I_V}
$ and $\tiny{\begin{pmatrix}c^2&0\cr 0&c^2
\end{pmatrix}_W^{\wedge3}\otimes I_V}
$.
By reordering some rows and columns in the matrices relative to the endomorphisms of $W$, we can work with
$\tiny{
 \begin{pmatrix}I_2&0&0\cr 0&\omega I_2&0 \cr 0&0&\omega^2 I_2
\end{pmatrix}}$
and $\tiny{
 \begin{pmatrix}0&0&I_2\cr I_2&0&0 \cr 0&I_2&0
\end{pmatrix}}.
$
So, a block matrix 
$A=(A_{ij})_{i,j=1,2,3}$ ($A_{ij}\in\rm{Mat}_{2\times2}(\F)$)
commutes with them if and only if $A_{11}=A_{22}=A_{33}$ and $A_{ij}=0$ if $i\ne j$. As
$0=\tr(A)=3\tr(A_{11})$, then there is a two-dimensional vector subspace $W'$ of $W$ such that
we can identify
$\fix(\mathcal{P}_4)\cap \sll(W)$ with $\sll(W')$. Besides all the elements in $\sll(V)$ remain fixed, hence
$\fix(\mathcal{P}_4)\cap\mathcal{L}_{\bar0}=\sll(W')\oplus\sll(V)\cong\frak{a}_1\oplus \frak{a}_1$.
It is now easy to check that  $\fix(\mathcal{P}_4)\cap\mathcal{L}_{\bar1}=S^3(W')\oplus V$,
turning out that $\fix(\mathcal{P}_4)\cong2\frak{a}_1\oplus V(3)\otimes V(1)$, which is a Lie algebra isomorphic to $\frak{g}_2$ (a well-known fact, see for instance \cite[Theorem~3.2]{modelosg2}).
\hfill  $\square$ \smallskip

As $\mathcal{P}_4$ must be conjugated to $\mathcal{P}_2$ by Proposition~\ref{pr_pasoacposarbitrarios}, we conclude that also the centralizer
$\cent_{ \mathfrak{G}}(\mathcal{P}_4)$ is a direct product of $\mathcal{P}_4$ with a copy of the group $G_2$.
Consequently, as $\Q_3=\mathcal{P}_4\times \mathcal{P}_3$ lives in this centralizer, and the subquasitorus
$\mathcal{P}_3\cong\Z_2^3$ is known to be necessarily a MAD-group of $G_2=\aut\frak{g}_2$ (see, for instance, \cite{g2}),
these arguments imply that $\Q_3$ is a MAD-group of $\aut\eseis$.

 In this occasion, we compute the type of the grading induced by $\Q_3$ without doing
  the simultaneous diagonalization (not difficult, but long)
 but taking into account Remark~\ref{re_sobrelasZ3cuadradogradsnotorales} applied to  $\mathcal{P}_4$.
 It is well-known that the $\Z_2^3$-grading on the octonion algebra  $\mathcal{C}$ is a grading of type $(8)$
 which induces one of type
 $(0,7)$ on $\der\mathcal{C}$ (each non-trivial component is a Cartan subalgebra).
 Hence the fine grading induced by $\Q_3$ has type $(0,7)+8(8,0)=(64,7)$.


\subsection{A $\Z_2^3\times\Z^2$-grading}\label{sub Z23Zala2gr}

We use here the notation in the previous subsection. Take $G_2=H_1^3$ and $G_3=H_2^3$, two order two automorphisms whose restrictions to
$\mathcal{M}_{\bar1}   $ are, of course,
 $$\begin{array}{l}
 G_2\vert_{\mathcal{M}_{\bar1}}=\begin{pmatrix}I_3&0\cr 0&-I_3
\end{pmatrix}_W^{\wedge3}\otimes\begin{pmatrix}1&0\cr 0&-1
\end{pmatrix}_V,\\
 \ \\
 G_3\vert_{\mathcal{M}_{\bar1}}= \begin{pmatrix}0&I_3\cr I_3&0
\end{pmatrix}_W^{\wedge3}\otimes\begin{pmatrix}0&1\cr 1&0
\end{pmatrix}_V,
 \end{array}
 $$
and $S_{\alpha,\beta} \in\aut \mathcal{M} $ the automorphism given by
$$
S_{\alpha,\beta} \vert_{\mathcal{M}_{\bar1}}=\begin{pmatrix}p_{\alpha,\beta}&0\cr 0&p_{\alpha,\beta}
\end{pmatrix}_W^{\wedge3}\otimes I_V,
$$
that is, $S_{\alpha,\beta} =\Psi'(p_{\alpha,\beta})$ with the notation in Equation~(\ref{eq_elpsiprima}).
Take, then,
$$
\Q_4:= \langle\{ G_1,G_2,G_3,S_{\alpha,\beta}\mid \alpha,\beta\in\F^*\}\rangle\cong\Z_2^3\times(\F^*)^2.
 $$
 It is clear that $\Q_4$ is a MAD-group of $\Int\eseis$, since $\Q_4$ lives in $\cent_{\aut \mathfrak{e}_6}(\mathcal{P}_3)= \mathcal{P}_3\times \Psi'(\SL(3))$, and the torus
$\span{p_{\alpha,\beta}\mid \alpha,\beta\in\F^*}\cong(\F^*)^2$ is   a maximal torus of $\SL(3)$ (and of $\PSL(3)=\Int\frak{a}_2$).

We can compute the type of the induced  fine grading on $\mathcal{M}$ by taking into consideration the $\Z_2^3$-grading induced by $\mathcal{P}_3$.
The fixed component $L_e=\fix \mathcal{P}_3$ is a Lie subalgebra of type $\frak{a}_2$, and the other $7$ components are, all of them,
$L_e$-modules isomorphic to the adjoint module summed with two trivial one-dimensional modules.
Thus the grading has one component of dimension $8$ and seven of dimension $10$.
From  here  it is easy to conclude that all the involved order two automorphisms are of type $2A$, that is, $\mathcal{P}_3$ is of type $A^7$.
 If we consider now the $\Z^2$-grading on $\mathcal{M}$ produced by  $\span{S_{\alpha,\beta} }$, it produces the  root decomposition on
 the identity component $\frak{a}_2$, which is of type $(6,1)$. And, on each of the homogeneous components it produces the weight decomposition,  the
 part fixed by the two-dimensional torus is the piece of dimension $2$ jointly with the two trivial submodules, so of dimension $4$ and such
 component is broken in $(6,0,0,1)$. Thus the grading on $\eseis$
  induced by $\Q_4$ is of type $(6,1,0,0)+7(6,0,0,1)=(48,1,0,7)$.


\subsection{A $\Z^6$-grading}\label{sub toromaximal}

Take as  $\Q_5$ a maximal torus of $\aut\eseis$, which induces a distinguished fine (group) grading, the Cartan-grading or the root decomposition,
which is a $\Z^6$-grading of type $(72,0,0,0,0,1)$, with fixed component a Cartan subalgebra and all the remaining components the corresponding one-dimensional root spaces.


\subsection{All the inner fine gradings}

As a corollary of Proposition~\ref{pr_contienepgrupoeltal}, which will be proved in the technical section, we  obtain one of the main results of this paper:

\begin{te}
The MAD-groups of $\aut\eseis$ contained in $\Int\eseis$ are $\Q_i$ for $i=1,\dots,5$.
\end{te}

\textbf{Proof.}
  If $A$ is a MAD-group different from a maximal torus (that is, if $A$ is not conjugated to $\Q_5$), it is non-toral and, according to Proposition~\ref{pr_contienepgrupoeltal},
$A$ contains a non-toral  subgroup
$V\leq\Int \frak{e}_6 $ isomorphic to either $V_2^3$ or   $V_3^{2a}$ or   $V_3^{2b}$, with the notations in Theorem~\ref{teo_E6-data}.
By Lemma~\ref{le_soloelP1}, Lemma~\ref{le_elP1P2}, Lemma~\ref{le_elP3} and Lemma~\ref{le_elP4}, $V$ must be conjugated to either
$\mathcal{P}_1$ or  $\mathcal{P}_2(\cong \mathcal{P}_4)$ or $\mathcal{P}_3$, and we can assume that $V$ is one of them.

\begin{itemize}
\item If $V = \mathcal{P}_3$, then $A\subset\cent(\mathcal{P}_3)= \mathcal{P}_3\times\Psi'(\SL(3))$, and the problem
reduces to calculate MAD-groups of $\Psi'(\SL(3))\cong\PSL(3)=\Int \frak{a}_2$.
There are four fine gradings on the algebra $\sll(3)$, with grading groups
$
\Z^2,\,\Z\times\Z_2,\,\Z_2^3,\,\Z_3^2;
$
that is,  there are four    MAD-groups of $\aut(\sll(3))\cong\PSL(3)\rtimes\Z_2$, but only two of them are inner,
produced by quasitori of $\PSL(3)$, namely,  $\Z_3^2$ and a two-dimensional torus.
(This result can be concluded from \cite{tipoA}, but the gradings are explicitly computed in \cite{sobreA2}.)
Hence the only possibilities for $A$ are  $\Q_3=\mathcal{P}_3\times \mathcal{P}_4$ and
$\Q_4=\mathcal{P}_3\times \span{S_{\alpha,\beta}}$.

\item If $V= \mathcal{P}_1$, then $A\subset\cent(\mathcal{P}_1)=\mathcal{P}_1\times\Psi(\SL(3))$, and again it
reduces to calculate MAD-groups of $\Psi(\SL(3))\cong\PSL(3)$, which are either the torus $\span{p_{\alpha,\beta}}\cong(\F^*)^2$ or a non-toral $ \mathbb{Z}_3^2$ (just $\span{b,c}$).
 In the first case we  obtain $\Q_2=\mathcal{P}_1\times \span{T_{\alpha,\beta}}$, 
and, in the second one, precisely $\Q_1$.

\item If $V= \mathcal{P}_2$, then $A\subset\cent(\mathcal{P}_2)= \mathcal{P}_2\times G_2$, and it
reduces to calculate the MAD-groups of $G_2=\aut \frak{g}_2$, which are
known to be the two-dimensional torus and $\Z_2^3$ (\cite{g2}). In the first case  we    again get
$\mathcal{P}_2\times T_2=\Q_2$ (see Lemma~\ref{le_elP1P2}).
In  the second case  (we can take $V=\mathcal{P}_4$),   $\Q_3=\mathcal{P}_3\times \mathcal{P}_4$ appears again.
\end{itemize}
This completes the
classification.
\hfill  $\square$ \smallskip

\section{Technical proofs for the inner gradings.}\label{sec_demostraciones}

The aim of this section is to prove  Proposition~\ref{pr_contienepgrupoeltal}, what we will do by means of computational tools inspired in \cite{f4}.
We use these computational techniques much less than there in \cite{f4}, because our computations will only make
use of the Weyl group of $\eseis$ and not of the explicit construction of the normalizer of a maximal torus.
In particular such computations are done
very easily with any mathematical software.
For most of the computations, not even we   need all the elements in the Weyl group but  it is enough to have representatives of its orbits, which can be found, for instance,  in
the web page
\cite{Atlas}. The auxiliar quasitori that we will need for our argumentation, $\tor{f}$, $\sor{f}$ and so on, are easily computed by hand.

\subsection{  Weyl group}\label{subsec_Weyl}

In order to describe the abstract Weyl group of $\eseis$, we must
begin by fixing a basis $\Delta=\{\a_i\mid i=1,\ldots,6\}$ of a
root system of $\eseis$. Its Dynkin diagram is

 
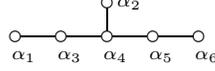
\begin{center}{\vbox{\begin{picture}(25,5)(4,-0.5)
\put(5,0){\circle{1}} \put(9,0){\circle{1}} \put(13,0){\circle{1}}
\put(17,0){\circle{1}} \put(21,0){\circle{1}}
\put(13,3){\circle{1}}
\put(5.5,0){\line(1,0){3}}
\put(9.5,0){\line(1,0){3}} \put(13.5,0){\line(1,0){3}}
\put(17.5,0){\line(1,0){3}} \put(13,0.5){\line(0,1){2}}
\put(4.7,-2){$\scriptstyle \alpha_1$} \put(8.7,-2){$\scriptstyle \alpha_3$}
\put(12.7,-2){$\scriptstyle \alpha_4$} \put(16.7,-2){$\scriptstyle \alpha_5$}
\put(20.7,-2){$\scriptstyle \alpha_6$} \put(13.9,2.6){$\scriptstyle \alpha_2$}
\end{picture}}}\end{center}

\noindent and its Cartan matrix is
\begin{equation}\label{matrizdeCartan}
\small\begin{pmatrix} 2 & 0&-1 & 0 & 0&0\cr 0&2&0&-1&0&0\cr
 -1 & 0&2 & -1&0&0\cr 0 & -1 & -1 & 2&-1&0\cr 0&0&0&-1&2&-1\cr
 0&0&0&0&-1&2
\end{pmatrix}.
\end{equation}

 Take the euclidean space
 $E=\sum_{i=1}^6\R\a_i$ with the inner product $(\ ,\ )$ described for instance in
\cite[Section 8]{Humphreysalg}.
The Weyl group of $\eseis$ is the subgroup $\W$ of $\GL(E)$
generated by the   reflections $s_i$ with $i=1,\dots,6$,
given by $s_i(x):=x-\span{x,\a_i}\a_i$,
 for  $\span{x,y}:=\frac{2(x,y)}{(y,y)}$ (so that the Cartan integers $\span{\a_i,\a_j}$ are just the entries
of the Cartan matrix).
 Identify  $\GL(E)$ to
$\GL(6,\R)$ by means of the matrices relative to the $\R$-basis
$\Delta$.

 We shall consider $\W\subset \GL(6,\R)$
lexicographically ordered. That is: first, for any two different couples
$(i,j), (k,l)$ such that $i,j,k,l\in\{1,\dots,6\}$, we define
$(i,j)<(k,l)$ if and only if either $i<k$ or $i=k$ and $j<l$; and
second, for any two different matrices $\si=(\si_{ij})$,
$\si'=(\si'_{ij})$ in $\W$,  we state  $\si<\si'$ if and only if
$\si_{ij}<\si'_{ij}$  
where $(i,j)$ is the least element (with the
previous order in the couples)  such that $\si_{ij}\ne\si'_{ij}$.
One possible way to compute the Weyl group with this particular
enumeration is provided by the following code implemented with
{\sl Mathematica}:
\smallskip

{\parindent=3cm\tt
W=Table[$s_i$,\{i,6\}];

a[L\_,x\_]:=Union[L,

\hskip 2cm Table[L[[i]].x,\{i,Length[L]\}],

\hskip 2cm Table[x.L[[i]],\{i,Length[L]\}]]

Do[W=a[W,$s_i$],\{i,6\}]\hskip 1cm \textrm{(6  times repeated)} }
\smallskip

We get a list of $51840=2^63^45$ elements in the table {\tt W} which
is nothing but the Weyl group $\W$ of $\eseis$. We are denoting by
$\si_i$ the $i$-th element of $\W$ with the lexicographical  order.

Recall from  \cite[pg.~75]{Humphreysalg} that any
$\sigma\in\W$ 
can be extended to an automorphism $\widetilde\sigma\in\Int\eseis $.  According to that theorem,
if $L=\eseis=\frak{h}\oplus(\oplus_{\a\in\Phi}L_\a)$ is the root decomposition relative to a Cartan subalgebra $\frak{h}$,
for any
choice $x_{\a_i}\in L_{\a_i}\setminus\{0\}$ and
$x'_{\sigma(\a_i)}\in L_{\sigma(\a_i)}\setminus\{0\}$ for
$i=1,\dots,6$, there is only one $\widetilde\sigma\in\aut\eseis $ such that
$\widetilde\sigma(t_{\a_i})=t_{\sigma(\a_i)}$ and
$\widetilde\sigma(x_{\a_i})=x'_{\sigma(\a_i)}$ for every $i=1,\dots,6$, where $t_\a$ is the only element in $\frak{h}$
such that $k(t_\a,\cdot )=\a$, for $k$ the Killing form.
For having fixed a precise family of extensions, consider all the choices $x_{\a_i},x'_{\sigma(\a_i)}$  in the base $B$ chosen as in Proposition~\ref{pr_pasoacposarbitrarios}. Thus we have extensions $\{\tilde\sigma_i\mid i\le51840\}\subset\Int\eseis$.
We are not going to make use of precise descriptions of these extensions.

Denote by $t_{x,y,z,u,v,w}$ the only automorphism of $\eseis$ which acts scalarly on $\frak{h},L_{\a_1},\dots,L_{\a_6},$ with eigenvalues
$\{1,x,y,z,u,v,w\}$
respectively. Take $\T=\{t_{x,y,z,u,v,w}\mid x,y,z,u,v,w\in\F^*\}$, which is a maximal torus of $\mathfrak{G}_0$.
Any other extension of $\sigma\in\W$ as in the above paragraph is equal to $\tilde \sigma t$ for some $t\in\T$.
Recall that
the Weyl group acts in this torus by means of
$\W\times\T\to\T$ given  by $\sigma\cdot t:=\tilde\si
t\tilde\si^{-1}$ for $\si\in\W$ and $t\in\T$. This action does not depend on the choice of the extension $\tilde\si$.
Thus $\si\cdot t_{x,y,z,u,v,w}=t_{x',y',z',u',v',w'}$
for
\begin{eqnarray}\label{ec_acciondelgrupodeWEyl}
 x'= & x^{a_{11}}y^{a_{12}}z^{a_{13}}u^{a_{14}}v^{a_{15}}w^{a_{16}}\cr
 y'= & x^{a_{21}}y^{a_{22}}z^{a_{23}}u^{a_{24}}v^{a_{25}}w^{a_{26}}\cr
 z'= & x^{a_{31}}y^{a_{32}}z^{a_{33}}u^{a_{34}}v^{a_{35}}w^{a_{36}}\cr
 u'= & x^{a_{41}}y^{a_{42}}z^{a_{43}}u^{a_{44}}v^{a_{45}}w^{a_{46}}\cr
 v'= & x^{a_{51}}y^{a_{52}}z^{a_{53}}u^{a_{54}}v^{a_{55}}w^{a_{56}}\cr
 w'= & x^{a_{61}}y^{a_{62}}z^{a_{63}}u^{a_{64}}v^{a_{65}}w^{a_{6
 6}}.
 \end{eqnarray}

 Take also $\No(\T)=\{f\in\aut\eseis\mid ftf^{-1}\in\T\,\forall t\in\T\}$ the normalizer of the torus and
 $\No_0(\T):=\No(\T)\cap\Int\eseis$,
 and consider the projection
  $\pi\colon\No_0(\T)\to \No_0(\T)/\T\cong\mathcal{W}$.
As   there does not exist a section of $\pi$ (see \cite{raros}),
we have an injection $\iota\colon\mathcal{W}\to\No_0(\T)$ ($\sigma\mapsto\tilde\sigma$)
but $\iota$ is not a group homomorphism. 

Let us consider our previous notations of Equation~(\ref{eq_eltsuperf}) in these new terms.
If $\eta\in\W$ and $s\in\T$, we denote by $\tor{\eta}:=\tor{\tilde\eta}=\{t\in\T\mid \eta\cdot t=t\}$, and by
 $Q(\eta,s):=Q(\tilde{\eta}s)$, that is, the quasitorus generated by $\tilde\eta s$ and $\tor{\eta}$.
We compute $\tor{\eta}$ for each representative $\eta$ of some orbit.
%
%
The $51840$ elements of the Weyl group $\W$ of $\eseis$ are distributed
in $25$ orbits (=conjugacy classes), whose representatives could be found  by using any matrix multiplication software.
We extracted such representatives from \cite{Atlas}, and we identify them to some elements in our ordered list
$\mathcal{W}$ in order to make possible to do computations with them.

\begin{center}
\begin{tabular}{|c|c|c|c|c|}
\hline Order& Representative & Orbit & Stabilizer &
Isomorphic to\cr \hline\hline
 $1$ &$40843$ & $1$ & $xyzuvw\ne0$ & $(\F^*)^6$\cr
  $2 $ &$19 $ & $270 $ & $w=u,v=\frac1{u^2xyz} $ & $(\F^*)^4 $\cr
  $ 2$ &$21 $ & $540 $ & $v=1, w=u, z=\frac1{u^2xy} $ & $(\F^*)^3 $\cr
  $ 2$ &$96 $ & $45 $ & $ z=\frac1{xy},u=w,u^2=v^2=1 $ & $(\F^*)^2 \times\Z_2^2$\cr
  $2 $ &$11323 $ & $36 $ & $x=1  $ & $(\F^*)^5 $\cr
  $4 $ &$2 $ & $3240 $ & $z=xy, u=\frac1{z},v=1,w=u $ & $(\F^*)^2 $\cr
  $ 4$ &$20 $ & $1620 $ & $u=w=1,v=\frac1{xyz} $ & $(\F^*)^3 $\cr
  $ 4$ &$75 $ & $540 $ & $u^2=v^2=1,y= \frac{v}{x},z=v,w=u $ & $\F^{* }\times\Z_2^2$\cr
  $4 $ &$ 140$ & $540 $ & $ x=u=w=1,v=\frac1{yz}$ & $(\F^*)^2 $\cr
  $8 $ &$ 1$ & $6480 $ & $u=v=z=w=1,y=\frac1{x} $ & $\F^{* }$\cr
  $3 $ &$292 $ & $480 $ & $z=y,v=x,w=u^3x^2y^2,(xyu)^3=1 $ & $(\F^*)^2 \times\Z_3$\cr
  $3 $ &$3819 $ & $80 $ & $x^3=y^3=z^3=1,u=y,v=y^2z,w=xy^2$ & $\Z_3^3$\cr
  $ 3$ &$ 4079$ & $240 $ & $x=1,w=\frac1{u^2v^2yz}  $ & $(\F^*)^4 $\cr
  $6 $ &$5 $ & $1440 $ & $u=\frac1{xy},v=\frac{xy}z,w=u $ & $(\F^*)^ 3$\cr
  $6 $ &$15 $ & $2160 $ & $y=\frac1{xu^2},z=v=1,w=u $ & $(\F^*)^2 $\cr
  $6 $ &$22 $ & $1440 $ & $u=v=w=1,z=\frac1{xy} $ & $(\F^*)^ 2$\cr
  $6 $ &$122 $ & $4320 $ & $ x^3=1,z=y,u=\frac1{y},v=x,w=ux^2$ & $\F^{* }\times\Z_3$\cr
  $6 $ &$124 $ & $720 $ & $x^3=1,y=z=u=1,v=x,w=x^2 $ & $\Z_3$\cr
  $6 $ &$195 $ & $1440 $ & $x^3=1,u^2=v^2=1,y=z=v,w=ux $ & $\Z_3\times\Z_2^2$\cr
  $6 $ &$435 $ & $1440 $ & $z=y=u,v=x^3y^5,w=xy,(xy^2)^3=1 $ & $\F^{* }\times\Z_3$\cr
  $9 $ &$121 $ & $5760 $ & $x^3=1,y=z=x^2,u=v=x,w=1 $ & $\Z_3$\cr
  $12 $ &$4 $ & $4320 $ & $z=u=v=t=1,y=\frac1{x} $ & $\F^{* }$\cr
  $12 $ &$218 $ & $4320 $ & $x^3=1,y=z=u=1,v=x,w=x^2 $ & $\Z_3$\cr
  $ 5$ &$3 $ & $5184 $ & $z=1,u=\frac1{xy},v=xy,w=u $ & $(\F^*)^2 $\cr
  $10 $ &$135 $ & $ 5184$ & $x=z=v=1,w=u,y=\frac1{u^2} $ & $\F^{* }$\cr
 \hline
\end{tabular}
\end{center}\begin{center}\begin{equation}\text{Table of representatives of $\W$}\label{tablainternos} \end{equation}\end{center}\vskip0.5cm

The second row of this table means the following: the element $\sigma_{19}$ has order $2$, its conjugacy class has $270$ elements,
and $\tor{\sigma_{19}}=\{t_{x,y,z,u,v,w}\in\T\mid w=u,v=\frac1{u^2xyz} \}\cong\F^{*4 }$. Such element is the 19th in our list $\W$, and
it appears explicitly in \cite{Atlas}.
For short we also denote   by $\tor{i}=\tor{\si_{i}}$ and   $Q(i,s)=Q( {\si_i},s)=Q(\widetilde{\si_i}s)$.

\begin{re}\label{re_losindicesnotorales}
According to Remark~\ref{re_cuandoQ(f)toralseguro}, if $Q( \eta,s)$ is non-toral and if
\begin{itemize}
\item $\eta$ has order $2$, then $\eta$ is conjugated to $\si_{96}$;
\item $\eta$ has order $4$, then $\eta$ is conjugated to $\si_{75}$;
\item $\eta$ has order $3$, then $\eta$ is conjugated to either $\si_{292}$ or $\si_{3819}$.
\end{itemize}
\end{re}

\noindent And, following again Remark~\ref{re_cuandoQ(f)toralseguro}, it is not possible that $Q( \eta,s)$ is non-toral
if $\eta$ has order $5$, so:

\begin{lm} \label{le_nohaydel5}
There is no non-toral $5$-group of $\aut\eseis$.
\end{lm}

\textbf{Proof.}
Suppose that there is  $Q\le\Int\eseis$ a non-toral $5$-group. Let $Q'$ be a non-toral minimal quasitorus contained in $Q$ (any $Q''\subsetneq Q'$ is toral). We can assume that $Q'\subset\No(\T)$ and that
$Q'\cap\T$ is maximal toral in $Q'$. Hence there are some  element $\eta\in\mathcal{W}$ of order a power of five and some $s\in\T$
such that $Q'\subset Q( {\eta},s)$. The contradiction appears since the last quasitorus is toral, and $Q'$ is non-toral.
\hfill  $\square$ \smallskip

\noindent (This lemma is  also consequence of the same result for the complex field, proved in \cite{Griess}, jointly with Proposition~\ref{pr_pasoacposarbitrarios}).

\begin{re}\label{re_proyecciones}
Again by looking at Table~\ref{tablainternos}, we observe the following useful fact:
For $f\in\No(\T)$, fix any subgroup $\hor{f}$ satisfying
the conditions of Lemma~\ref{le_comoeselTsuperf} (that is, $\tor{f}=\sor{f}\times\hor{f} $). Then denote by $\wp_f\colon\tor{f}\to\hor{f}$ the projection.
Now, if $Q$ is a non-toral subquasitorus of $Q(f)$ such that $\pi(f)$ is   in the orbit of neither $\si_{3819}$ nor $\si_{195}$,
then $\wp_f(Q\cap\T)=\hor{f}$.  

Indeed, the quasitorus $\span{f} \times\wp_f(Q\cap\T)$ is toral taking into account   that there are no non-toral $2$-groups with less than $3$ factors. By Lemma~\ref{le_toroporalgotoral},   $\span{f}\times\sor{f}\times\wp_f(Q\cap\T)$ is also toral, as well as the quasitorus $Q$ which is contained in it.
\end{re}

\subsection{MAD-groups of $\Int\eseis$ in computational terms}

\begin{pr}\label{pr_MADsenterminoscomputacionales}
The MAD-groups described in Section~\ref{sec descripcionesinternos} can be described in these terms as follows.
\begin{itemize}
\item  The quasitorus $Q(3819,\id)\cong\Z_3^4$ is conjugated to $\Q_{1}$.
\item  The quasitorus $Q(292,\id)\cong(\F^*)^2\times\Z_3^2$ is conjugated to $\Q_{2}$.
\item  The quasitorus $Q(195,\id)\cong\Z_3^2\times\Z_2^3$ is conjugated to $\Q_{3}$.
\item  The quasitorus $Q(96,\id)\cong(\F^*)^2\times\Z_2^3$ is conjugated to $\Q_{4}$.
\item The quasitorus $Q(\id)=\T$ is the maximal torus conjugated to $\Q_5$.
\end{itemize}
\end{pr}

\textbf{Proof.}
We use the notations in Section~\ref{sec descripcionesinternos}.
\begin{itemize}
\item[$\Q_1$)] As  $\span{F_2,F_3,F_4}\cong\Z_3^3$ is toral as in Remark~\ref{re_unZ3cubotoral} but  $\Q_1$ is a non-toral quasitorus, there is a maximal torus
(we can conjugate to choose such torus equal to $\T$) such that $\span{F_2,F_3,F_4}\subset \T$ and $F_1\in\No(\T)$ by Lemma~\ref{le_puedoajustarlaparatedeltoro}.  Moreover, $F_1\in\No_0(\T)$ because $F_1\in\Int\eseis$.
As $F_1$ has order $3$, also $\pi(F_1)$ has order $3$ (not $1$ because $\Q_1$ would be  toral), so that we can conjugate without changing the torus to get
$\pi(F_1)=\si_{j}$ with $j\in\{ {292}, {3819}, {4079}\}$ (recall that if $\sigma$ and $\tau$ are conjugate in $\W$, then some element in $\pi^{-1}(\sigma)$ is conjugated to some element in $\pi^{-1}(\tau)$ in $\No(\T)$). The   possibility   $j=4079 $ is ruled out as in Remark~\ref{re_losindicesnotorales}.
Moreover, $\Q_1=\span{F_2,F_3,F_4,F_1}\subset Q(F_1)$ (since $\span{F_2,F_3,F_4}\subset\tor{F_1}$) and $\mathcal{Q}_1$ is a MAD-group, so that $\mathcal{Q}_1=Q(F_1)$. Thus $j\ne292$, since $Q(292)\cong(\F^*)^2\times\Z_3^2$.
Therefore $F_1=\widetilde\si_{3819}s$ for some $s\in\T$ 
and $\Q_1=Q(3819,s)$, which is conjugated to $Q(3819,\id)$ as in Remark~\ref{re_paramover}, since $\tor{3819}$ is finite.

\item[$\Q_2$)] $F_1$ is a toral element because it is an inner automorphism, so that $\langle F_1, T_{\alpha,\beta}\mid \alpha,\beta\in\F^*\rangle$
is also toral  by Lemma~\ref{le_toroporalgotoral} and, 
as before, we can assume that $\langle F_1, T_{\alpha,\beta}\mid \alpha,\beta\in\F^*\rangle\subset\T$ and that $F_2=\widetilde\si_{j}s\in\No_0(\T)$ for some  $j\in\{ {292}, {3819}
\}$ and $s\in\T$. It is clear that $j=292$, because a two-dimensional torus is not contained in $\tor{3819}\cong\Z_3^3$,
 so that
$\Q_2=Q(292,s)$, which is conjugated to $Q(292,\id)$ since $\tor{292}\cap\jor{292}=\{t_{x,u,u,u,x,x^2u^2}\mid x^3=u^3=1\}\cong\Z_3^2$ is finite.

\item[$\Q_3$)] Now note that $\span{H_2,G_1}$ is toral (by Lemma~\ref{le_algunodelosfactoresnotoral}, since it is isomorphic to $\Z_2^2\times\Z_3$), so that we can assume that  $\span{H_2,G_1}\subset\T$ and $H_1\in \No_0(\T)$.
    Hence
     $\mathcal{Q}_3=\span{H_2,G_1,H_1}\subset Q(H_1)$, and, as $\Q_3$ is maximal, then $\Q_3=Q(H_1)$. That fact forces
     $\Q_3$ to be $Q(195,s)$ for some $s\in\T$,  which is conjugated to $Q(195,\id)$ as in Remark~\ref{re_paramover}, since $\tor{195}$ is finite.

\item[$\Q_4$)] Finally, $\{G_1,G_2\}$ is toral (two $2$-factors are always toral)
and we can assume that $\span{G_1,G_2,S_{\alpha,\beta}}\subset\T$
and that $G_3\in\No_0(\T)$ projects in some element of the orbit of $\si_{96}$, by Remark~\ref{re_losindicesnotorales}. Thus $\Q_4$ is conjugated to some $Q(96,s)$ and hence to $Q(96,\id)$ because $\tor{96}\cap\jor{96}=\{t_{ yz,y,z,u,v,u}\mid y^2=z^2=u^2=v^2=1\}\cong\Z_2^4$ is finite.
\end{itemize} \hfill  $\square$ \smallskip


\subsection{Order of the extensions}\label{subsec_ordenesdelasextensiones}

Note that the order of $f\in\No_0(\T)$ is multiple of the order of $\pi(f)\in\W$, and that both numbers could not coincide. For instance, in
\cite[Remark~1]{f4} it is observed that $\widetilde\sigma_3\in\aut\f4$ has order $8$, while its projection $\sigma_3$ on the Weyl group of $\f4$ has order $4$. Moreover, any element in $\pi^{-1}(\sigma_3)$ has also order $8$, as in Remark~\ref{re_sobreordendelaextension}. None of them has the same order than its projection.

\begin{re}\label{re_afinandoelSsuperf}
 We can extend the results in Lemma~\ref{le_comoeselTsuperf}  a little bit. If $f\in\No_0(\T)$ verifies that its projection $\pi(f)$ has order $r'$  (in general, $r'$ divides the order of $f$),
 then $\hor{f}\subset\{t\in \T\mid t^{r'}=1_G\}$ and the torus
 $\sor{f}=\{(tf)^{r'}(f)^{-r'}\mid t\in\T\}=\{\pi_{i=0}^{r'-1}f^i\cdot t\mid t\in\T\}$.
 \end{re}

\noindent In the previous example about $\f4$, what happens is $(\widetilde\sigma_3s)^8\in\sor{3}=\id$ but $\id\ne(\widetilde\sigma_3s)^4\in\tor{3}$ for every $s$ in the corresponding torus.
\smallskip

It is difficult in general to know the order of the extension of a concrete element of the Weyl group
 by applying only
the isomorphism theorem in
\cite[Section~14.2]{Humphreysalg},
 but in this case we have extra-information extracted from Proposition~\ref{pr_MADsenterminoscomputacionales}.
 Using, also, the results in Subsection~\ref{subsec_detallestecnicos},
  we conclude that there are elements $s_{292}\in\tor{292}$ and $s_{96}\in\tor{96}$ such that
 the order of     $\widetilde\si_{96}s_{96}$ is $2$, the order of the extensions $\widetilde\si_{292}s_{292}$ and $\widetilde\si_{3819}$ are $3$ and
the order of the extension $\widetilde\si_{195}$ is $6$, and, up to conjugation,
$\mathcal{P}_2\cong\mathcal{P}_4$ is conjugated to the set of order three elements in $Q(195)$
and
$$
\mathcal{P}_1\cong \span{ \widetilde\si_{292}s_{292}}\times\hor{292},\qquad
\mathcal{P}_3\cong \span{ \widetilde\si_{96}s_{96}}\times\hor{96},
$$
for $\hor{292}=\span{t_{1,1,1,\omega,11}}$ and $\hor{96}=\{t_{1,1,1,u,v,u}\mid u^2=v^2=1\}$, that we fix for the rest of this section.

We obtain the same conclusions by reading the Appendix, in which natural representatives of those such extensions are constructed.

For technical purposes, note that, according to Remark~\ref{re_afinandoelSsuperf},
$(\widetilde\si_{292}t)^3\in\sor{292}$ and   $(\widetilde\si_{96}t)^2\in\sor{96}$ for any $t\in\T$, since there exist extensions
of $\sigma_{96}$ and of $\sigma_{292}$ of orders $2$ and $3$ respectively.

\subsection{On the elementary $p$-group.}\label{subsec_contienepgrupoeltal}

\ \smallskip

\textbf{Proof of Proposition~\ref{pr_contienepgrupoeltal}.}
Now suppose that $Q$ is in the conditions of the Proposition~\ref{pr_contienepgrupoeltal}.
Then $Q=P\times \prod_{i}P_{p_i}$, where $P$ is a torus and $P_{p_i}$ are   $p_i$-groups, at least one of them non-toral (by Lemma~\ref{le_algunodelosfactoresnotoral}) for some $p_i\in\{2,3\}$.
We want to prove that either $P_2$ contains a non-toral $\Z_2^3$-subgroup or $P_3$ contains a non-toral $\Z_3^2$-subgroup.
By Lemma~\ref{le_puedoajustarlaparatedeltoro}, we can assume that $Q$ is contained in $\No_0(\T)$ for some maximal torus $\T$ in such a way that $Q\cap\T$ is maximal-toral in $Q$:
that is, if $Q\cap\T\subset  Q'\subset Q$ with $Q'$  toral, then $Q\cap\T=Q'$.
Observe that $P_{p_i}\cap\T$ is maximal-toral in $P_{p_i}$:
otherwise certain $h\in P_{p_i}\setminus\T$ would verify that $\span{P_{p_i}\cap\T,h}$ would be toral, and then $\span{Q\cap\T,h}$ would be toral too, by \cite[Corollary~1]{f4spin}.
In particular $P\times P_5\times P_7\times...\subset\T$, $P_2=P_2\cap\T\cdot\span{f_1,\dots,f_n}$ with each $f_j\in\No_0(\T)$ of order   a power of $2$ and $P_3=P_3\cap\T\cdot\span{g_1,\dots,g_m}$ with each $g_j\in\No_0(\T)$ of order   a power of $3$.
Besides, as $Q$ is a MAD-group, then $Q\cap\T=\tor{f_1}\cap\dots\cap\tor{f_n}\cap\tor{g_1}\cap\dots\cap\tor{g_m}$.
Of course, we can take   $f_i\notin \span{Q\cap\T,f_1,\dots,f_{i-1}}$, since otherwise
$\span{Q\cap\T,f_1,\dots,f_{i-1}}=\span{Q\cap\T,f_1,\dots,f_{i}}$. Hence $f_i\notin \span{ \T,f_1,\dots,f_{i-1}}$
and $\pi(f_i)$ does not belong to the group generated by $\{ \pi(f_1),\dots,\pi(f_{i-1})\}$.

As $\span{f_i,P_2\cap\T}$ is a non-toral subquasitorus of $Q(f_i )$,
we know by Remark~\ref{re_losindicesnotorales} that $\pi(f_i)$ (an order 2, 4 or 8 element) belongs to the orbit of either
$\si_{96}$ or $\si_{75}$. The same argument shows that $\pi(g_i)$ (an order 3 or 9 element) belongs to the orbit of either
$\si_{292}$, $\si_{3819}$ or $\si_{121}$. We can discard the  possibilities $75$ and $121$:
\begin{itemize}
\item In the first case, note that $\tor{\sigma_{75}^2}\cong(\F^*)^4$, so that $\span{P_2\cap\T,f_i^2}$ is toral contained in $P_2$, a contradiction with the choice of $P_2\cap\T $ as maximal-toral   in $P_2$.
\item In the second case, we can change the element  $g_i$ of order $9l$  by its conjugated element $\widetilde\si_{121}s$ with $s\in\T$ (we
conjugate by means of an element in the normalizer, so
that $Q\cap\T$ is still contained in $\T$). As $Q\cap\T\subset\tor{g_i}\cong\Z_3$, then $ Q\cap\T=\tor{g_i}$ (since jointly with $g_i$ is non-toral)
and
 $Q\subset \cent_{ \mathfrak{N}_0(\T)}(Q(121,s))$. This centralizer is known to coincide with itself (a direct computation with the computer)
so that by maximality $Q=Q(121,s)$. We have got a contradiction because this set is not a MAD-group:
Inspired in \cite[Proposition~7]{f4}, there is $h\in\aut\eseis$ such that $hg_i h^{-1}\in\T$ and $hth^{-1}\in\No(\T)$
for $t=t_{\omega,\omega^2,\omega^2,\omega,\omega,1}$   the generator of $ \tor{121}$.
Thus $hQh^{-1}\subset Q(k,s')$ for certain $\si_k$ of order $3$ and $s'\in\T$,
but in no case the quasitorus  $Q(k,s')$ is isomorphic to $\Z_3\times\Z_9$ (according to Table~\ref{tablainternos}) and the contained is strict.
\end{itemize}

We have proved in particular that the order of any $\pi(f_i)$ is $2$ and the order of any $\pi(g_i)$ is $3$.

Observe the following technical fact (which saves a lot of computations).
There are $113$ order $2$ elements $\si_{j}$ (including $j=96$) in $\W$ commuting with  $\si_{96}$. Only $13$ of them verifies that
$\si_j$ is in the orbit of  $\si_{96}$. And all these $13$ elements verify that $\si_{j}\si_{96}$ is not in the orbit of $\si_{96}$.
This proves that necessarily $n\le1$.

$\bigstar$
If $m\le1$, then by maximality $Q\in\{Q(f_1),Q(g_1),Q(f_1g_1)\}$. The only possible non-toral
$Q(j,s)$ for some $\si_j$ of order $2,3$ or $6$ are those ones with $j\in\{96,292,3819,122,124,195,435\}$ by Lemma~\ref{le_toroporalgotoral}.
 If $j\in\{96,292,3819,195 \}$, then $Q$ is conjugated to some quasitorus in the list $ \{\Q_i\mid i=1,2,3,4\}$ by Proposition~\ref{pr_MADsenterminoscomputacionales}.
In such a case we are done because all these quasitori contain one of the required elementary $2$-groups or $3$-groups (since $\mathcal{P}_1\subset\Q_1\cap\Q_2$ and $\mathcal{P}_3\subset \Q_3\cap\Q_4$).
If $j\in\{122,435\}$,
then $\sigma_j^3$ is not conjugated to $\sigma_{96}$,
so that $\span{(\widetilde\si_js)^3,\tor{j}}$ is toral, which is
a contradiction with the choice of $Q\cap\T$.
And if $j=124$, although $\sigma_j^3$ is   conjugated to $\sigma_{96}$,
again $\span{(\widetilde\si_js)^3,\tor{j}}$ is toral. To see it, it is enough to note that $\tor{j}\cong  \Z_3$,
and, once we have changed $ (\widetilde\si_js)^3$ with $ \widetilde\si_{96}s'$, the image of the non-toral quasitorus by the projection $\wp_{96}$ should be $\hor{96}$, as in Remark~\ref{re_proyecciones}, but it is identity.

$\bigstar$ Finally suppose that $m\ge2$.
Again we can check that
there are $26$ order $3$ elements $\si_{j}$ (including $\si_{292}$ and $\si_{292}^2$) in $\W$ commuting with  $\si_{292}$, $8$ of them belonging to the orbit of $\si_{3819}$, $12$ in the orbit of $\si_{292}$ and $6$ of them in the orbit of  $\si_{4079}$, the  \lq\lq toral orbit\rq\rq (knowing the orbit is an easy computation, it is enough in this case to compute the stabilizer). But what it is useful is that either $\si_{j}$ or $\si_{j} \si_{292}$ or $\si_{j} \si_{292}^2$ is not conjugated to $ \si_{292}$. Besides $\tor{292}\cap\tor{j}\cong\mathbb{Z}_3^2$ in all such cases when $\si_{j}$ conjugated to $ \si_{3819}$.

Hence we can assume that $\pi(g_1)= \si_{3819}$. In particular $g_1$ has order $3$ and $Q\cap\T\subset \tor{3819}\cong\Z_3^3$.
If $Q\cap\T$ is isomorphic to $\Z_3$, we are done because $ \span{Q\cap\T,g_1}$ is a non-toral elementary $3$-group of rank $2$.
Also $Q\cap\T$ is not isomorphic to $\Z_3^3$,   because  in such a case $ Q\supsetneq Q(g_1)$, but $Q(g_1)$ is a MAD-group.
Hence $Q\cap\T\cong\Z_3^2$. If $g_2$ has order $3$,
then $Q\supset Q'=\span{g_1,g_2}\cdot\tor{g_1}\cap\tor{g_2}\cong\Z_3^4$. According again to Theorem~\ref{teo_E6-data}
and Proposition~\ref{pr_pasoacposarbitrarios},
 $Q'$ is either of type $V_3^{4a}\equiv\Q_1\supset\mathcal{P}_1$ or of type $V_3^{4b}$, which is identified with the order three
 elements in $\Q_2$, so that it also contains $\mathcal{P}_1$.
 Otherwise, $g_2^3\in\T\cap P_3\cong\mathbb{Z}_3^2$, hence $g_2$ has order $9$ (and $\pi(g_2)$ is necessarily conjugated to $\si_{292}$).
In particular there is $t\in\T\cap P_3$ of order $3$ such that $P_3\supset P_3'=\span{g_1,g_2,t}\cong\mathbb{Z}_3^2\times\mathbb{Z}_9$, which is non-toral.
As $g_2$ is an inner automorphism of order $9$, there is $h\in\aut\eseis$ such that $hg_2 h^{-1}\in\T$ and $hg_1 h^{-1},hth^{-1}\in\No_0(\T)$. Thus $\pi(hg_1 h^{-1})$ and $\pi(hth^{-1})$ are conjugated to $\si_{292}$, because both contain some $\mathbb{Z}_9$ in the stabilizer and $P_3'$ is non-toral. We can move again with another automorphism, this time in $\No(\T)$, such that  $\pi(hg_1 h^{-1})=\si_{292}$ and $\pi(hth^{-1})=\si_{j}$. As in the paragraph above there is some $l=0,1,2$ such that $\si_{292}^l\si_j$ is conjugated to either $\si_{3819}$
or $\si_{4079}$. In the first case,
 the order $9$ element
$hg_2 h^{-1}\in\tor{292}\cap\tor{j}=\tor{\si_{292}}\cap\tor{\si_{292}^l\si_j}\cong\mathbb{Z}_3^2$, a contradiction.
In the second case, $\span{g_2,g_1^lt}\cong\mathbb{Z}_3\times \mathbb{Z}_9$ would be toral, so that we could conjugate $P_3'$ to a subgroup of some $Q(k,s)$. But $P_3'$ cannot fill $Q(k,s)$, because according to Table~\ref{tablainternos}, there are no $Q(k,s)$ isomorphic to $\mathbb{Z}_3^2\times \mathbb{Z}_9$. Thus $P_3'\ne Q$ and either $m>2$ or $n>0$.
If $m>2$, then $\span{\pi(g_1),\pi(g_2),\pi(g_3)}$ is just the whole set of elements of order divisor of 3 which commute with $\si_{292}$ (the $26$ elements described before commute among them). But $\tor{g_1}\cap\tor{g_2}\cap\tor{g_3}=\{t_{\omega,\omega,\omega,\omega,\omega,\omega}\}\cong\mathbb{Z}_3$, so $Q\cap\T\ne\mathbb{Z}_3^2$.
And if $n\ne0$, $\pi(f_1)$ commutes with $\span{\si_{292},\si_j}$ for one of the $8$ $\si_j$'s commuting with $\si_{292}$  and in the orbit of $\si_{3819}$ (in fact, the number of candidates to  $\si_j$ can be reduced with considerations about orbits fixing $\si_{292}$  or also taking into account that several of them generate the same groups). For some cases there does not exist an order two element in $\W$  commuting with
 $\span{\si_{292},\si_j}$, and, in all the cases in which there exists such element, $\tor{292}\cap\tor{j}\cap\tor{\pi(f_1)}\cong\mathbb{Z}_3$.

This finishes the proof of Proposition~\ref{pr_contienepgrupoeltal}.
\hfill  $\square$  \medskip

Now we can exhibit an example of a non-toral quasitorus of $\Int\eseis$ which does not contain a non-toral elementary $p$-group
of $\aut\eseis$.
The quasitorus
$$
Q=\span{\widetilde\sigma_{292}s_{292},t_{\xi,\xi,\xi,\xi^4,\xi,\xi^7}   },
$$
isomorphic to $\mathbb{Z}_3\times\mathbb{Z}_9$ ($\xi^3=\omega^2$), satisfies such condition, since  $t_{\xi,\xi,\xi,\xi^4,\xi,\xi^7} \in
\tor{292}\setminus\sor{292}$ but $( t_{\xi,\xi,\xi,\xi^4,\xi,\xi^7} )^3\in\sor{292}$.

\section{Description of the outer gradings}\label{sec descripcionesexternos}

\subsection{Outer automorphisms of finite order}\label{sub_autexternos}

The outer automorphisms of finite order $m$ (necessarily even)
can be obtained from the affine diagram 

\begin{center}{$E_6^{(2)}$:  \hspace{0.6cm}
\begin{picture}(25,5)(4,-0.5)
\put(5,0){\circle{1}} \put(9,0){\circle{1}} \put(13,0){\circle{1}}
\put(17,0){\circle{1}} \put(21,0){\circle{1}}
\put(5.5,0){\line(1,0){3}}
\put(9.5,0){\line(1,0){3}} \put(13.5,0.3){\line(1,0){3}} \put(13.5,-0.1){\line(1,0){3}}
 \put(17.5,0){\line(1,0){3}} \put(14.3,-0.5){$<$}
 \put(4.7, 1.7){$\scriptstyle 1$} \put(8.7, 1.7){$\scriptstyle 2$}
\put(12.7, 1.7){$\scriptstyle 3$} \put(16.7, 1.7){$\scriptstyle 2$}
\put(20.7, 1.7){$\scriptstyle 1$}
\put(4.7,-2){$\scriptstyle \alpha_0$} \put(8.7,-2){$\scriptstyle \alpha_1$}
\put(12.7,-2){$\scriptstyle \alpha_2$} \put(16.7,-2){$\scriptstyle \alpha_3$}
\put(20.7,-2){$\scriptstyle \alpha_4$}
\end{picture}
}\end{center}\medskip

\noindent by assigning weights $\bar p=(p_0,\dots,p_4)$ ($p_i\in\Z_{\ge0}$) such that
$2(p_0+2p_1+3p_2+2p_3+p_4)=m$. Obviously the only possibilities for outer order two
automorphisms (up to conjugation) are obtained when $\bar p=(1,0,0,0,0)$ and $\bar p=(0,0,0,0,1)$,
choices which provide outer automorphisms with fixed subalgebras

\begin{center}{$\f4$ from
\begin{picture}(25,5)(4,-0.5)
 \put(5,0){\circle*{1}}
\put(9,0){\circle{1}} \put(13,0){\circle{1}}
\put(17,0){\circle{1}} \put(21,0){\circle{1}}
\put(5.5,0){\line(1,0){3}}
\put(9.5,0){\line(1,0){3}} \put(13.5,0.3){\line(1,0){3}} \put(13.5,-0.1){\line(1,0){3}}
 \put(17.5,0){\line(1,0){3}} \put(14.3,-0.5){$<$}
\end{picture}
and $\frak{c}_4$ from
\begin{picture}(25,5)(4,-0.5)
 \put(5,0){\circle{1}}
\put(9,0){\circle{1}} \put(13,0){\circle{1}}
\put(17,0){\circle{1}} \put(21,0){\circle*{1}}
\put(5.5,0){\line(1,0){3}}
\put(9.5,0){\line(1,0){3}} \put(13.5,0.3){\line(1,0){3}} \put(13.5,-0.1){\line(1,0){3}}
 \put(17.5,0){\line(1,0){3}} \put(14.3,-0.5){$<$}
\end{picture}.
}\end{center}\smallskip
We will say that an order two outer automorphism $F\in 2C$ (respectively $F\in 2D$) if its fixed subalgebra is isomorphic
to $\f4$ (respectively to $\frak{c}_4$) (\cite{Griess} does not use these automorphisms,
 so that the   notation which we have taken here is simply compatible with Subsection~\ref{subsec_autinternos}).

Note that there are $3$ conjugacy classes of outer automorphisms of order $4$, corresponding to $(1,0,0,0,1)$, $( 0,1,0,0,0)$ and $( 0,0,0,1,0)$.
The fixed subalgebras are of types $\frak{c}_3+ Z$, $\frak{b}_3\oplus\frak{a}_1$ and $\frak{a}_3\oplus\frak{a}_1$ respectively, and hence of dimensions $22$, $24$ and $18$, respectively. This implies that the conjugacy class of  an outer order 4 automorphism is distinguished only by the dimension of the fixed subalgebra.
An  element which will be useful for us is $\Upsilon_1$ any (necessarily outer) automorphism of order $4$ whose fixed  subalgebra is of type $\frak{a}_3\oplus\frak{a}_1$.

 Recall that $\mathfrak{G}=\aut\eseis=\Int\eseis\cup F\Int\eseis$ for any $F\in \mathfrak{G}\setminus \mathfrak{G}_0$.
In this section we study the maximal quasitori of $\aut\eseis$ not contained in the identity component  $\mathfrak{G}_0 =\Int\eseis$. First we describe those ones containing some automorphism of the class $2C$ and next those
ones containing some automorphism of the class $2D$ (some of them coincide).
Afterwards we consider the case when the MAD-group of  $ \mathfrak{G}$ does not contain any outer order two automorphism, although is not contained in
$ \mathfrak{G}_0$. In this case we will prove that the MAD-group   contains necessarily an automorphism conjugated to $\Upsilon_1$.

\subsection{Gradings based on a $\frak{f}_4$-model.}\label{sub_losquevienendef4}

Let $\mathcal{J}$ be the Albert algebra, and $\der\mathcal{J}$ its derivation algebra, which is  a Lie algebra of type $\frak{f}_4$.
Denote by $\mathcal{ J}_0$ the set of zero trace elements of the Albert algebra.
Take $\mathcal{N }:=\mathcal{ J}_0\oplus\der\mathcal{ J}$ with the product given by
\begin{itemize}
\item the restriction of the bracket to $\der\mathcal{ J}$ is the usual bracket;
\item if $d\in\der\mathcal{ J}$ and $x\in \mathcal{ J}_0$, take $[d,x]:=d(x)\in \mathcal{ J}$;
\item if  $x,y\in  \mathcal{ J}_0$, take $[x,y]:=[R_x,R_y]\in\der\mathcal{ J}$, where $R_x$ denotes the multiplication operator in $\mathcal{ J}$.
\end{itemize}
It is well-known that $\mathcal{N}$ is a  Lie algebra   of type $\eseis$.
Consider $G_4$ the order two automorphism producing the grading $\mathcal{ N}_{\bar0}:=\der\mathcal{ J}$ and $  \mathcal{ N}_{\bar1}:=\mathcal{ J}_0$.

It is also  well-known that every automorphism of the Albert algebra can be extended to an automorphism of $\eseis$. Namely, if $f\in\aut\mathcal{ J}$, take $f^{\bullet}\in\aut\mathcal{ N}$  given by $f^{\bullet}(d)=\Ad f(d):=fdf^{-1}$ if $d\in\der\mathcal{ J}$, and $f^{\bullet}(x)=f(x)$ if   $x\in\mathcal{ J}_0$.
Moreover, $\cent_{\aut\mathfrak{e}_6}(G_4)=\{f^{\bullet},f^{\bullet} G_4\mid f\in\aut\mathcal{ J}\}\cong\aut\mathcal{ J}\times \Z_2$.
Indeed, if $\varphi\in\cent_{\aut \mathfrak{e}_6}(G_4)$, then
$\varphi(\mathcal{ N}_{\bar0})\subset \mathcal{ N}_{\bar0}$,
so that $\varphi\vert_{\mathcal{ N}_{\bar0}}\in\aut\der\mathcal{J}=\Ad(\aut\mathcal{J})$
and there is $f\in \aut\mathcal{J}$ such that $\varphi\vert_{\mathcal{ N}_{\bar0}}=\Ad f=f^{\bullet}\vert_{\mathcal{ N}_{\bar0}}$.
Now $\varphi(f^\bullet)^{-1}\vert_{\mathcal{ N}_{\bar1}}\in\hom_{\mathcal{ N}_{\bar0}}(\mathcal{ N}_{\bar1},\mathcal{ N}_{\bar1})$
so 
there is certain $\alpha\in\F^*$ such that $\varphi(f^\bullet)^{-1}\vert_{\mathcal{ N}_{\bar1}}=\alpha\id$ by Schur's Lemma.
As $[\mathcal{N}_{\bar1},\mathcal{N}_{\bar1}]= \mathcal{N}_{\bar0}$, then $\alpha^2=1$ and $\varphi(f^\bullet)^{-1}\in\{\id,G_4\}$.

Hence any MAD-group of $\aut\mathcal{ N}$ containing $G_4$ is the direct product of a MAD-group of $\aut\mathcal{ J}$ (its copy by means of $\bullet$) with $\langle G_4\rangle$, and conversely, any direct product of a MAD-group of $\aut\mathcal{ J}$
 with $\langle G_4\rangle$ is a  MAD-group of $\aut\mathcal{ N}$.

The MAD-groups of $\aut \mathcal{ J}$ are completely described in \cite{f4}. According to it and with its  notations, there are four MAD-groups, described by
\begin{itemize}
\item $\{t_{x, y, z ,u }\mid x,y,z,u\in\F^*\}\cong(\F^*)^4$, which produces fine gradings on $\mathcal{ J}$ and $\der\mathcal{ J}$ of types $(24,0,1)$ and $(48,0,0,1)$ respectively.
\item $\{t_{x, y ,z, u }\mid x^2=y^2=z^2=u^2=1\}\times\langle \widetilde\sigma_{405}\rangle\cong \Z_2^5$, which produces fine gradings on $\mathcal{ J}$ and $\der\mathcal{ J}$ of types $(24,0,1)$ and $(24,0,0,7)$ respectively.
\item $\{t_{x ,y ,xy, u }\mid x^2=y^2=1,u\in\F^*\}\times\langle \widetilde\sigma_{105}\rangle\cong \Z_2^3\times \F^*$, which produces fine gradings on $\mathcal{ J}$ and $\der\mathcal{ J}$ of types $(25,1)$ and $(31,0,7)$ respectively.
\item $\langle \{t_{\omega,1,\omega^2,\omega^2},t_{1,\omega,\omega,1},\widetilde\sigma_{15}\}\rangle\cong \Z_3^3$, which produces fine gradings on $\mathcal{ J}$ and $\der\mathcal{ J}$ of types $(26)$ and $(0,26)$ respectively.
\end{itemize}
 Besides $\dim\mathcal{ J}_e=3,3,2,1$ respectively, so the type of the graded subspace $\mathcal{ J}_0$ is $(24,1)$ in the first and second cases, and $(26)$ in the third and fourth ones (taking into account that $1\in\mathcal{ J}$ belongs always to $\mathcal{ J}_e$). Consequently, the MAD-groups of $\aut\mathcal{ N}$ containing $G_4$ (equivalently, the MAD-groups of $\aut\eseis$ containing an order two automorphism fixing a subalgebra of type $\frak{f}_4$) are:
\begin{itemize}
\item $\Q_6:=\{t_{x, y, z ,u }^{\bullet}\mid x,y,z,u\in\F^*\}\times\langle G_4\rangle\cong\Z_2\times(\F^*)^4$, which produces a fine grading on $\mathcal{ N}$ of type $(72,1,0,1)$.
\item $\Q_7:=\{t_{x, y ,z, u }^{\bullet}\mid x^2=y^2=z^2=u^2=1\}\times\langle G_4,\widetilde\sigma_{405}^{\bullet}\rangle\cong \Z_2^6$, which produces a fine grading on $\mathcal{ N}$ of type $(48,1,0,7)$.
\item $\Q_8:=\{t_{x ,y ,xy, u }^{\bullet}\mid x^2=y^2=1,u\in\F^*\}\times\langle G_4,\widetilde\sigma_{105}^{\bullet}\rangle\cong \Z_2^4\times \F^*$, which produces a fine grading on $\mathcal{ N}$ of type $(57,0,7)$.
\item $\Q_9:=\langle \{G_4,t_{\omega,1,\omega^2,\omega^2}^{\bullet},t_{1,\omega,\omega,1}^{\bullet},\widetilde\sigma_{15}^{\bullet}\}\rangle\cong \Z_2\times\Z_3^3$, which produces a fine grading on $\mathcal{ N}$ of type $(26,26)$.
\end{itemize}

It is easy to give an automorphism conjugated to $G_4$ in terms of the model  in Equation~(\ref{eq_modeloAdams}).
If we take $G_4'$ the automorphism interchanging $V_1$ with $V_2$, then the fixed subalgebra is
$\{f_1+f_2\mid f\in\sll(V)\}\oplus \sll(V_3)\oplus \{(v_1\otimes v_2+v_2\otimes v_1) \otimes v_3\mid v_i\in V\}
\oplus S^2(V^*)\otimes V_3^*$, which is isomorphic to $\frak{f}_4$ (see, for instance, \cite{modelosf4}).
With this terminology, $\span{G_4',F_1,F_3,F_4}\cong\Z_2\times\Z_3^3$ is a MAD-group conjugated to $\Q_9$.
(Thus there are two fine gradings refining the Jordan grading).
Observe that
this MAD-group only contains one order two element, in particular it does not contain automorphisms of type $2D$.

\subsection{Gradings based on a $\frak{c}_4$-model.}\label{sub_losquevienendec4}

Take $G_5\in 2D$. There are some $8$-dimensional vector space $V$ and a non-degenerate symplectic bilinear form
$b\colon V\times V\to\F$ such that the subalgebra of $\mathcal{L} =\eseis$ fixed by $G_5$ is (isomorphic to)
$\frak{sp}(V,b)=\{f\in\End V\mid b(f(x),y)+b(x,f(y))=0\quad\forall x,y\in V\}$, a simple Lie algebra of type
$\frak{c}_4$. Recall from \cite[Chapter~8]{Kac} that, if $\mathcal{L}_{\bar0}\oplus
\mathcal{L}_{\bar1}$
is the $\Z_2$-grading induced by $G_5$, then
$\mathcal{L}_{\bar1}$ is an irreducible $\mathcal{L}_{\bar0}$-module. But the only irreducible $\frak{c}_4$-module of dimension $78-36=42$
is $V( \lambda_4)$ ($\lambda_i$'s the fundamental weights as in \cite{Humphreysalg}). A suitable way of describing it is as a submodule of $\wedge^4V$, which has dimension $\binom84=70$.
The decomposition of this module into its irreducible summands is $\wedge^4V\cong V( \lambda_4)\oplus V( \lambda_2)\oplus V(0)$.
Thus, if we consider the contraction
$$
\begin{array}{llcl}
c\colon  &\bigwedge^4V&\to &\bigwedge^2V \  \  (\cong V(\lambda_2))\\
&v_1\wedge v_2\wedge v_3 \wedge v_4  &\mapsto&\sum_{\tiny
{\begin{array}{l}\sigma\in S_4\\
\sigma(1)<\sigma(2)\\
\sigma(3)<\sigma(4)\end{array}}}
(-1)^\sigma b(v_{\sigma(1)},v_{\sigma(2)})
v_{\sigma(3)}\wedge v_{\sigma(4)}
\end{array}
$$
for $S_4$ the group of permutations of $\{1,2,3,4\}$, its kernel $\ker c$
is isomorphic to $\mathcal{L}_{\bar1}+\F$, with the action on $\F$ trivial.    
This construction seems not to be natural, but one only has to recall that the Lie algebra of type $\frak{e}_7$
can be modeled as   $\sll(V)\oplus \wedge^4V$ and that $\eseis$ lives here.
The main reason for use it is that it is quite easy to extend the automorphisms
of $\frak{c}_4$ until the whole $\eseis$.
Recall that $\aut\mathfrak{c}_4\cong\SP(V)$ and take the map
$$
\SP(V)=\{f\in\End V\mid b(f(x), f(y))=b(x,y)\ \forall x,y\in V\}\rightarrow\aut\mathcal{L},\
f\mapsto f^\diamondsuit
$$
given by
$f^\diamondsuit(g)=f^{-1}gf$ if $g\in \mathcal{L}_{\bar0}=\frak{sp}(V,b)$
and $f^\diamondsuit(v)=\sum_{ }f(v_{i_1})\wedge f(v_{i_2})\wedge f(v_{i_3}) \wedge f(v_{i_4}) $
 if $v=\sum_{ } v_{i_1}\wedge v_{i_2}\wedge v_{i_3} \wedge v_{i_4} \in\mathcal{L}_{\bar1}\subset \ker c$.
It is a   computation very similar to that one in Subsection~\ref{sub_losquevienendef4} that $\cent_{\aut \mathfrak{e}_6}(G_5)=\{f^\diamondsuit\mid f\in\SP(V)\}\cdot\span{G_5}\cong\SP(V)\times\Z_2$.
This implies that every MAD-group of $\aut\eseis$ containing $G_5$ (that is, containing some automorphism of the isomorphy
class $2D$)
is of the form $\{f^\diamondsuit\mid f\in Q\}\cdot\span{G_5}$ for some MAD-group $Q$ of $\SP(V)$.
There are $7$ MAD-groups of $\aut\frak{c}_4$, according to \cite{LGII,Albclasicas}. The induced fine gradings can also be extracted from \cite{Ivan2}, although
in such paper there is one missing  grading.
Although we know that we will obtain just $7$ MAD-groups of $\aut\eseis$ by means of this procedure, we don't know a priori how many of these
MAD-groups have appeared before (equivalently, how many of these MAD-groups contain an outer order two automorphism of type $2C$).
Thus we are going to recall the descriptions of these quasitori, and extend each automorphism of $\frak{c}_4 $ to $\eseis$ to get the complete
simultaneous diagonalizations.
We will make use of the notations of \cite{LGII} for giving the MAD-groups of $\SP(V)$. According to it, the MAD-groups are
$$
\begin{array}{c}
 \Xi_1=T_{8,0}^{(0)};\
 \Xi_2=T_{4,0}^{(1)}(I_4,I_4);\
 \Xi_3=T_{2,2}^{(1)}(I_4,I_4);\
 \Xi_4=T_{0,4}^{(1)}(I_4,I_4);\\
 \Xi_5=T_{0,2}^{(2)}((\tiny{\begin{pmatrix}1&0\cr0&\imag\end{pmatrix}},I_2),(I_2,I_2));\\
 \Xi_6=T_{2,0}^{(2)}((I_2,I_2),(I_2,I_2));\
 \Xi_7=T_{0,1}^{(3)}((1,1),(1,1),(1,1)),
 \end{array}
$$
for $\imag\in\F$ such that $\imag^2=-1$.
We try to avoid that the reader has to dive into the details of such paper by providing the exact descriptions of the automorphisms
in the $\Xi_i$'s.
We will use the so called \emph{Pauli's matrices}, given by
$$
\theta_1=\begin{pmatrix}0&1\cr1&0\end{pmatrix}, \quad
\theta_3=\begin{pmatrix}1&0\cr0&-1\end{pmatrix},\quad
\theta_2=\theta_3\theta_1=\begin{pmatrix}0&1\cr-1&0\end{pmatrix}.
$$
As usual, if $A=(a_{ij})\in{\mathop{\rm Mat}}_{m\times n}(\F)$ and $B\in{\mathop{\rm Mat}}_{p\times q}(\F)$,
the Kronecker product $A\otimes B$ denotes the block-matrix
$$
 A\otimes B=\tiny{\begin{pmatrix}
 a_{11}B&\dots&a_{1n}B \cr
 \vdots&\ddots&\vdots\cr
 a_{m1}B&\dots&a_{mn}B\end{pmatrix}}\in{\mathop{\rm Mat}}_{mp\times nq}(\F) .
$$
\begin{itemize}
\item $\Xi_1$. This is the case of the maximal $4$-dimensional torus of $\SP(V)$, which produces a $\Z^4$-grading on $\frak{c}_4$
of type $(32,0,0,1)$.

\item $\Xi_2$. We fix some basis $\mathcal{B}=\{w_1,\dots,w_8\}$ of $V$ such that the matrix of the form $b$
relative to $\mathcal{B}$ is  the skew-symmetric matrix $I_2\otimes \theta_1\otimes \theta_2$ (take into account that
$(A\otimes B)^t=A^t\otimes B^t$).
Take in $\SP(V,b)$ the automorphisms of $V$ whose related matrices in the basis $\mathcal{B}$ are
$$\Xi_2= \span{r_{\alpha,\beta},I_4\otimes \theta_1,I_4\otimes \theta_3\mid \alpha,\beta\in\F^*}\cong(\F^*)^2\times\Z_2^2,$$
for $r_{\alpha,\beta}=\diag\{\alpha,\frac1\alpha,\beta,\frac1\beta\}\otimes I_2$. The induced $\Z^2\times\Z_2^2$-grading on $\frak{c}_4$ is of type $(28,4)$.

\item  $\Xi_3$. We fix some basis   of $V$ such that the matrix of the form $b$
relative to it is  the skew-symmetric matrix $\diag\{I_2,\theta_1\}\otimes \theta_2$, and take
$$
\Xi_3= \span{r_{1,\alpha},I_4\otimes \theta_1,I_4\otimes \theta_3,\diag\{-I_2,I_2,I_2,I_2\}\mid  \alpha\in\F^*}\cong \F^* \times\Z_2^3.
$$
The induced $\Z \times\Z_2^3$-grading on $\frak{c}_4$ is of type $(27,0,3)$.

\item $\Xi_4$. We fix some basis  of $V$ such that the matrix of the form $b$
relative to it is  $I_4\otimes\theta_2$, and take
$$
\Xi_4= \span{ I_4\otimes \theta_1,I_4\otimes \theta_3,\diag\{-1,1,1,1\}\otimes I_2,\diag\{1,-1,1,1\}\otimes I_2,\diag\{ 1,1,-1,1\}\otimes I_2 }\cong  \Z_2^5.
$$
The induced $\Z_2^5$-grading on $\frak{c}_4$ is of type $( 24,0,0,3)$.

\item  $\Xi_5$. We fix a basis   of $V$
relative to  which the matrix of the form $b$ is  $\diag\{\theta_2\otimes I_2,\theta_2\otimes\theta_1\}$, and take
$$
\Xi_5= \span{\tiny{\begin{pmatrix}1&0\cr0&\imag\end{pmatrix}}\otimes I_2\otimes \theta_3  , I_2\otimes \theta_1\otimes I_2,I_2\otimes \theta_3\otimes I_2,I_4\otimes \theta_1 }\cong\Z_4\times  \Z_2^3.
$$
The induced $\Z_4\times  \Z_2^3$-grading on $\frak{c}_4$ is of type $( 24,6 )$.

\item  $\Xi_6$. We fix some basis
relative to which the matrix of the form $b$ is    $  \theta_1 \otimes \theta_2\otimes I_2$, and take
$$
\Xi_6= \span{\diag\{ \alpha,\frac1\alpha\}\otimes I_4,I_4\otimes \theta_1,I_4\otimes \theta_3,I_2\otimes \theta_1\otimes I_2,I_2\otimes \theta_3\otimes I_2\mid  \alpha\in\F^*}\cong \F^* \times\Z_2^4.
$$
The induced $\Z \times\Z_2^4$-grading on $\frak{c}_4$ is of type $( 36 )$.

\item $\Xi_7$. We fix some basis
relative to which the matrix of the form $b$ is    $  \theta_2\otimes I_4$, and take
$$
\Xi_7=  \span{ I_4\otimes \theta_1,I_4\otimes \theta_3,I_2\otimes \theta_1\otimes I_2,I_2\otimes \theta_3\otimes I_2,\theta_1\otimes I_4,\theta_3\otimes I_4}\cong  \Z_2^6.
$$
The induced $ \Z_2^6 $-grading on $\frak{c}_4$ is of type $( 36 )$.
\end{itemize}
Take into consideration now that $G_5\cdot\diag\{-I_2,I_2,I_2,I_2\}^\diamondsuit$ is an automorphism of the class $2C$, and hence that
we do not get anything new from $\Xi_1$, $\Xi_3$ and $\Xi_4$.


Extending the automorphisms and   making the simultaneous diagonalization is a tedious task, although straighforward. Let us provide some details of the first of
our significant cases, $\Xi_2$. Denote by $w_{ijkl}=w_i\wedge w_j\wedge w_k\wedge w_l\in \wedge^4V$, so that $\{w_{ijkl}\mid
1\le i<j<k<l\le8\}$ is a basis of $\wedge^4V$. Denote by $L_{(i,j,k)}=\{x\in\ker c\mid
r_{\alpha,\beta}(x)=ix,I_4\otimes \theta_1(x)=jx,I_4\otimes \theta_3(x)=kx\}$ the homogeneous components
of the grading produced by $G_5$ and $\Xi_2$ restricted to the odd part.
Then all the homogeneous components are one-dimensional (for instance, $L_{(\alpha^2\beta^2,1,1)}=
\span{w_{1256}}$ and so on) except for
one four-dimensional component
$$
\begin{array}{lr}
L_{(1,1,1)}=&\span{w_{ 1368}+w_{ 2457},w_{ 1458}+w_{1467 }+w_{2358 }+w_{2367 },\\
&w_{ 1357}+w_{ 2468},w_{1234 }-w_{1467 }-w_{ 2358}+w_{5678 }   }
\end{array}
$$
and three two-dimensional components
$$
\begin{array}{l}
L_{(1,-1, 1)}=\span{ w_{ 1357}-w_{ 2468},w_{ 1368}-w_{ 2457}},\\
L_{(1, 1,-1)}=\span{  w_{1358 }+w_{ 1367}+w_{ 2458}+w_{2467 },w_{1457 }+w_{2357 }+w_{2368 }+w_{1468 }   },\\
L_{(1,-1,-1)}=\span{  w_{1358 }+w_{ 1367}-w_{ 2458}-w_{2467 },w_{1457 }+w_{2357 }-w_{2368 }-w_{1468 }  }.
\end{array}
$$
The type of the induced grading in $\eseis$ is hence
$(28,4,0,0)+(32,3,0,1)=(60,7,0,1)$.

In the same way, we check that the type of the grading induced by  $\Xi_6$  restricted to $\mathcal{L}_{\bar1}=\ker c$
is $(37,0,0,0,1)$ because all the homogeneous components have dimension $1$ except for one of dimension $5$:
$$\begin{array}{lr}
\fix(\Xi_6)\cap\mathcal{L}_{\bar1}=&\span{
w_{1368 }+w_{2457 },
w_{ 1458}+w_{2367 },
w_{1357 }-w_{1467 }-w_{ 2358}+w_{2468 },\\
&w_{1256 }+w_{3478 },
w_{1278 }-w_{ 1467}-w_{ 2358}+w_{3456 }
}.
\end{array}
$$
Observe that in this case we are not dealing with the same contraction $c$ than for $\Xi_2$, although we use the same notation.

Finally, it is trivial to compute the types of the gradings induced by  $\Xi_5$  and $\Xi_7$ restricted to $\mathcal{L}_{\bar1}=\ker c$,
which are, respectively, $(24,7,0,1)$ and $(36,0,0,0,0,1)$.

Consequently, the MAD-groups of $\aut\mathcal{ L}$ containing $G_5$ (equivalently, the MAD-groups of $\aut\eseis$ containing an order two automorphism fixing a subalgebra of type $\frak{c}_4$) but not conjugated to any MAD-group in Subsection~\ref{sub_losquevienendef4} are:
\begin{itemize}
\item $\Q_{10}:=\{f^\diamondsuit\mid f\in\Xi_2\}\cdot\span{G_5}\cong(\F^*)^2\times\Z_2^3$, which induces a  fine grading of type $( 60,7,0,1  )$.
\item $\Q_{11}:=\{f^\diamondsuit\mid f\in\Xi_5\}\cdot\span{G_5}\cong \Z_4\times\Z_2^4$,  which induces a  fine grading of type $(48,13,0,1   )$.
\item $\Q_{12}:=\{f^\diamondsuit\mid f\in\Xi_6\}\cdot\span{G_5}\cong\F^* \times\Z_2^5$,  which induces a  fine grading of type $( 73,0,0,0,1  )$.
\item $\Q_{13}:=\{f^\diamondsuit\mid f\in\Xi_7\}\cdot\span{G_5}\cong \Z_2^7$,  which induces a  fine grading of type $( 72,0,0,0,0,1  )$.
\end{itemize}

\subsection{A MAD-group containing outer automorphisms but without outer order two automorphisms.}\label{sub_laZ4cuborara}

As was mentioned in Introduction, there are simple Lie algebras $L$ with MAD-groups of $\aut L$ not contained in $\Int L$ but without outer involutions. In this section we will describe how this situation occurs  again for $L$ of type $\eseis$.

We are going to describe a $\Z_4^3$-fine grading on $\eseis$.
Take $\Upsilon_1$ any  outer  automorphism of order $4$  which fixes a subalgebra of type $\frak{a}_3\oplus\frak{a}_1$, as in Subsection~\ref{sub_autexternos}.

\begin{re}
We can find an automorphism of this conjugacy class with our descriptions in Subsection~\ref{sub_losquevienendef4}.
For instance,  the automorphism $t=t_{-\imag, \imag,-1,\imag}$ produces a $\Z_4$-grading on $\mathcal{J}$
with components of dimensions $9,6,6$ and $6$, and $\Ad t$ produces a $\Z_4$-grading on $\der\mathcal{J}$
with components of dimensions $12,14,12$ and $14$. Thus the subalgebra fixed by
$G_4 t^{\bullet}\in \Q_6$ is
the sum of the elements of  $\der\mathcal{J}$ fixed by $\Ad t$ with the elements of $\mathcal{J}_0$
antifixed by $t$, whose dimension is $12+6=18$, so that $G_4 t_{-\imag, \imag,-1,\imag}^{\bullet}$ is conjugated to $\Upsilon_1$.
And in terms of the notations in Section~\ref{sec descripcionesinternos}, it corresponds to the automorphism $G_4''T_{1,\imag}$
if $G_4''$ is the automorphism interchanging $V_1$ with  $V_3$ in Equation~(\ref{eq_modeloAdams}).
\end{re}

Consider now the $\Z_4$-grading on $\eseis$ induced by $\Upsilon_1$, $\mathcal{L}=\mathcal{L}_{\bar0}\oplus\mathcal{L}_{\bar1}\oplus
\mathcal{L}_{\bar2}\oplus
\mathcal{L}_{\bar3}$. Recall that $\mathcal{L}_{\bar i}$ must be an $\mathcal{L}_{\bar0}$-irreducible module (if $i\ne0$), that is,
a tensor product of an irreducible $\frak{a}_3$-module with an  irreducible $\frak{a}_1$-module.
If we take also into account
that $\dim\hom_{\mathcal{L}_{\bar 0}}(\mathcal{L}_{\bar i}\otimes\mathcal{L}_{\bar j},\mathcal{L}_{\bar i+\bar j})=1$, then
the unique possibility for the decomposition of $\mathcal{L}$ as a sum of $\mathcal{L}_{\bar0}$-modules is
$$
\frak{a}_3\oplus\sll(V)\oplus\,V(2\lambda_1)\otimes V\,\oplus\, V(2\lambda_2)\otimes\F\,\oplus\,V(2\lambda_3)\otimes V
$$
 for $V$ a two-dimensional vector space, and $\{\lambda_i\mid i=1,2,3\}$    the set of fundamental weights for $\frak{a}_3$.
 The dimension  of each non-identity component is $20$.

Take $W$ a four-dimensional vector space, so that $\mathfrak{a}_3\cong\sll(W)$. The natural module $W$ is isomorphic to
$V(\lambda_1)$ and its second symmetric power $S^2(W)$ is isomorphic to $V(2\lambda_1)$.
Their dual $\sll(W)$-modules, $W^*$ and $S^2(W^*)$, are respectively of types $V(\lambda_3)$ and $V(2\lambda_3)$.
Finally  $\bigwedge^2W$ is of type $V(\lambda_2)$ and its second symmetric power
$S^2(\bigwedge^2W)\cong V(2\lambda_2)\oplus V(0)$. So consider $S^2(\bigwedge^2W)'$       its only non-trivial submodule
(of type $V(2\lambda_2)$). We have an isomorphism of $\sll(W)\oplus\sll(V)$-modules between $\eseis$ and
$$
\mathcal{N}:=\sll(W)\oplus\sll(V)\oplus\,S^2(W)\otimes V\,\oplus\, S^2(\bigwedge^2W)'\otimes\F\,\oplus\,S^2(W^*)\otimes V^*
$$
which endows $\mathcal{N}$ with a Lie algebra structure $\mathbb{Z}_4$-graded such that
$\mathcal{N}_{\bar1}=S^2(W)\otimes V$, $\mathcal{N}_{\bar2}=S^2(\bigwedge^2W)'\otimes\F=[\mathcal{N}_{\bar1},\mathcal{N}_{\bar1}]$ and
$\mathcal{N}_{\bar3}=S^2(W^*)\otimes V^*=[\mathcal{N}_{\bar2},\mathcal{N}_{\bar1}]$.
Thus $\Upsilon_1$ can be considered as the automorphism of $\mathcal{N}\cong\eseis$ producing this $\mathbb{Z}_4$-grading.

Now, for each $f\in\End W$ and  $g\in\End V$, denote by $\exte(f\otimes g)\in\End \mathcal{N}_{\bar1}$
 the map given by
 $\exte(f\otimes g)(w\cdot w'\otimes v)=f(w)\cdot f(w')\otimes g(v)$,
 for all $w,w'\in W$ and $v\in V$, where $\cdot$ denotes the symmetric product.

Define the automorphisms $\Upsilon_2$ and $\Upsilon_3$ as the only automorphisms of $ \mathcal{N}\,(\cong\eseis)$ whose restrictions
to $\mathcal{N}_{\bar1}$ are:
$$
\begin{array}{l}
\Upsilon_2\vert_{ \mathcal{N}_{\bar1} }=\exte\left( \tiny{\begin{pmatrix}
 0&0&0 & 1  \cr
 1&0&0&0 \cr
 0&1&0&0 \cr
 0&0&1&0
 \end{pmatrix}}_W\otimes \begin{pmatrix}
 0&1 \cr
 1&0\end{pmatrix}_V\right)\\
 \Upsilon_3\vert_{ \mathcal{N}_{\bar1} }=\exte\left( \tiny{\begin{pmatrix}
 1&0&0 & 0  \cr
 0&\imag&0&0 \cr
 0&0&-1&0 \cr
 0&0&0&-\imag
 \end{pmatrix}}_W\otimes \begin{pmatrix}
 1&0 \cr
 0&-1\end{pmatrix}_V\right)
 \end{array}
$$
where we have chosen $\{w_0,w_1,w_2,w_3\}$ and $\{v_0,v_1\}$ basis of $W$ and $V$ respectively
and we have identified the endomorphisms of $W$ and $V$ with their matrices in such bases.
Thus $\{w_j\cdot w_k\otimes v_l\mid 1\le j\le k\le 4,\,l=0,1\}$ is a basis of $\mathcal{N}_{\bar1}$,
and
$\Upsilon_2(w_j\cdot w_k\otimes v_l)=w_{j+1}\cdot w_{k+1}\otimes v_{l+1}$
($j$ and $k$ summed modulo $4$ and $l$ summed modulo $2$)
and
$\Upsilon_3(w_j\cdot w_k\otimes v_l)=(\imag)^{j+k+2l}w_{j }\cdot w_{k }\otimes v_{l }$.
It is a straightforward computation that the extensions  $\Upsilon_2$ and $\Upsilon_3$ are Lie algebra automorphisms.

Now $\Upsilon_2\Upsilon_3(w_j\cdot w_k\otimes v_l)=(\imag)^{j+k+2l}w_{j +1}\cdot w_{k +1}\otimes v_{l+1 }=
\Upsilon_3\Upsilon_2(w_j\cdot w_k\otimes v_l)$, so that $\Upsilon_2$ and $\Upsilon_3$ commute.
Consider, then,  the quasitorus of $\aut\eseis$ given by
$$
\bold{\Q_{14}}:=\span{\Upsilon_1,\Upsilon_2,\Upsilon_3}\cong\Z_4^3
$$
Let us compute the simultaneous diagonalization of  $ \mathcal{N}$ relative to $\Q_{14}$. In Proposition~\ref{pr_elZ4cuboenterminoscomputacionales}
we will prove that it is a maximal quasitorus and hence the induced grading is fine.
Again the common diagonalization is a question of patience.
Denote by $L_{(\bar i,\bar j, \bar k)}=\{x\in\mathcal{N}\mid \Upsilon_1(x)=ix, \Upsilon_2(x)=jx,\Upsilon_3(x)=kx\}$.
Of course $\langle\Upsilon_2\vert_{\fix \Upsilon_1},\Upsilon_3\vert_{\fix \Upsilon_1}\rangle$ produce the non-toral $\Z_4^2$-grading on $\frak{a}_3$ with
 trivial identity component and all the remaining $15$ components of dimension $1$, and produce the non-toral $\Z_2^2$-grading on $\frak{a}_1$
 given by Pauli's matrices. Thus $\dim L_{(\bar 0 ,\bar 0 , \bar  0)}=0$,
 $\dim L_{(\bar 0 ,\bar  2, \bar  0)}=\dim L_{(\bar  0,\bar  0, \bar  2)}=\dim L_{(\bar  0,\bar 2 , \bar 2 )}=2$
 and $\dim L_{(\bar 0 ,\bar i , \bar  j)}=1$ for the remaining $i,j$. So the type of the restriction to
 $\mathcal{N}_{\bar0}\cong\mathfrak{a}_3\oplus\mathfrak{a}_1$ is $(12,3)$. Now consider the restriction to
 $\mathcal{N}_{\bar1}$. Observe that  $\Upsilon_2^2(x=w_0\cdot w_0\otimes v_k)\ne w_0\cdot w_0\otimes v_k$, so that
 $\Upsilon_2$ acts with eigenvalue $ \varepsilon\in\{1,i,-1,-i\}$
 in $x+\varepsilon^3\Upsilon_2(x)+\varepsilon^2\Upsilon_2^2(x)+\varepsilon\Upsilon_2^3(x)$, and the same happens to
 $x=w_0\cdot w_1\otimes v_k$. On the contrary, $\Upsilon_2^2(x=w_0\cdot w_2\otimes v_k)= w_0\cdot w_2\otimes v_k$,
 and hence $\Upsilon_2$ acts with eigenvalue $\pm{1}$ in $x\pm\Upsilon_2(x)$.
 Hence $\dim L_{(\bar 1 ,\bar 0 , \bar  0)}=\dim L_{(\bar 1 ,\bar  2, \bar  0)}=\dim L_{(\bar  1,\bar  0, \bar  2)}=\dim L_{(\bar  1,\bar 2 , \bar 2 )}=2$
 and $\dim L_{(\bar 1 ,\bar i , \bar  j)}=1$ for the remaining $i,j$, and the type of the restriction to $\mathcal{N}_{\bar1}$ is $(12,4)$.
 It is not difficult to work with the other two components to arrive at the conclusion that both $\mathcal{L}_{\bar2}$ and
 $\mathcal{L}_{\bar3}$ break into $12$ one-dimensional components and $4$ two-dimensional components, therefore the type of the grading induced by $\Q_{14}$ in $\eseis$ is $(12,3)+3(12,4)=(48,15)$.
 More precisely, the identity component is trivial, $\dim L_{(\bar i ,\bar j , \bar  k)}=2$ if $j,k\in\{0,2\}$ and the remaining homogeneous componenets are one-dimensional.
 In particular, $\dim\fix\span{\Upsilon_1^2,\Upsilon_2,\Upsilon_3}=\dim L_{(\bar 2 ,\bar 0 , \bar  0)}=2<6$, and hence the quasitorus
 $\span{\Upsilon_1^2,\Upsilon_2,\Upsilon_3}\cong\mathbb{Z}_{2}\times\mathbb{Z}_{4}^2$ is non-toral. Observe also that $ \fix\span{ \Upsilon_2,\Upsilon_3}$ is just a Cartan subalgebra
 and hence $\span{ \Upsilon_2,\Upsilon_3}$ is contained in a torus. These considerations will be useful to argument with $\Q_{14}$.

\subsection{All the outer fine gradings}

As a corollary of Propositions~\ref{pr_elZ4cuboenterminoscomputacionales} and~\ref{pr_sinohayextorden2eselQ14} in the next technical section, we will obtain   the other main result  of the paper:

\begin{te}
The MAD-groups of $\aut\eseis$ not contained in $\Int\eseis$ are, up to conjugation,  $\Q_i$ for $i=6,\dots,14$.
\end{te}

\textbf{Proof.}
  If $A$ is a MAD-group of $\aut\eseis$ such that $A \nsubseteq \Int\eseis$,
  then either $A$ has an outer automorphism of type $2C$,
  and in such case $A$ is conjugated to some $\Q_i$ with $i\in\{6,7,8,9\}$;
  or $A$ has an outer automorphism of type $2D$,
  and in such case $A$ is conjugated to some $\Q_i$ with $i\in\{6,7,8,10,11,12,13\}$;
  or $A$ has not any outer automorphism of order two,
  and in such case $A$ is conjugated to   $\Q_{14}$ by Proposition~\ref{pr_sinohayextorden2eselQ14}.
  The quasitorus  $\Q_{14}$ appears in this list by Proposition~\ref{pr_elZ4cuboenterminoscomputacionales}.
\hfill  $\square$ \smallskip


\section{Technical proofs for the outer gradings.}\label{sec_demostracionesouter}

\subsection{Extended Weyl group}

The diagram automorphism interchanging $\a_3$ with $\a_5$, and  $\a_1$ with $\a_6$, has matrix relative to
$\Delta$:
$$\sigma=\tiny{\begin{pmatrix}
 0&0&0 & 0 & 0 & 1\cr
 0&1&0&0&0&0\cr
 0&0&0&0&1&0\cr
 0&0&0&1&0&0\cr
 0&0&1&0&0&0\cr
 1&0&0&0&0&0
 \end{pmatrix}}
 $$
and the extended Weyl group is   $\mathcal{V}:=\mathcal{W}\cup \mathcal{W}\sigma\cong\mathcal{W}\rtimes\Z_2$, which is
the set of automorphisms of the root system.


We will make an extensive use of those orbits in $\mathcal{V}$ whose representatives have   order a power of $2$. We find $10$ new orbits, and we summarize the information related to them in the following table:

\begin{center}
\begin{tabular}{|c|c|c|c|c|}
\hline Order& Representative   & Stabilizer  &
Isomorphic to\cr \hline\hline
 $2$ &$\si\equiv\eta_1$   & $t_{x,y,z,u,z,x}$ & $(\F^*)^4 $\cr
 $2$ &$\si\sigma_{555}\equiv\eta_2$   & $ t_{\frac{z}{w^2y^4},y,z,\frac1z,y,w}$ & $(\F^*)^3$\cr
 $2$ &$\si\sigma_{458}\equiv\eta_3$   & $t_{x,vx,z,x,v,\frac{\alpha}{zv^2}}\mid \tiny{x^2=\alpha^2=1} $ & $(\F^*)^2\times\Z_2^2$\cr
 $2$ &$\si\sigma_{2402}\equiv\eta_4$   & $ t_{x,y,z,u,xy,\frac\alpha z}\mid \tiny{x^2=y^2=u^2=\alpha^2=1}$ & $\F^*\times\Z_2^4$\cr
 $2$ &$-\id=\eta_5$   & $t_{x,y,z,u,v,w}\mid \tiny{x^2=y^2=z^2=u^2=v^2=w^2=1} $ & $\Z_2^6$\cr
 $4$ &$\si\sigma_{15}\equiv\mu_1$   & $t_{x,y,\frac1{x^3y},x,x^3y,\frac1{x^2y}} $ & $(\F^*)^2$\cr
 $4$ &$\si\sigma_{52}\equiv\mu_2$   & $t_{x,\frac{zx^2}{w^2},z,\frac{w^2}{x^3z^2},\frac{zx^2}{w^2},w} $ & $(\F^*)^3$\cr
 $4$ &$\si\sigma_{460}\equiv\mu_3$   & $t_{x,y,x,x,xy,\frac{\a}{y^2}}\mid \tiny{x^2=\alpha^2=1} $ & $\F^{* }\times\Z_2^2$\cr
 $4$ &$\si\sigma_{484}\equiv\mu_4$   & $t_{x,y,x,x,xy,x^2}\mid \tiny{x^4=y^4=1} $ & $\Z_4^2$\cr
 $8$ &$\si\sigma_{17}$   & $ t_{x,\frac1{x^2},\frac1x,x,x,1}$ & $\F^*$\cr
 \hline
\end{tabular}
\end{center}\begin{center}\begin{equation}\text{ Table of representatives of $\W\sigma$ of order a power of $2$}\label{tablaexternos} \end{equation}\end{center}
We employ the notations $\eta_i,\mu_i$ for the order $2$ and $4$ representatives, in order to shorten the use of indices.
The only restriction about the scalars in the third column when nothing is said is that they are non-zero.\vskip0.5cm

This time the projection $\pi\colon\No(\T)\to \No(\T)/\T\cong\mathcal{V}$ is an extension of the projection $\pi$ considered in
Subsection~\ref{subsec_Weyl}. For each $\nu\in\mathcal{W}\sigma$  we   choose some element $\widetilde{\nu }\in\pi^{-1}(\nu )$
of minimum order among the elements in $\pi^{-1}(\nu )$. (Perhaps our choice of $\tilde\nu$ when $\nu\in\mathcal{W}$ does not coincide with that one in  Subsection~\ref{subsec_Weyl}, where we made a more concrete election not based in the order, but this does not interphere with our next arguments).

\begin{re}\label{re_extensionesdeoreden2}
A consequence of the Appendix, and of our elections of extensions, is that all those extensions $\widetilde\eta_j$ have order $2$. It will extremely
 useful in the next proofs that this fact (the existence of some order two extensions) jointly with Remark~\ref{re_afinandoelSsuperf}  imply that $(\widetilde\eta_js)^2\in\sor{\eta_j}$ for all $s\in\T$ and for all $j$.
\end{re}

\subsection{The $\Z_4^3$-grading in computational terms}

Recall that in Subsection~\ref{sub_laZ4cuborara} we found a quasitorus $\Q_{14}=\span{\Upsilon_1,\Upsilon_2,\Upsilon_3}$ isomorphic
as abstract group
to $\Z_4^3$ satisfying the following conditions:
\begin{itemize}
\item $\Upsilon_1$ is an outer automorphism;
\item $\span{\Upsilon_2,\Upsilon_3}$ is toral;
\item $\span{\Upsilon_1^2,\Upsilon_2,\Upsilon_3}\subset\Int\eseis$ is non-toral.
\end{itemize}

A first consequence is that

\begin{pr}\label{pr_elZ4cuboenterminoscomputacionales}
$\Q_{14}$  is  conjugated
to $Q(\widetilde\mu_4)$. Moreover, $\Q_{14}$ is a MAD-group.
\end{pr}

 \textbf{Proof.}
By Lemma~\ref{le_puedoajustarlaparatedeltoro}, we can assume that $\span{\Upsilon_2,\Upsilon_3}\subset\T$ and that
$\Upsilon_1\in\No(\T)$ (by conjugating, if necessary).
So $\pi(\Upsilon_1)$ is an element in $\V\setminus \W$ of order a divisor of $4$. But such order cannot be $2$, since in such a case
we could suppose that $\Upsilon_1=\widetilde\eta_js$ for some $j=1,\dots,5$ and $s\in\T$ and then $\Upsilon_1^2\in\T$, so that $\span{\Upsilon_1^2,\Upsilon_2,\Upsilon_3}\subset\T$ would be  toral.
Thus there is some $j=1,\dots,4$ and some $s\in\T$ such that
$\Upsilon_1=\widetilde\mu_js$ (again after conjugating) and $\span{\Upsilon_2,\Upsilon_3}\subset\tor{\mu_j}$.
Note that $j\ne3$ because there is no subgroup isomorphic to $\Z_4^2$ contained in
$\tor{\mu_3}\cong\F^*\times\Z_2^2$.
 Moreover, $j\ne1,2$ because in such cases $\tor{\mu_j}$ is a torus, so that $\span{\Upsilon_1^2,\Upsilon_2,\Upsilon_3}\subset\span{\Upsilon_1^2,\tor{\mu_j}}$ would be toral (we can apply Lemma~\ref{le_toroporalgotoral} because $\Upsilon_1^2$ is inner). Hence $\Q_{14}\subset Q(\widetilde\mu_4s)$ and they must coincide ($\widetilde\mu_4s$ has order $4$ and both
quasitori are then isomorphic to $\Z_4^3$). Besides $Q(\widetilde\mu_4s)\cong Q(\widetilde\mu_4)$ by Remark~\ref{re_paramover}, since $\tor{\mu_4}$ is finite.

Now we are going to check that $Q(\widetilde\mu_4 )$ is its own centralizer.  
Let us take $f\in\cent_{\aut \mathfrak{e}_6}(Q(\widetilde\mu_4))$. 
Consider $Z=\cent_{\aut \mathfrak{e}_6}(\tor{\mu_4})$.
There exists $\T'$ a maximal torus of $Z$ such that $Q(\widetilde\mu_4 )\cup\{f\}$ is contained in the normalizer
of $\T'$. As all the maximal tori of $Z$
are conjugated, there is $p\in Z$ such that $p\T'p^{-1}=\T$, hence
$ptp^{-1}=t$ for all $t\in\tor{\mu_4}$ and $p(Q( \widetilde\mu_4 )\cup\{f\})p^{-1}\subset\No(\T)$.
Take into account that $\{\nu\in\sigma\W\mid \nu\cdot t=t\quad\forall t\in\tor{\mu_4}\}=\{\mu_4,\mu_4^3\}$
and therefore there are $l\in\{1,3\}$ and $s\in\T$ such that
 $p\widetilde\mu_4p^{-1}=\widetilde\mu_4^ls$.
 Hence
  $pQ(\widetilde\mu_4)p^{-1}=Q(\widetilde\mu_4^ls)$.
  As in Remark~\ref{re_paramover}, there is $s'\in\T$
  such that $s'(\widetilde\mu_4^ls)(s')^{-1}\in\widetilde\mu_4^l\tor{\mu_4}$, so that $\Ad p'$ for $p'=s'p\in\aut\eseis$ leaves invariant
  $Q(\widetilde\mu_4)$ (and fixes $\tor{\mu_4}$ pointwise).
  As $p'fp'^{-1}$ also belongs to $\No(\T)$ and commutes with $\tor{\mu_4}$,
  and taking into account that
  $\{\nu\in \W\mid \nu\cdot t=t\quad\forall t\in\tor{\mu_4}\}=\{\id,\mu_4^2\}$,
then
$p'fp'^{-1}=\widetilde\mu_4^rs''$ for some $r\in\{0,1,2,3\}$ and $s''\in\T$.
But $p'fp'^{-1}$ commutes with $p'\widetilde\mu_4p'^{-1}\in\widetilde\mu_4^l\tor{\mu_4} $ and with $\tor{\mu_4}$, so that  $p'fp'^{-1}$ commutes with $\widetilde\mu_4$ 
and hence
$s''$ also does, in other words $s''\in\tor{\mu_4}$.
This means that
  $p'fp'^{-1}=\widetilde\mu_4^rs''\in  Q(\widetilde\mu_4 )=p'Q(\widetilde\mu_4 )p'^{-1}$ and hence that
 $f\in Q(\widetilde\mu_4 )$.
\hfill  $\square$ \smallskip

In particular,  the  outer automorphism $\widetilde\mu_4$ has order $4$.

It is clear  that $\Q_{14}$ does not contain any order two outer automorphism, since any outer automorphism in
$Q(\mu_4,\id)$ belongs to the set $\{\widetilde\mu_4s,\widetilde\mu_4^3s\mid s\in\tor{\mu_4}\}$ and its  square
does not belong to $\T$ (in particular, its square is not $\id$). 

Finally, a last property satisfied by $\Q_{14}$ is the following:
\begin{itemize}
 \item $\span{\Upsilon_1^2,\Upsilon_2^2,\Upsilon_3}\subset\Int\eseis$ is  toral.
 \end{itemize}

 \begin{re}\label{re_esemododeencontrarZ4cubonoessorpresa}
This fact implies that there are $j=1,\dots,5$, $\sigma_i$ in the orbit of $\sigma_{96}$ and $s,s'\in\T$ such that
$\Q_{14}$ is conjugated to $\langle \widetilde{\eta}_js,\widetilde{\sigma}_{i}s',\tor{\sigma_{i}}\cap\tor{\eta_{j}}\rangle$.
This is only remarkable for technical purposes.
\end{re}

To prove such property, an easy computation says that $\tor{\mu_4}=\span{t_1,t_2}$ for
$t_1=t_{\imag,1,\imag,\imag,\imag,-1}$ and $t_2=t_{1,\imag,1,1,\imag,1}$.
If we take $\hor{\mu_4^2}=\{t_{1,1,\alpha,1,\alpha\beta,\alpha}\mid \alpha^2=\beta^2=1\}$,
then $\wp_{\mu_4^2}(t_{x,y,\frac\alpha{x},x,\frac{\alpha\beta}{xy},\alpha})=t_{1,1,\alpha, 1,\alpha\beta,\alpha}$
and hence $\wp_{\mu_4^2}(t_1)=t_{1,1,-1,1,-1,-1}$, $\wp_{\mu_4^2}(t_1^2)=\id$, $\wp_{\mu_4^2}(t_2)=t_{1, 1,1,1,-1,1}$ and $\wp_{\mu_4^2}(t_2^2)=\id$.
Thus, the projection by $\wp_{\mu_4^2}$ of any proper quasitorus of $\tor{\mu_4 }$ is not the whole $\hor{\mu_4^2}$, in particular of $\span{ \Upsilon_2^2,\Upsilon_3}$,  which is equivalent to the fact that $\span{\Upsilon_1^2,\Upsilon_2^2,\Upsilon_3}\subset\Int\eseis$ is  toral.

Here another example of the situation described in Subsection~\ref{sub_estructuraMADs} appears: $\span{\Upsilon_1^2,\Upsilon_2,\Upsilon_3}\subset\Int\eseis$ is a non-toral quasitorus isomorphic to $\mathbb{Z}_2\times\mathbb{Z}_4^2$
which does not contain an elementary non-toral $2$-group.

\subsection{MAD-groups without order two outer automorphisms}\label{subsec_soloelZ43}

\begin{lm}\label{re_dosextbuenasdeorden2queconmutan}
There are two order two commuting automorphisms in $\No(\T)$ whose projections on the extended Weyl group
are $ \eta_3$
 and $   \sigma_{11127}$, respectively.
\end{lm}

\textbf{Proof.}
Start with $\Q_{11}$ the MAD-group isomorphic to $\mathbb{Z}_2^4\times\mathbb{Z}_4$ obtained after combining
an automorphism of type $2D$ with  a copy of $\Xi_5$, the  MAD-group of $\aut\mathfrak{c}_4$ isomorphic to $\mathbb{Z}_2^3\times\mathbb{Z}_4$ described in Subsection~\ref{sub_losquevienendec4}. Of course this copy of $\Xi_5$ is a non-toral quasitorus not only of $\aut\mathfrak{c}_4$ but of $\aut\eseis$. Besides the subquasitorus $\span{\tiny{\begin{pmatrix}1&0\cr0&\imag\end{pmatrix}}\otimes I_2\otimes \theta_3  ,I_2\otimes \theta_3\otimes I_2,I_4\otimes \theta_1 }\cong\Z_4\times  \Z_2^2$ is toral. This implies that, after conjugating,  we can find two order two commuting automorphisms $f$ and $g$, being $f$ inner, in $\No(\T)$ such that $\tor{f}\cap\tor{g}\cong\Z_4\times  \Z_2^2$. As
  $\pi(g)$ cannot be the identity, this implies the existence  of $j=1,\dots,5$, $\sigma_i$ in the orbit of $\sigma_{96}$ and $s,s'\in\T$ such that
$\Q_{11}$ is conjugated to $\langle \widetilde{\eta}_js,\widetilde{\sigma}_{i}s',\tor{\sigma_{i}}\cap\tor{\eta_{j}}\rangle$.
With the help of a computer, we study the elements in the orbit of $\sigma_{96}$ (there are only $45$) which commute with each of the $\eta_j$'s and divide in orbits, getting the following possibilities:
$$
\begin{array}{l}
j=1,\,i=25470,\tor{\sigma_{i}}\cap\tor{\eta_{j}}\cong  \mathbb{Z}_2^4 ,\\
j=1 ,\,i=2416 ,\tor{\sigma_{i}}\cap\tor{\eta_{j}}\cong  \mathbb{Z}_2^2\times \F^* ,\\
j=2 ,\,i=11127 ,\tor{\sigma_{i}}\cap\tor{\eta_{j}}\cong \mathbb{Z}_2^3  ,\\
j=2 ,\,i=11104 ,\tor{\sigma_{i}}\cap\tor{\eta_{j}}\cong  \mathbb{Z}_2\times \F^* ,\\
j=3 ,\,i= 11127,\tor{\sigma_{i}}\cap\tor{\eta_{j}}\cong  \mathbb{Z}_2^2\times\mathbb{Z}_4 ,\\
j=3 ,\,i= 11104,\tor{\sigma_{i}}\cap\tor{\eta_{j}}\cong  \mathbb{Z}_2^2\times \F^*,  \\
j=4 ,\,i=11007 ,\tor{\sigma_{i}}\cap\tor{\eta_{j}}\cong  \mathbb{Z}_2^2, \\
j= 5,\text{ any } i  ,\tor{\sigma_{i}}\cap\tor{\eta_{j}}\subset \tor{\eta_{5}} \cong\mathbb{Z}_2^6 .
\end{array}
$$
\hfill  $\square$ \smallskip

\begin{pr}\label{pr_sinohayextorden2eselQ14}
If $Q$ is a MAD-group of $\aut\eseis$, not contained in $\Int\eseis$, such that there is not an order two outer automorphism in $Q$,
then $Q$ is conjugated to $\Q_{14}$.
\end{pr}

\textbf{Proof.}
Take $Q$ a MAD-group of $\aut\eseis$, not contained in $\Int\eseis$, such that there is not an order two outer automorphism in $Q$.
Take $\T$ a maximal torus of $\aut\eseis$ such that $Q$ is contained in the normalizer of such torus,
$\No(\T)$, and $Q\cap\T$ is maximal toral in $Q$.
Hence there are indices $i,i_1,\dots,i_l\in\{1,\dots, 51480\}$ and toral elements $s,s_1,\dots,s_l\in\T$ such that
$
Q\cap\T=\tor{\sigma\sigma_i}\cap \tor{\sigma_{i_1}}\cap\dots\cap \tor{\sigma_{i_l}}
$
and
$$
Q=\langle  \widetilde{\sigma}\widetilde{\sigma}_is\rangle\cdot\langle \widetilde{\sigma}_{i_1}s_1,\dots,\widetilde{\sigma}_{i_l}s_l\rangle\cdot
Q\cap\T,
$$
and such that the quasitorus generated by $Q\cap\T\,\cup \{\widetilde\sigma_{i_j}s_j\}$ is non-toral for all $j$
and 
$\widetilde{\sigma}_{i_j}s_j\notin\langle \widetilde{\sigma}_{i_1}s_1,\dots,\widetilde{\sigma}_{i_{j-1}}s_{j-1}\rangle\cdot
Q\cap\T$.

We can assume that  $\widetilde{\sigma}\widetilde{\sigma}_is$ is an outer automorphism with order minimum  in $Q$.
This  order is $2^hm$ for $m$ an odd number, but note that $m=1$, because otherwise $(\widetilde{\sigma}\widetilde{\sigma}_is)^m$
would be outer with order $2^h$. Besides $h>1$, by hypothesis.

Take $r$ the order of $\sigma\sigma_i$. It divides $2^h$, so that $r\in\{2,4,8\}$ according to Table~(\ref{tablaexternos}).

$\star$
If $r=8$, we can assume that $i=17$ and so $Q\cap\T\subset\tor{\sigma\sigma_{17}}\cong\F^*$. By Lemma~\ref{le_toroporalgotoral},
the quasitorus $\langle(\widetilde{\sigma}\widetilde{\sigma}_{17}s)^2, Q\cap\T\rangle$
 is toral ($(\widetilde{\sigma}\widetilde{\sigma}_{17}s)^2 $ is an inner automorphism), and, according to our choice of $\T$, we have
$(\widetilde{\sigma}\widetilde{\sigma}_{17}s)^2\in Q\cap\T$. Hence $(\sigma\sigma_{17})^2=\id$, which is a contradiction.

$\star$
Suppose now that $r=4$. Hence we can assume that $\sigma\sigma_i\in\{\mu_1,\mu_2,\mu_3,\mu_4\}$. If $\sigma\sigma_i$ were $\mu_1$ or
$\mu_2$, we would obtain   a contradiction as in case $r=8$, since $\tor{\mu_1}\cong(\F^*)^2$ and  $\tor{\mu_2}\cong(\F^*)^3$ so that
$\span{(\widetilde{\sigma}\widetilde{\sigma}_{i}s)^2, Q\cap\T}$ would be in the conditions of Lemma~\ref{le_toroporalgotoral}. 
If $\sigma\sigma_i=\mu_3$, as $\tor{\mu_3^2}=\{t_{x,y,\frac1x,x,v,w}\mid x,y,v,w\ne0\}\cong(\F^*)^4$, then
$\langle(\widetilde{\sigma}\widetilde{\sigma}_{i}s)^2, Q\cap\T\rangle$ is contained in the toral quasitorus
$\langle(\widetilde{\sigma}\widetilde{\sigma}_{i}s)^2, \tor{\mu_3^2}\rangle$.
Hence we can assume that $\sigma\sigma_i=\mu_4$.
The case $l\ne0 $ never happens: as
$ \langle\widetilde{\sigma}_{i_1}s_1\rangle\cdot
Q\cap\T\subset\langle \widetilde{\sigma}_{i_1}s_1\rangle\cdot \tor{\mu_4}$ is non-toral, and $\tor{\mu_4}\cong\Z_4^2$,
then $\sigma_{i_1}$ is in the orbit of either $\sigma_{96}$ or $\sigma_{75}$. The first possibility
does not occur because the only element in the orbit of $\sigma_{96}$ which commutes with $\mu_4$ is $\mu_4^2$.
The second one is also impossible because $\sigma_{75}^2$ is not conjugated to
$\sigma_{96}$ (then  $ \langle(\widetilde{\sigma}_{i_1}s_1)^2\rangle\cdot
Q\cap\T$ would be toral and we would get a contradiction by the same arguments of the case $r=8$).
Thus $l=0$,
  $Q=Q(\widetilde\mu_4s)$, and,
 as $\tor{\mu_4}\cong\Z_4^2$ is finite, then $Q$ is conjugated to $Q(\widetilde\mu_4 )$ by Remark~\ref{re_paramover}.

Our purpose now is to check that this quasitorus $Q(\widetilde\mu_4)$ of type $\Z_4^3$ is the only possible MAD-group satisfying the required conditions. By Remark~\ref{re_esemododeencontrarZ4cubonoessorpresa}, we have to expect for its apparition in the case $r=2$.

$\star$
Thus suppose that $r=2$ and that (perhaps by changing 
the  element $s$ in the torus) an outer automorphism in $Q$ of minimum order $2^h$ is $\widetilde\eta_{j}s$ for some $j\in\{1,\dots,5\}$.
Observe first some useful facts:
\begin{itemize}
\item[a)] $l\ne0$.

Otherwise $Q=Q(\widetilde\eta_j s)=\langle\widetilde\eta_js,\tor{\eta_j}\rangle$, but this quasitorus contains outer automorphisms of order just $2$, 
a contradiction:
    Indeed, the automorphism $(\widetilde\eta_js)^2$ belongs to $\sor{\eta_j}$ by Remark~\ref{re_extensionesdeoreden2}, and it
    has order $2^{h-1}$ (a multiple of $2$), so that $\sor{\eta_j}\ne\id$ is a non-trivial torus. A torus contains a square root of each of its elements, so there is $s'\in\sor{\eta_j}\subset\tor{\eta_j}\subset Q$ such that
$(s')^2=(\widetilde\eta_j s)^2$, so that $\widetilde\eta_js(s')^{-1}$ is an outer order two automorphism in $Q$.

\item[b)] $(\widetilde\eta_js)^2\in\sor{\eta_j}$ but it is not contained in any subtorus of $Q\cap \T$. In particular, $j\ne5$ since
$\sor{\eta_5}=\id$ ($\tor{\eta_5}=\Z_2^6$).

\item[c)] There is no $k\in\{1,\dots,l\}$ such that $\sigma_{i_k}$ has order three.

In order to check it, take into account two facts.
First, the quasitorus $\langle \widetilde{\sigma}_{i_k}s_k,\tor{\eta_j}\cap\tor{\sigma_{i_k}}\rangle$ is non-toral since it   contains
$\langle \widetilde{\sigma}_{i_k}s_k,Q\cap\T\rangle$. This implies, by Lemma~\ref{le_toroporalgotoral}, that $\tor{\eta_j}\cap\tor{\sigma_{i_k}}$ is the direct product of a torus with
a finite group, say $P$. This $P$ contains a non-trivial $3$-group, since $\span{\widetilde{\sigma}_{i_k}s_k,P}$ is non-toral but $\span{(\widetilde{\sigma}_{i_k}s_k)^3,P}\subset\T$ is toral. Second, $(\widetilde{\eta}_js)^{2^{h-1}}$ is an order two element in
$\sor{\eta_j}\cap\tor{\sigma_{i_k}}$. Now we check that these conditions  do not happen for any value of $j=1,2,3,4$.

For $j=1$, there are $80$ order three elements in
$\mathcal{W}$ commuting with $\eta_1$, but they can divided in three orbits under  conjugation by some element which preserves $\eta_1$, with representatives $2920$, $12406$ and $3826$. We can look only at these representatives because if $p\in\mathcal{W}$ such that $p\eta p^{-1}=\eta$ and $p\nu p^{-1}=\nu'$, then $\tor{\eta}\cap\tor{\nu}\cong\tor{\eta}\cap\tor{\nu'}$.
But $\tor{\eta_1}\cap\tor{\sigma_{2920}}\cong(\F^*)^2$ and $\tor{\eta_1}\cap\tor{\sigma_{12406}}\cong(\F^*)^2$ have not  direct factors $\Z_3$, and $\tor{\eta_1}\cap\tor{\sigma_{3826}}\cong\Z_3^2$ has not order two elements.

For $j=2$, there are $8$ order three  elements (one is $3026$) in
$\mathcal{W}$ commuting with $\eta_2$, all of them conjugated preserving fixed $\eta_2$. Note that $\tor{\eta_2}\cap\tor{\sigma_{3026}}=\{t_{x,x^{-1},x^{-3},x^3,x^{-1},1}\mid x\ne0\}\cong\F^*$.

Again for $j=3$, there are $8$ order three  elements (one is $4796$) in
$\mathcal{W}$ commuting with $\eta_3$, all of them conjugated preserving fixed $\eta_3$. We see that $\tor{\eta_3}\cap\tor{\sigma_{4796}}=\{t_{1,y,z,1,y,\frac{1}{y^2z}}\mid y,z\ne0\}\cong(\F^*)^2$.

Finally, for $j=4$, there are $80$ order three  elements in
$\mathcal{W}$ commuting with $\eta_4$, but they can grouped in two conjugation orbits fixing $\eta_4$, with representatives $3839$, $4079$.
Now we check that
$\tor{\eta_4}\cap\tor{\sigma_{ 3839}}=\{t_{1,1,z,1,1,\frac1z }\mid z\ne0 \}\cong\F^*$ and
$\tor{\eta_4}\cap\tor{\sigma_{ 4079}}=\{t_{1,y,z,u,y,\frac yz}\mid y^2=u^2=1,z\ne0  \}\cong\F^*\times\Z_2^2$.

\item[d)] The matrix $\sigma_{i_k}$ has order just $2$ for all $k\in\{1,\dots,l\}$ and it belongs to the orbit of $\sigma_{96}$.

Indeed, if $\sigma_{i_k}$ has order  $5$ or $10$, then $\langle \widetilde{\sigma}_{i_k}s_k,\tor{\eta_k} \rangle$ is toral,
and, if the order $\sigma_{i_k}$ is multiple of $3$,
then one of its powers $\sigma_{i_k}^m$ has order just
three, and the arguments in   item c) work for  $\sigma_{i_k}^m$. Thus $\sigma_{i_k}$ has order  a power of $2$.
As in   Remark~\ref{re_losindicesnotorales},  the element $\sigma_{i_k}$ is conjugated to either $\sigma_{96}$ or $\sigma_{75}$. But in the
latter case, we reason again that $\langle (\widetilde{\sigma}_{i_k}s_k)^2,Q\cap\To\rangle$ is toral since $\tor{(\sigma_{75})^2}\cong(\F^*)^4$, which is a contradiction with the maximal-torality of $Q\cap\T$.
\end{itemize}

With items a), b) and d) in mind, we proceed to a detailed analysis of the possible cases.

$\clubsuit$ Case $j=4$. There are $15$ elements in the orbit of $\sigma_{96}$ which commute with $\eta_4$, but all of them are conjugated
by means of an element of $\mathcal{W}$ which fixes $\eta_4$. So we can assume that $\sigma_{i_1}$ is any of them, for instance,
$\sigma_{11007}$.

As $\langle\widetilde{\sigma}_{i_1}s_1,\tor{\sigma_{i_1}}\cap\tor{\eta_{4}}\rangle$ is non-toral (it contains
$\langle\widetilde{\sigma}_{i_1}s_1,Q\cap\T\rangle$), then by Lemma~\ref{le_toroporalgotoral}
$\langle\widetilde{\sigma}_{i_1}s_1,\wp_{i_1}(\tor{\sigma_{i_1}}\cap\tor{\eta_{4}})\rangle$ is non-toral too
  and
     $\wp_{i_1}(\tor{\sigma_{i_1}}\cap\tor{\eta_{4}})=\hor{i_1}$ as in Remark~\ref{re_proyecciones}.
Note that
$$\begin{array}{l}
\tor{11007}=\{t_{x,x,\frac uz,z,x,w}\mid  u^2=x^2=1,z,w\ne0\}\cong(\F^*)^2 \times\Z_2^2,\\
\sor{11007}=\{t_{1,1,z,\frac1z,1,w}\mid z,w\ne0\}\cong(\F^*)^2,
\end{array}
$$
so that a complement satisfying the conditions in Lemma~\ref{le_comoeselTsuperf} is, for instance,
$$
\hor{11007}=\{t_{x,x,u,1,x,1}\mid u^2=x^2=1\}\cong\Z_2^2.
$$
Now $\tor{\eta_4}\cap\tor{\sigma_{ 11007}}=\{t_{1,1,z,u,1,w }\mid z^2=u^2=w^2=1  \}\cong\Z_2^3$,
so we have to compute
 $\wp_{11007}(t_{1,1,z,u,1,w })$. As $t_{1,1,z,u,1,w }=t_{1,1,\frac1u,u,1,w}t_{1,1,uz,1,1,1}$, with the first
 factor in $\sor{\sigma_{ 11007}}$ and the second one in $\hor{\sigma_{ 11007}}$,
 then $\wp_{11007}(t_{1,1,z,u,1,w })=t_{1,1,uz,1,1,1}$.
So that we obtain a contradiction since  $\wp_{11007}(\tor{\sigma_{11007}}\cap\tor{\eta_{4}})\cong\Z_2$  (roughly speaking, although
$\tor{\sigma_{11007}}\cap\tor{\eta_{4}}$ contains three $\Z_2$'s, it only contains one of the two \emph{bads} required $\Z_2$'s).

$\clubsuit$ Case $j=2$. There are $7$ elements in the orbit of $\sigma_{96}$ which commute with $\eta_2$, divided in two orbits of
 $\mathcal{W}$ under conjugation by an element fixing $\eta_2$, with representatives, for instance, $11127$ and $11104$.
If $\sigma_{i_1}= \sigma_{11127}$,  then $\tor{\eta_2\sigma_{11127}}\cong\F^*\times\Z_2^4$, so that $\eta_2\sigma_{11127}$ is conjugated to $\eta_4$
and such case has already been studied.
If $\sigma_{i_1}= \sigma_{11104}$,  then $\tor{\sigma_{11104}}\cap\tor{\eta_{2}}=\{t_{x,y,x,x,y,\frac1{y^2}}\mid x^2=1,y\ne0\}\cong\F^*\times\Z_2$, so that  $\langle\widetilde{\sigma}_{11104},\tor{\sigma_{11104}}\cap\tor{\eta_{2}}\rangle$ is toral by Remark~\ref{re_cuandoQ(f)toralseguro}, and again we have found a contradiction.

$\clubsuit$ Case $j=1$. There are $13$ elements in the orbit of $\sigma_{96}$ which commute with $\eta_1$, divided in two orbits of
 $\mathcal{W}$ under conjugation by an element fixing $\eta_1$, with representatives, for instance, $25470$ and $2416$.
 A simple computation shows us that
$\tor{\eta_1\sigma_{25470}}\cong \Z_2^6$ and that
$\tor{\eta_1\sigma_{2416}}=\{t_{x,y,z,u,\frac{y}x,\frac{y}z}\mid y^2=u^2=1,x,z\ne0\}\cong(\F^*)^2\times \Z_2^2$,
hence $\eta_1\sigma_{25470}$ is conjugated to $\eta_5$ and $\eta_1\sigma_{2416}$ is conjugated to $\eta_4$, and both cases have been previously considered.

$\clubsuit$ Case $j=3$. There are $5$ elements in the orbit of $\sigma_{96}$ which commute with $\eta_3$, divided in two orbits of
 $\mathcal{W}$ under conjugation by an element fixing $\eta_3$, $O_1=\{k_1=10850,k_2=11104\}$ and $O_2=\{k_3=23234,k_4=11127,k_5=28154\}$.
 Thus we can take $i_1\in\{k_2,k_4  \}$.
 It is important to note that $l=1$, because
 there are no $k_i ,k_j,k_m\in O_1\cup O_2$ such that $\sigma_{k_i}\sigma_{k_j}=\sigma_{k_m}$, that is, $\sigma_{k_i}\sigma_{k_j}$ is not in the orbit
 of $\sigma_{96}$.

In the first subcase,
$(\widetilde{\eta}_3s)^2\in\sor{\eta_3}\cap\tor{10850}=\{t_{1,y,\frac1{y^2},1,y,1}\mid y\ne0\}\cong\F^*$, a contradiction with the item b), since $\sor{\eta_3}\cap\tor{10850}\subset Q\cap\T$ ($l=1$).
Consider then  that $i_1=11127$.
For our convenience, we take the extensions $\widetilde{\eta}_3$ and $\widetilde\sigma_{11127} $ such that not only
$\widetilde{\eta}_3^2=\widetilde\sigma_{11127}^2=\id$ but besides $\widetilde{\eta}_3\widetilde\sigma_{11127} =\widetilde\sigma_{11127} \widetilde\eta_3$. Note that we can make such choice by Lemma~\ref{re_dosextbuenasdeorden2queconmutan}.
As $l=1$,  there are $s,s'\in\T$    such that
$Q=\langle \widetilde{\eta}_3s,\widetilde{\sigma}_{11127}s',\tor{\sigma_{i_1}}\cap\tor{\eta_{3}}\rangle$.
Take $f_1=\widetilde{\eta}_3s$, $f_2=\widetilde{\sigma}_{11127}s'$ and $f_3=t_{-1,\imag,1,-1,-\imag,1}\in \tor{\sigma_{11127}}\cap\tor{\eta_{3}}=\{t_{y^2,y,z,y^2,y^3,w}\mid y^4=z^2=w^2=1\}\cong\Z_4\times\Z_2^2$.

Our first aim is to check first that $f_1$ and $f_2$ are order 4 automorphisms (with $f_1^2\ne f_2^2$), what implies  that $Q=\span{f_1,f_2,f_3}$.
Note that if we could prove that:
\begin{itemize}
\item[i)] $f_1$ is an outer automorphism,
\item[ii)] $\span{f_2,f_3}$ is toral,
\item[iii)] $\span{f_1^2,f_2,f_3}\subset\Int\eseis$ is non-toral,
\end{itemize}
hence,
we could apply the same arguments than in Proposition~\ref{pr_elZ4cuboenterminoscomputacionales}
to conclude that $Q$ is conjugated to $Q(\widetilde\mu_4 )$ and then it would be also conjugated  to $\Q_{14}$.

Thus, our second aim is to prove items ii) and iii).
In order to make such comprobations about torality, take $\hor{11127}=\{t_{x,1,z,x,1,1}\mid x^2=z^2=1\}$.
The corresponding projection is hence given by
$$
\begin{array}{l}
\wp_{11127}\colon\tor{11127}\to\hor{11127}\\
\wp_{11127}(t_{x,y,z,x,\frac1y,w} )= t_{x,1,z,x,1,1},
\end{array}
$$
 if $x^2=z^2=1$, $y,w\ne0$.
We will check that $\wp_{11127}(f_2^2)=\id$, which guarantees that ii) is verified,  since
$\wp_{11127}(\span{f_2^2,f_3})=\span{\wp_{11127}( f_3)}=\span{t_{-1,1,1,-1,1,1}}\ne\hor{11127}\cong\Z_2^2$.
And finally we will check that $ \wp_{11127}(f_1^2)\notin\span{ t_{-1,1,1,-1,1,1}}$, so that
$\wp_{11127}(f_1^2,f_3)=\hor{11127}$ and $\span{f_1^2,f_2,f_3}\subset\Int\eseis$ is non-toral.

 Begin by noting that
 $$
 \begin{array}{l}
 (\widetilde{\eta}_3s)^2\in \sor{\eta_3}\cap\tor{11127}=\{t_{1,y,z,1,y,z}\mid y^2=z^2=1\}\cong\Z_2^2,\\
 (\widetilde{\sigma}_{11127}s')^2\in \tor{\eta_3}\cap\sor{11127}=\{t_{1,y,1,1,y,w}\mid y^2=w^2=1\}\cong\Z_2^2,
 \end{array}
 $$
 so that $f_1$ has order $4$ (not $2$ by hypothesis)
 and $ f_2$ has order either $2$ or $4$.
 Besides $f_1^2\ne t_{1,-1,1,1,-1,1}=f_3^2$, because otherwise
    $f_1/f_3$ would be an outer order two automorphism in $Q$. If $f_1^2=t_{1,-1,-1,1,-1,-1}$,
then $(f_1f_3)^2=t_{1, 1,-1,1, 1,-1}$, so that we can assume that $f_1^2=t_{1,1,-1,1,1,-1}$.
In particular,
  $\wp_{11127}(f_1^2=t_{1,1,-1,1,1,-1})= t_{1,1,-1,1,1,1}\notin\span{ t_{-1,1,1,-1,1,1}} $.

Besides this implies that there are $a,b,c,d\in\F^*$ such that $s=t_{a,b,c,\frac{-1}{ac^2},\frac{-ac^2}{b},d}$.
Now take   the element $r_1=t_{u,\frac{b}{uc},1,\frac1{uc},1,w}$ for some $u,w\in\F^*$ such that $u^2=-a$ and
$w^2=\frac{-adc^3}{b^2}$. It is a straightforward computation that $r_1(\widetilde{\eta}_3s)r_1^{-1}=\widetilde{\eta}_3t_{-1,1,1,1,1,1}$,
so that we can assume that $s=t_{-1,1,1,1,1,1}$.
Now the commutativity of $\widetilde{\eta}_3s$ with $\widetilde{\sigma}_{11127}s'$
means that $s'=t_{x,y,z,x,xy,\frac{\a \imag}{y^2z}}$ for some $x,\a\in\{\pm1\}$ and $y,z\in \F^*$.
Then $f_2^2=t_{1,\frac1x,1,1,x,-x^6\a^2}=t_{1,\frac1x,1,1,x,-1}\ne\id$, so that $f_2$ has order $4$ and in any case $f_2^2\ne f_1^2$.
Besides the projection $\wp_{11127}(t_{1,\frac1x,1,1,x,-1})=\id$, what finishes the proof.
\hfill  $\square$ \smallskip


In fact the power of these kind of arguments is strong, and we could have described all the MAD-groups of $\aut\eseis$ in computational terms:
  It is possible to take good choices of the extensions $\widetilde\sigma_{25470},\widetilde\sigma_{3826},\widetilde\sigma_{11104},\widetilde\sigma_{11127}$
  (the two latter ones as in Remark~\ref{re_dosextbuenasdeorden2queconmutan} and in Lemma~\ref{re_dosextbuenasdeorden2queconmutan}, respectively) such that
\begin{equation}\label{eq_losMADsversioncomputacional}
\begin{array}{ll}
\Q_6\cong Q(\widetilde\eta_1),\qquad&
\Q_{7}\cong \span{\widetilde\eta_1,\widetilde\sigma_{25470} ,\tor{\eta_{1}}\cap\tor{\sigma_{25470}} } ,\\
\Q_{10}\cong Q(\widetilde\eta_3),&
\Q_{9}\cong \span{\widetilde\eta_1,\widetilde\sigma_{3826} ,\tor{\eta_{1}}\cap\tor{\sigma_{3826}} },\\
\Q_{12}\cong Q(\widetilde\eta_4),&
\Q_{8}\cong \span{\widetilde\eta_3,\widetilde\sigma_{10850} ,\tor{\eta_{3}}\cap\tor{\sigma_{11104}} },\\
\Q_{13}\cong Q(\widetilde\eta_5),&
\Q_{11}\cong \span{\widetilde\eta_3,\widetilde\sigma_{11127} ,\tor{\eta_{3}}\cap\tor{\sigma_{11127}} }.
\end{array}
\end{equation}
The proof can be made with analogous arguments to those ones in the proof of Proposition~\ref{pr_sinohayextorden2eselQ14}.
In fact, with those reasonings
we could   have proved that $\{\Q_i\mid i=1,\dots,14\}$ are all the MAD-groups up to conjugation without using the
 description of the non-toral elementary $p$-groups. The paper would have then been
practically self-contained, but, for evident reasons (to avoid most of the such unpleasant computational arguments) we have preferred
a combined option.

\section*{Appendix}

In this section we  will provide natural descriptions of $\widetilde\sigma$  for all the representatives of the orbits of order two elements in $\mathcal{V}$. This may help the reader to have a better understanding of the situation.
All the computations here are made by hand.

But this section has also a practical objective: to be able to state Remark~\ref{re_extensionesdeoreden2} (a key piece in the proof of Proposition~\ref{pr_sinohayextorden2eselQ14}) thanks to the existence of outer order two automorphisms which project on $\eta_j$ for all $j=1,\dots,5$.

We do not want to construct explicit expressions of  these extensions by following the lines in \cite[Section~14.2]{Humphreysalg}.
We prefer the following  procedure, also constructive: Choose a maximal torus, take some distinguished automorphisms and compute their projections on
$\No(T)/T$.
Take into account that the only computation of $\tor{f}$ is enough to distinguish the orbit of the element (among the elements in $\mathcal{W}$ and also in
$\si\mathcal{W}$).

We work again with the model described in Equation~(\ref{eq_modeloAdams}). We choose $\{E_1,E_2,E_3\}$ a basis of $V$ and $\{e_1,e_2,e_3\}$ its dual basis
(of $V^*$). We call $e_{ij}^k$ the element in $\sll(V_k)$ which maps  $E_j\in V_k$ into $E_i\in V_k$ (and $V_l$ into $0$ if $k\ne l$).
Take $\frak{h}$ the abelian subalgebra spanned by $ \{h_1=e_{11}^1-e_{33}^1,h_2=e_{22}^1-e_{33}^1,h_3=e_{11}^2-e_{33}^2,
h_4=e_{22}^2-e_{33}^2,h_5=e_{11}^3-e_{33}^3,h_6=e_{22}^3-e_{33}^3\} $, which is a Cartan subalgebra of $\mathcal{L}$.
Now, if $h=\sum_{i=1}^6 w_ih_i$ is an arbitrary element in $\frak{h}$, then
$[h,E_i\otimes E_j\otimes E_k]=(\alpha_{1,i}+\alpha_{2,j}+\alpha_{3,k}) E_i\otimes E_j\otimes E_k$ for all $i,j,k\in\{1,2,3\}$,
for the scalars $\alpha_{l,1}=w_{2l-1}$, $\alpha_{l,2}=w_{2l}$ and  $\alpha_{l,3}=-w_{2l-1}-w_{2l}$, $l\in\{1,2,3\}$. Define
$\alpha_i\colon\frak{h}\to\F$
by
$$
\begin{array}{ll}
\alpha_1(h)=w_1-w_2,\qquad &  \alpha_4(h)=-w_1-w_2-w_3-w_4-w_5-w_6 ,\\
\alpha_2(h)=2w_3+w_4,   & \alpha_5(h)=w_5+2w_6, \\
\alpha_3(h)=w_1+2w_2,& \alpha_6(h)=w_5-w_6.
\end{array}
$$
It is a simple computation that $\{\alpha_i\}_{i=1}^6$ is a set of simple roots of $\Phi$, the root system relative to $\frak{h}$.
Moreover, the root spaces corresponding to the simple roots are
$$
\begin{array}{lll}
L_{\alpha_1 }=\langle e_{12}^1 \rangle, &  L_{\alpha_3 }=\langle e_{ 23}^1 \rangle,  & L_{\alpha_5 }=\langle e_{23 }^3 \rangle, \\
L_{\alpha_2 }=\langle e_{13 }^2 \rangle,&L_{\alpha_4 }=\langle E_3\otimes E_3\otimes E_3 \rangle,   & L_{\alpha_6}=\langle e_{ 12}^3 \rangle.
\end{array}
$$
If $\mathfrak{T}$ is the torus of the automorphisms fixing pointwise $\frak{h}$, then an automorphism $f\in\aut\eseis$ belongs to $\No(\mathfrak{T})$
if and only if $f$ permutes the root spaces. Precisely, if there is $\eta\in\GL(E)$ for $E=\sum_{i=1}^6\F\a_i$ such that $f(L_\a)\subset L_{\eta(\a)}$
for all $\a\in\Phi$, then $\pi(f)=\eta\in\mathcal{V}$. Identify $\GL(E)$ with $\GL(6)$ by mapping each $\eta$ to its matrix relative to the basis $\{\a_i\}_{i=1}^6$.

If   $f,g,h\in\GL(3)\equiv\GL(V)$,    we call $f\otimes g\otimes h$ the unique automorphism of $\eseis$ whose action in $\mathcal{L}_{\bar1}$ is $u\otimes v\otimes w\mapsto f(u)\otimes g(v)\otimes h(w)$  (thus $\psi(A)$ in Subsection~\ref{sub Z34gr} coincides with $A\otimes A\otimes A$).

Take $\varphi_1= \tiny{\begin{pmatrix}1&0&0\cr 0&0&1 \cr 0&1&0
\end{pmatrix}}\otimes I_3\otimes I_3\in\aut\eseis$, which obviously is an order two automorphism.
As $\varphi_1(e_{12}^1)=e_{13}^1\in L_{\alpha_1 +\alpha_3}$,  $\varphi_1(e_{23}^1)=e_{32}^1\in L_{-\alpha_3}$ and $\varphi_1(E_3\otimes E_3\otimes E_3)=E_2\otimes E_3\otimes E_3\in L_{\alpha_3 +\alpha_4}$, that means that $\pi(\varphi_1)=\varrho_1$ for $\varrho_1$ the automorphism of the root
system mapping $\alpha_1$ into $\alpha_1+\alpha_3$, $\alpha_3$ into $-\alpha_3$, $\alpha_4$ into $\alpha_3+\alpha_4$ and $\alpha_2$,
$\alpha_5$ and $\alpha_6$ into themselves, whose related matrix is:
$$
\varrho_1=\tiny{\begin{pmatrix}1&0&1&0&0&0\cr 0&1&0&0&0&0 \cr 0&0&-1&0&0&0 \cr
0&0&1&1&0&0 \cr 0&0&0&0&1&0 \cr 0&0&0&0&0&1
\end{pmatrix}}.
$$
Hence $\tors{\varrho_1}=\{t\in\mathfrak{T}\mid \varrho_1\cdot t=t\}=\{t_{x,y,1,u,v,w}\mid x,y,u,v,w\in\F^*\}\cong(\F^*)^5$ and
$\pi( \varphi_1)$
is conjugated to $\sigma_{11323}$. In particular $\varphi_1\in\Int\eseis$.

Take $\varphi_2= I_3\otimes \tiny{\begin{pmatrix}1&0&0\cr 0&0&1 \cr 0&1&0
\end{pmatrix}}\otimes I_3\in\aut\eseis$. By symmetry with the previous case, $\varphi_2$ is also an inner automorphism.
It satisfies
$$
\begin{array}{l}
\varphi_2(e_{13}^2)=e_{12}^2\in L_{\alpha_1 +2\alpha_2+2\alpha_3+3\alpha_4+2\alpha_5+\alpha_6},\\
\varphi_2(E_3\otimes E_3\otimes E_3)=E_3\otimes E_2\otimes E_3\in L_{-\alpha_1-\alpha_2-2\alpha_3-2\alpha_4-2\alpha_5-\alpha_6},
\end{array}
$$
and it does not move the other simple root spaces, so that
 $\varrho_2=\pi(\varphi_2)$ is the element in the Weyl group
$$
\varrho_2=\tiny{\begin{pmatrix}1&0&0&0&0&0\cr 1&2&2&3&2&1 \cr 0&0&1&0&0&0 \cr
-1&-1&-2&-2&-2&-1 \cr 0&0&0&0&1&0 \cr 0&0&0&0&0&1
\end{pmatrix}}.
$$
It is a hand-computation that $\tors{\varrho_2}$ is again isomorphic to $ (\F^*)^5$,  but that $\tors{\varrho_1\varrho_2}$ is   isomorphic to $ (\F^*)^4$, hence implying that $\varphi_1\varphi_2$ is an order two automorphism  ($\varphi_1$ commutes with $\varphi_2$)
such that $\pi(\varphi_1\varphi_2)$ is conjugated to ${\sigma}_{19}$.

Third take $\varphi_3=I_3\otimes  I_3\otimes \tiny{\begin{pmatrix}1&0&0\cr 0&0&1 \cr 0&1&0
\end{pmatrix}}\in\Int\eseis$ and check that it only moves the simple root vectors
 $\varphi_3(e_{12}^3)=e_{13}^3\in L_{\alpha_5 +\alpha_6}$,  $\varphi_3(e_{23}^3)=e_{32}^3\in L_{-\alpha_5}$ and $\varphi_3(E_3\otimes E_3\otimes E_3)=E_3\otimes E_3\otimes E_2\in L_{\alpha_4 +\alpha_5}$, so that
 $$
 \pi(\varphi_3)=\varrho_3=\tiny{\begin{pmatrix}1&0&0&0&0&0\cr 0&1&0&0&0&0 \cr 0&0&1&0&0&0 \cr
0&0&0&1&1&0 \cr 0&0&0&0&-1&0 \cr 0&0&0&0&1&1
\end{pmatrix}}.
 $$
Hence $t_{x,y,z,u,v,w}\in\tors{\varrho_1\varrho_2\varrho_3}$ whenever $z=v=1=xyu^3w$, that is, $\tors{\varrho_1\varrho_2\varrho_3}\cong(\F^*)^3$
and  $\pi(\varphi_1\varphi_2\varphi_3)$ is conjugated to $ \sigma_{21}$. Note that there is no confusion with the orbit of $\eta_2$, because
all the $\varphi_i$'s until this moment are inner, so that $\varphi_1\varphi_2\varphi_3 $ is also.

Take $\varphi_4$ the order two automorphism of $\eseis$ mapping $E_i\otimes E_j\otimes E_k$ into $e_i\otimes e_j\otimes e_k\in (V^*)^{\otimes3}$.
Note that it maps any element in $ \sll(V_i)$  into   the opposite of its transpose. Hence $\varrho_4=\pi(\varphi_4)=-\id=\eta_5$ because it applies
each $\alpha_i$ into $-\alpha_i$.
 Note that $\varphi_4$ is an outer automorphism, more precisely, $\varphi_4\in 2D$ since $\dim\fix\varphi_4=\frac{72}2$.

Another order two outer automorphism, which will be denoted by $\varphi_5$, is that one interchanging $V_1$ with $V_3$ (which fixes a subalgebra isomorphic to $\frak{f}_4$,
as we saw in Subsection~\ref{sub_losquevienendef4}). As it applies $e_{ij}^1$ into $e_{ij}^3$, its
projection $\varrho_5=\pi(\varphi_5)$ interchanges $\alpha_1$ with $\alpha_6$ and  $\alpha_3$ with $\alpha_5$, not moving $\alpha_2$ and
$\alpha_4$. Hence $\varphi_5\in 2C$ is an order two outer automorphism with  $\pi(\varphi_5)=\sigma$.

Observe that $ \varphi_4\varphi_5$ is an inner order two automorphism ($\varphi_4$ and $\varphi_5$ commute) whose projection is in the orbit of $\si_{96}$.

The remaining descriptions of representatives of the outer orbits can be found by combining the previous elements. Notice that
$\varphi_5$ commutes with $\varphi_2$ (not with $\varphi_1$ or $\varphi_3$) and $\varphi_4$ commutes with all the $\varphi_i$'s.
Hence, by looking again at the $\tors{j}$'s we get that
$$
\begin{array}{l}
\pi(\varphi_5\varphi_2)\sim \eta_2,\\
\pi(\varphi_4\varphi_2)\sim \eta_4,\\
\pi(\varphi_4\varphi_1\varphi_2)\sim \eta_3.
\end{array}
$$


It is   interesting to observe that the MAD-groups $\Q_1$ and $\Q_2$ are just contained in the normalizer of this maximal torus $\mathfrak{T}=\{f\in\aut\eseis\mid f\vert_{\frak{h}}=\id\}$. More  precisely,
$F_1=t_{\omega^2,\omega,\omega^2,1,\omega^2,\omega^2}$, $\pi(F_2)=\sigma_{51529}$ (in the $292$-orbit), $F_3=t_{1,1,1,\omega,1,1}$
and $\pi(F_4)=\sigma_{30245}$ (in the $3819$-orbit) and $T_{\alpha,\beta}=t_{\frac{\alpha}{\beta},\alpha^2\beta,\alpha\beta^2,
\frac1{\alpha^3\beta^3},\alpha\beta^2,\frac\alpha\beta}$.

\begin{re}\label{re_dosextbuenasdeorden2queconmutan}
Note that again most of   the MAD-groups can be described in these natural terms starting from  Equation~\ref{eq_modeloAdams}.
If we denote by $Q(f;g)$ the quasitori generated by $f$, $g$ and $\tor{f}\cap\tor{g}$ for each pair of commuting automorphisms $f,g\in\No(\T)$, then
$$
\begin{array}{ll}
\Q_1=Q(F_4),      &  \Q_{8}\cong Q(\varphi_4\varphi_1\varphi_3;\varphi_4\varphi_5 ) , \\
\Q_2=Q(F_2) ,     &   \Q_{9}\cong Q(F_4\varphi_5 ), \\
\Q_4\cong Q(\varphi_4\varphi_5),      &   \Q_{10}\cong Q(\varphi_4\varphi_1\varphi_3 ) ,\\
\Q_5\cong Q(\id)  ,    &  \Q_{12}\cong Q(\varphi_4\varphi_1 ) , \\
\Q_6\cong Q( \varphi_5)    ,  &  \Q_{13}\cong Q(\varphi_4 ) . \\
\Q_7\cong Q(\varphi_4;\varphi_5) ,     &
\end{array}
$$

Observe that, although $\varphi_4\varphi_1\varphi_3$ and
$\varphi_4\varphi_5$ are two order two commuting automorphisms such that $\pi(\varphi_4\varphi_1\varphi_3)\sim \eta_3$
and $\pi(\varphi_4\varphi_5) \sim \sigma_{96}\sim\sigma_{11127}$, this does not imply the existence of order two commuting extensions
$\widetilde\eta_3$ and $\widetilde\sigma_{11127}$, as was necessary  in Subsection~\ref{subsec_soloelZ43}, because the conjugating element has to be the same. In other words, what have been found are two   order two commuting extensions
$\widetilde\eta_3$ and $\widetilde\sigma_{11104}$, what is checked by  computing the stabilizer.

\end{re}

\section*{Thanks}
We would like to thank  Professors C\'{a}ndido Mart\'{\i}n Gonz\'{a}lez and Alberto Elduque Palomo for their 
continuous support. The stays of the first author in the University of Zaragoza have been very fruitful for this work.


\end{document}